%% file: main.arxiv.tex
\documentclass[a4paper,10pt]{article}
\usepackage[a4paper, margin=2.5cm]{geometry}

\usepackage[table]{xcolor}
\usepackage[strict]{changepage}
\usepackage[colorinlistoftodos]{todonotes}
\usepackage[T1]{fontenc}
\usepackage[utf8]{inputenc}
\usepackage{graphicx}
\usepackage{pdflscape}
\usepackage{fancyvrb}
\usepackage{amssymb}
\usepackage{amsmath}

\usepackage{minted}
\usepackage{listings}

\usepackage{caption}
\usepackage{multirow}
\usepackage{tabularx}
\usepackage{longtable}
\usepackage{booktabs}
\usepackage{calc} 
\usepackage{tikz}
\usetikzlibrary{positioning, arrows.meta, trees}
\usepackage[%
    autostyle=true,%
    babel,%
]{csquotes}
\usepackage{authblk}

\usepackage{url}
\usepackage[%
    naturalnames=false,%
    breaklinks=true,%
    colorlinks=true,
    frenchlinks=false,%
    linkcolor=blue,%
    citecolor=blue,%
    debug=true,%
    filecolor=magenta,%
    urlcolor=magenta,%
    linktoc=all,
]{hyperref}                 
\usepackage{cleveref}

\newcommand{\nt}[1]{\ensuremath{\langle\mathrm{#1}\rangle}}
\newcommand{\kw}[1]{\texttt{#1}}
\newcolumntype{Y}{>{\centering\arraybackslash}X}

\newcommand{\nd}{\noindent}
\newcommand{\I}{\ensuremath{\mathcal{I}}}
\newcommand{\highi}[1]{\ensuremath{\displaystyle {#1}^\I}}

\title{Fuzzy OWL~2 Reasoning: A Re-Engineered Python Framework}

\author[1,2]{Fernando Bobillo}
\author[3]{Giuseppe Filippone\thanks{Corresponding author}}
\author[3]{Gianmarco La Rosa}
\author[4]{Umberto Straccia}
\author[3,5]{Marco Elio Tabacchi}

\affil[1]{University of Zaragoza, Zaragoza, Spain}
\affil[2]{Aragon Institute of Engineering Research (I3A), Zaragoza, Spain}
\affil[3]{Dipartimento di Matematica e Informatica, Università degli Studi di Palermo, Palermo, Italy}
\affil[4]{CNR-ISTI, Pisa, Italy}
\affil[5]{Istituto Nazionale di Ricerche Demopolis, Palermo, Italy}

\date{}

\begin{document}

\maketitle


\begin{abstract}
In many real-world domains, knowledge is inherently vague or imprecise --- features that classical ontology languages based on crisp Description Logics (DLs) are unable to capture effectively. This shortcoming poses particular challenges for applications in the Semantic Web and Explainable Artificial Intelligence (XAI), where robust reasoning over vague or graded information is essential. Fuzzy ontologies address this limitation by enriching DLs with fuzzy logic, enabling the expression of partial truth and supporting more nuanced modelling of real-world knowledge. 

This article presents \texttt{fuzzy-dl-owl2}, a complete re-engineering in Python of the \texttt{fuzzyDL} reasoner and the \texttt{Fuzzy OWL~2} framework. 
The former is an expressive fuzzy DL reasoner, while the latter allows for defining fuzzy ontologies within OWL~2, a W3C standard to formalise ontologies. 

Our contribution addresses several shortcomings of the original software, including semantic inconsistencies, rigid architectural design, and limited solver integration. The re-implementation features a modular class hierarchy tailored for extensibility, supports a broader range of Mixed-Integer Linear Programming (MILP) solvers (including open-source alternatives), and corrects IRI 
ambiguities arising from overlapping ontological elements.
Furthermore, a dedicated Python library (\texttt{pyowl2} \cite{pyowl2}) has also been developed to handle OWL~2 annotations in a standards-compliant manner, improving interoperability with existing Semantic Web tooling and resolving IRI ambiguities. The resulting framework offers a portable, extensible, and theoretically grounded platform for reasoning with fuzzy ontologies, suitable for both research and deployment in vague-aware systems.  
Performance tests have also been conducted that show improved execution times w.r.t.~the original Java implementation.
The source code and full documentation are publicly available to facilitate community adoption and further development.

\vspace{1em}
\noindent \textbf{Keywords:} Fuzzy Ontologies, \texttt{Fuzzy OWL~2}, Fuzzy Description Logics, \texttt{fuzzyDL} Reasoner, Python Library

\end{abstract}

    \section{Introduction}\label{sec:intro}
    \nd In recent decades, there has been an increasing interest in ontology languages that are capable of handling vagueness, which is often a requirement in real-world applications. In fact, despite the considerable success of ontologies, classical, i.e., crisp ontology languages are inadequate in addressing vagueness, a trait inherent to most real-world domains~\cite{uncertainty}.

    The modelling of vague information is facilitated by fuzzy set theory and fuzzy logic~\cite{zadeh1965fuzzy}. Consequently, fuzzy ontologies have been effectively applied in diverse areas, including information retrieval~\cite{Calegari2008,Knappe2007,Simou2008} and the Semantic Web~\cite{CostaLaskeyLukasiewicz2008,sanchez2006fuzzy,Foga2006}. The aforementioned applications underscore the pragmatic utility of integrating fuzziness in ontology-based systems.
    
    \emph{Description Logics} (DLs)~\cite{Baader_Calvanese_McGuinness_Nardi_Patel-Schneider_2007} are a family of logics that have been tailored for the representation of structured knowledge. The aforementioned offer a Tarski-style declarative semantics and a well-defined set of constructors (indicated by sequences such as $\mathcal{ALC}$, $\mathcal{SROIQ}$, etc.). This enables precise control over expressivity and computational complexity. DLs have underpinned ontology languages for decades, most notably the W3C standard for expressing ontologies OWL~2~\cite{Grau2008,OWL2}, which is grounded in the expressive DL $\mathcal{SROIQ}(\mathbf{D})$.
    
    Researchers have proposed the use of \emph{fuzzy DLs} to model vagueness explicitly, in recognition of the limitations of crisp DLs. Since the seminal introduction of fuzzy DLs by J. Yen in 1991~\cite{Yen1991}, there has been extensive theoretical work which has expanded the landscape of fuzzy logics and DLs (see~\cite{Borgwardt2017,CrossChen2018,LUKASIEWICZ2008291} for reviews on this topic). These formalisms extend classical DLs with fuzzy concept constructors and truth degrees.
    
    A variety of reasoning engines nowadays support fuzzy DLs, including \texttt{fuzzyDL} \cite{BobilloStraccia2016}, DeLorean~\cite{BOBILLO2012258}, and FiRE~\cite{Stoilos2006}. In particular, \texttt{fuzzyDL}, which is widely recognised as the first mature fuzzy DL reasoner, supports fuzzy $\mathcal{SHIF}(\mathbf{D})$ (a fragment of OWL~2) under multiple fuzzy logics and offers a Protégé~\cite{protege,ProtegePaper} plug-in in the \texttt{Fuzzy OWL~2} format~\cite{BOBILLO20111073}. DeLorean further improves the landscape by providing a complete reduction from \texttt{Fuzzy OWL~2} to crisp OWL~2, enabling the reuse of standard DL reasoners such as Hermit~\cite{HermiT}. However, each reasoner utilises its personally unique syntax, emphasising the necessity for a standardised representation.
    
    A promising approach is given by \texttt{Fuzzy OWL~2}, which encodes fuzzy axioms via OWL~2 annotation properties, thus remaining compatible with existing tools such as Protégé. Annotation-based encoding also facilitates the effective utilisation of crisp DL tools, allowing them to operate effectively in the presence of fuzzy components, although not using the fuzzy part. Alternative strategies may be considered to involve the extension of OWL with fuzzy constructors, as previously proposed in~\cite{Gao2005562,Stoilos2010656,stoilos2007extending}.

    Despite these efforts, a widely accepted standard for representing fuzzy ontologies is still unrealised. Therefore, it is imperative to bridge this gap to facilitate interoperability, streamline deployment, and leverage existing OWL-based workflows without compromising theoretical rigour or computational tractability. It is essential that future work converges on a consensus representation, building on annotation-based approaches such as \texttt{Fuzzy OWL~2} or careful OWL language extensions. Moreover, it is crucial to ensure broad tool support, ideally through extensions to OWL~2 or the standardisation of fuzzy DL profiles.

    In this work, we describe \texttt{fuzzy-dl-owl2}, an extension and a porting in Python of \texttt{fuzzyDL}~\cite{BobilloStraccia2016} and \texttt{Fuzzy OWL~2}~\cite{BOBILLO20111073}, originally developed by U. Straccia and F. Bobillo. In particular, the following aspects have been addressed in the porting process:
    \begin{itemize}
        \item the resolution of some semantic problems;
        \item the support of multiple Mixed-Integer Linear Programming (MILP) solvers;~\footnote{\texttt{fuzzyDL} supports only Gurobi.} and
        \item the re-engineering of the code to provide both a better architecture and ease of extension.
    \end{itemize}

    \nd In the following, we proceed as follows. In \Cref{sec:background}, we provide a brief description of fuzzy DLs and \texttt{Fuzzy OWL~2} language to give the reader the necessary background notions. 
    In \Cref{sec:related}, we provide the related works on the topic.
    In \Cref{sec:porting}, we describe the process of re-engineering in the Python language of \texttt{fuzzyDL} and \texttt{Fuzzy OWL~2}, along with the resolutions of some problems, while in \Cref{sec:architecture} we describe the architecture of \texttt{fuzzy\_dl\_owl2}. 
    In \Cref{sec:applications}, a concise illustration of the library's practical application is shown.  
    In Section~\ref{sec:benchmarks} a comparison between the Java and Python codebases through a series of benchmark experiments is presented.
    Finally, in \Cref{sec:conclusions}, we provide the conclusions of this work and illustrate some topics for future work.
    We refer the reader to the appendixes for further details on the grammar (\ref{sec:grammar}) and on the parser (\ref{sec:appendix-b}).

    \section{Background}\label{sec:background}

\nd This section provides the necessary background on fuzzy DLs and \texttt{Fuzzy OWL~2}.

    \subsection{Fuzzy Description Logic}\label{sec:fuzzy-dl}

    \nd \emph{Description Logics} (DLs) are a family of formal languages used for knowledge representation~\cite{Baader_Calvanese_McGuinness_Nardi_Patel-Schneider_2007}. \emph{Fuzzy Description Logics} (Fuzzy DLs, or FDLs)~\cite{Borgwardt2017,CrossChen2018,LUKASIEWICZ2008291} extend classical DLs using fuzzy set theory and fuzzy logic~\cite{zadeh1965fuzzy} to handle vagueness and imprecision. These are characteristics commonly found in real-world knowledge, where many domains cannot be accurately described using only crisp, binary classifications; instead, a degree of membership or truth is often more appropriate. Examples of such concepts are: a young and rich person, a big animal, a nice beach, heavy rain, a hot day, a good, dry, but expensive wine, etc.
    
\begin{figure}[!hbt]
        \centering
        \includegraphics[width=0.9\linewidth]{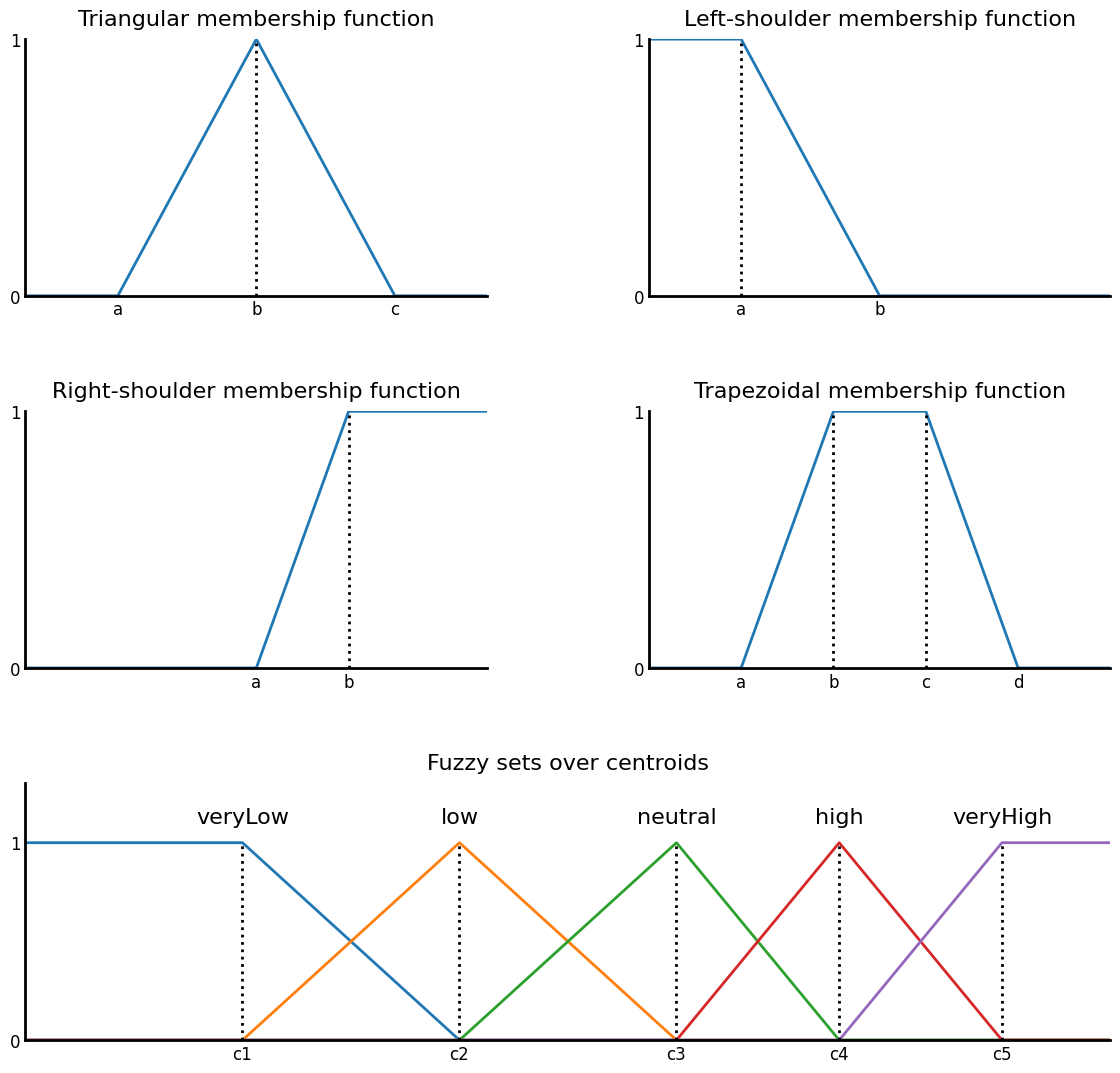}
        \caption{Membership functions used for fuzzy modelling: triangular function $\mathit{tri}(a,b,c)$, left-shoulder function $\mathit{ls}(a,b)$, right-shoulder function $\mathit{rs}(a,b)$, trapezoidal function $\mathit{trz}(a,b,c,d)$, and fuzzy linguistic sets defined over the centroids $c_1,\ldots,c_5$.}
        \label{fig:muf}
    \end{figure}

    \texttt{fuzzyDL} is a reasoning system designed to work with fuzzy DLs. It supports expressive constructs from classical DLs and enhances them with fuzzy set theory, allowing for the representation of concepts with a degree of membership. It operates over a discretised domain of truth values\footnote{A discretised domain is needed to ensure decidability~\cite{Penaloza2015}.} and includes support for concrete data types such as integers and strings. The system allows users to define fuzzy concepts using common membership functions (e.g., triangular, trapezoidal, and shoulder functions, see \Cref{fig:muf}). Let us recall that one easy method to define the membership functions is to uniformly partition the values of an attribute, for instance, the price of hotel rooms, into 3, 5, or 7 fuzzy sets. Another popular approach may consist of using the \emph{c-means} fuzzy clustering algorithm (see, e.g.~\cite{Bezdek81}) with 3, 5, or 7 clusters, where the fuzzy membership functions are shaped as triangular or shoulder functions built around the centroids of the clusters~\cite{Cardillo24,Cardillo22,Datil} (see \Cref{fig:muf}). In \cite{Huitzil20}, instead, aggregation operators are employed to combine (define the consensus among) fuzzy membership functions proposed by a group of experts. 
    %

    %
    %
    
    A peculiarity of the \texttt{fuzzyDL}  system is that it incorporates a reasoning engine based on a hybrid approach combining tableau algorithms and \emph{Mixed Integer Linear Programming} (MILP)~\cite{salkin1989foundations}, which enables it to compute satisfiability, subsumption, best entailment degree and other inference tasks in the presence of fuzzy axioms.
    
    The system also includes a query language for computing membership degrees, optimising variable values, and even performing \emph{defuzzification}~\cite{Klir95}. It also supports the import of OWL ontologies to facilitate the integration with existing Semantic Web frameworks.

    \subsection{Fuzzy OWL~2} \label{sec:fuzzy-owl-2}
    \nd In~\cite{BOBILLO20111073}, the authors address the need to represent vague or imprecise knowledge within the standard Semantic Web languages. Traditional ontology languages, such as OWL and OWL~2~\cite{OWL2,OWL2Overview}, have not been designed to handle fuzzy information. But, rather than proposing an entirely new fuzzy extension of OWL~2, the authors present a method to incorporate fuzzy information directly into OWL~2 using its \emph{annotation property mechanism}. This allows fuzzy ontologies to be created and edited with existing tools while maintaining compatibility with standard OWL~2 infrastructure.

    The approach focuses on defining fuzzy concepts, roles, axioms, and datatypes through XML annotations, which specify degrees of membership, fuzzy modifiers, and membership functions. They propose a structured format for these annotations using the \texttt{fuzzyLabel} annotation property, enabling fuzzy semantics to be attached to otherwise crisp OWL elements without altering the underlying logic of OWL~2. The paper details how different types of fuzzy elements--like modified concepts, weighted concepts, and fuzzy roles--can be syntactically encoded using this annotation mechanism. 

    In subsequent papers, the authors proposed small extensions of \texttt{Fuzzy OWL~2} to support new features, e.g., aggregation operators~\cite{ASOC2013}.
    
    The authors also developed a parser to translate the annotated OWL~2 files into the languages supported by the
    \texttt{fuzzyDL} system (described in \Cref{sec:fuzzy-dl}) and the DeLorean reasoner.
    
    \section{Related Works}\label{sec:related}
    \nd For some time, the de facto standard API to manage classical ontologies within Java applications is the OWL API~\cite{OWLAPI}. It makes it possible to load, create, and update OWL ontologies, but also to interact with those ontology reasoners which implement the OWL API interfaces, such as HermiT~\cite{HermitReasonerWebsite,HermiT} and Pellet~\cite{PelletReasonerGithub,SIRIN200751}. In order to manage fuzzy ontologies, \texttt{Fuzzy OWL~2} provides parsers to translate OWL~2 annotations encoding fuzzy knowledge into different formats, and \texttt{fuzzyDL}  reasoner provides a Java API to query ontologies.
    
    However, the management of ontologies in Python has not received similar attention so far. To the best of our knowledge, this is the first effort towards a full Python support for reasoning with and representing fuzzy ontologies in OWL~2. To manage classical ontologies in Python, it is worth mentioning the Python library \texttt{owlready2}, which, however, does not completely address OWL~2 and has a quite different usage objective.
    Specifically, in \texttt{owlready2} it is possible to insert basic annotations, but not more complex ones such as those associated with \texttt{General Class Axioms}, i.e., axioms of the form \texttt{C is a subclass of B}, where \texttt{C} is not a named class, but a complex concept built using class constructs such as the intersection of classes or a union of classes. 
    
    Apart from the fact that \texttt{owlready2} does not support fuzzy ontologies, the main differences between \texttt{owlready2} and our approach are that it 
    \begin{enumerate}
        \item[(i)] treats ontology classes, instances, and properties directly as standard Python objects, which we do not do by purpose; and
        \item[(ii)] integrates the HermiT and Pellet reasoners (without the need for a separate Java bridge), whereas we rely on third-party MILP solvers.
    \end{enumerate}

    \section{Porting to Python}\label{sec:porting}
    
    \nd In this section, we present our Python port of the \texttt{fuzzyDL} reasoner and of the \texttt{Fuzzy OWL~2} language infrastructure, in order to improve efficiency, portability, extensibility, and, ultimately, usability. We refer the reader to \ref{sec:appendix-b} for a complete description of the objects generated by the parser for each \texttt{fuzzyDL} directive.
    For the system architecture description, we refer the reader to Section~\ref{sec:architecture}.

    The entire source code of the project is available from the  \texttt{fuzzy-dl-owl2} GitHub repository \cite{fuzzydlowl2_github}, while the documentation is available from the \texttt{Read the Docs} site \cite{fuzzydlowl2_docs}. For further details and all the UML diagrams, the reader is referred to the section entitled \texttt{API Reference} in the official documentation \cite{fuzzydlowl2_docs}. Specifically, each module and class is accompanied by a UML diagram located at the top of the page and a description of the relevant module or class.

    The online documentation was generated using the \texttt{sphinx-autoapi} Python package  for automatic API generation. For more details, we refer the reader to the directory \texttt{docs} of the \texttt{GitHub} repository \cite{fuzzydlowl2_github}.

    The UML diagrams presented in the figures have been generated using the \texttt{pyreverse} Python package. It is noteworthy that all private~\footnote{In the context of Python object-oriented programming, it is important to note that there are no private methods or attributes. Conventionally, these elements are designated with at least an underscore at the beginning of the name. For instance, the method \texttt{foo} of the class \texttt{Foo} can be defined with the name \texttt{\_\_foo} if it is a \texttt{private} method. However, it should be noted that this method will still be accessible through the instances of the class \texttt{Foo} with the name \texttt{\_Foo\_\_foo}. For more details, we refer the reader to \cite{python_private_name_mangling}.} methods and attributes have been excluded from the UML diagrams. Abstract methods, on the other hand, will be indicated in italics.

    \subsection{FuzzyDL}

    \nd In this section, we discuss the main changes in our Python implementation related to the Java-based \texttt{fuzzyDL} reasoner.

    A major change comes from the re-engineering of \texttt{fuzzyDL} 's \texttt{Concept} class. 
    In the original Java implementation, the class \texttt{Concept} contains most of the logic for the implemented DL concepts. While this facilitated the use of that class in other Java classes, such as the class \texttt{KnowledgeBase}, it also includes by far too many details, i.e., attributes, even if these are not necessary for the concept's usage itself. 
    
\begin{figure}
        \centering
        \includegraphics[width=0.67\linewidth]{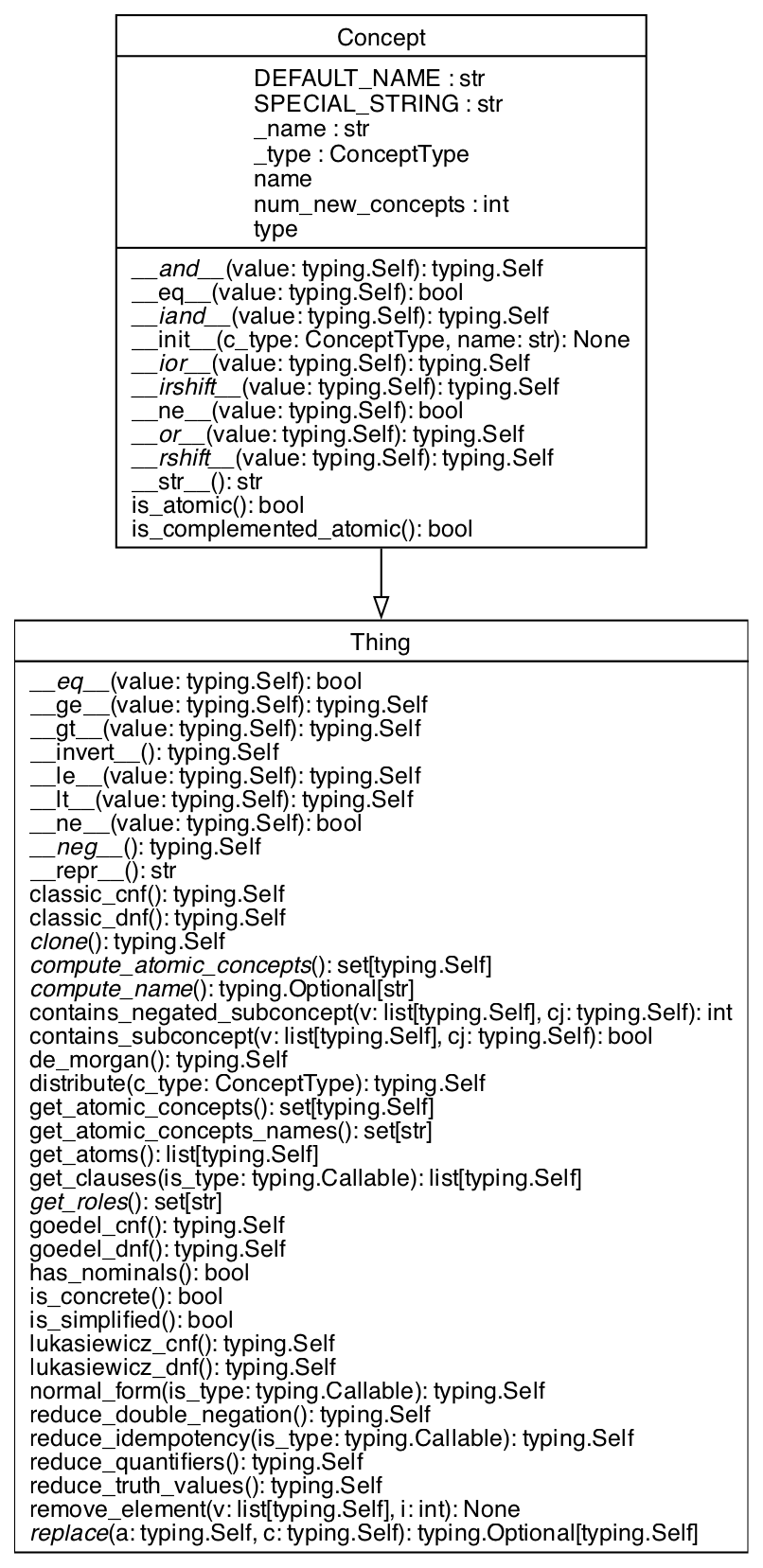}
        \caption{UML Class Diagram of the class \texttt{Concept}.}
        \label{fig:concept-class}
    \end{figure}

\begin{figure}
        \centering
        \caption{UML Class Diagram of the class \texttt{OperatorConcept}.}
        \label{fig:operator-concept-class}
        \includegraphics[width=0.45\linewidth]{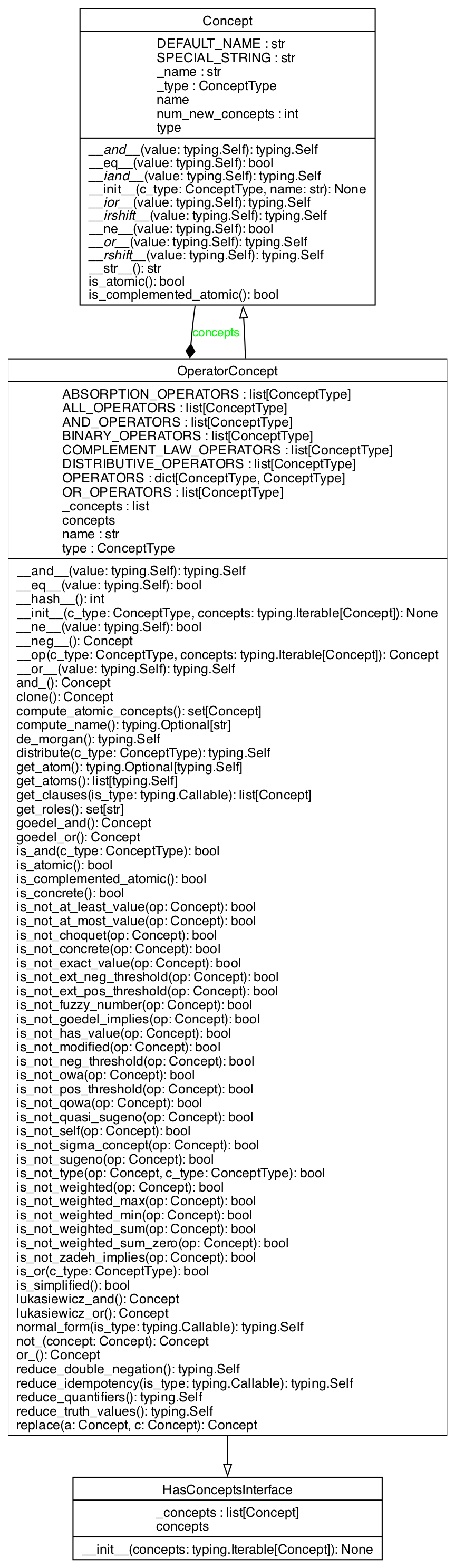}
    \end{figure}

    For example, a \emph{concept conjunction} only needs the list of concepts it logically connects, whereas the \emph{weighted concept}, i.e., the scalar product between a real constant and a concept, only needs the pair (scalar, concept). The original implementation, however, associates a list of concepts (empty if not needed) and a scalar weight (zero if not needed) with each concept implemented, even if the latter is not used. For this reason, in our current version in Python, the class \texttt{Concept} only contains some basic attributes and methods (see \Cref{fig:concept-class}), having specific sub-classes as illustrated in \Cref{fig:concept-class-diagram}. For instance, the class \texttt{OperatorConcept} (see \Cref{fig:operator-concept-class}) is defined as a subclass of \texttt{Concept} and is aimed to represent a set of concepts connected by logical connectives such as \texttt{AND}, \texttt{OR} and \texttt{NOT}. 

    Another major difference with the original implementation is the decoupling between the various classes and the class \texttt{KnowledgeBase}. Specifically, the original Java implementation contained several class methods with the \texttt{KnowledgeBase} class as a parameter in order to use its methods or, for instance,  to add MILP constraints in the form of equations or inequalities.
    Instead, in our implementation, we moved all the logic, i.e., attributes and methods, to the class \texttt{KnowledgeBase} to centralise its management and avoid circular class imports -- a problem to be avoided in a language like Python. 
    
\begin{figure}
        \centering
        \caption{UML Class Diagram of the class \texttt{IndividualHandler}.}
        \label{fig:individual-handler-class}
        \includegraphics[width=0.85\linewidth]{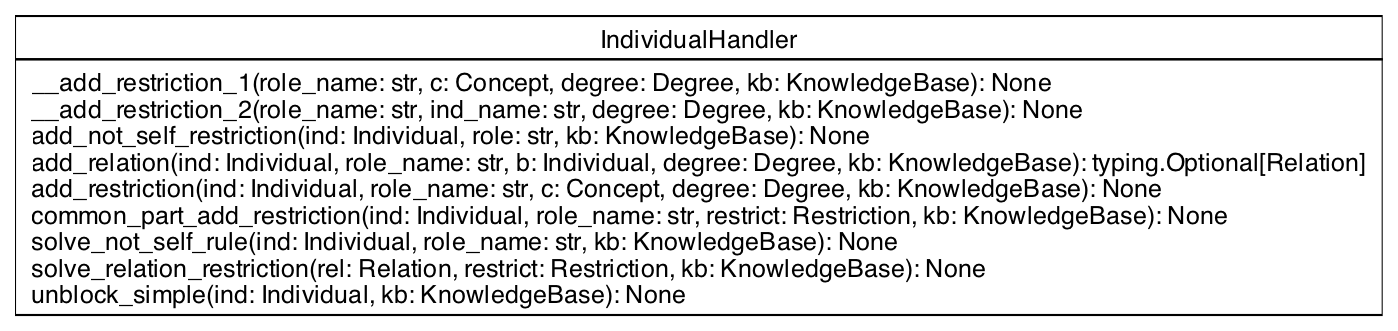}
    \end{figure}

\begin{figure}
        \centering
        \caption{UML Class Diagram of the class \texttt{CreatedIndividualHandler}.}
        \label{fig:created-individual-handler-class}
        \includegraphics[width=0.95\linewidth]{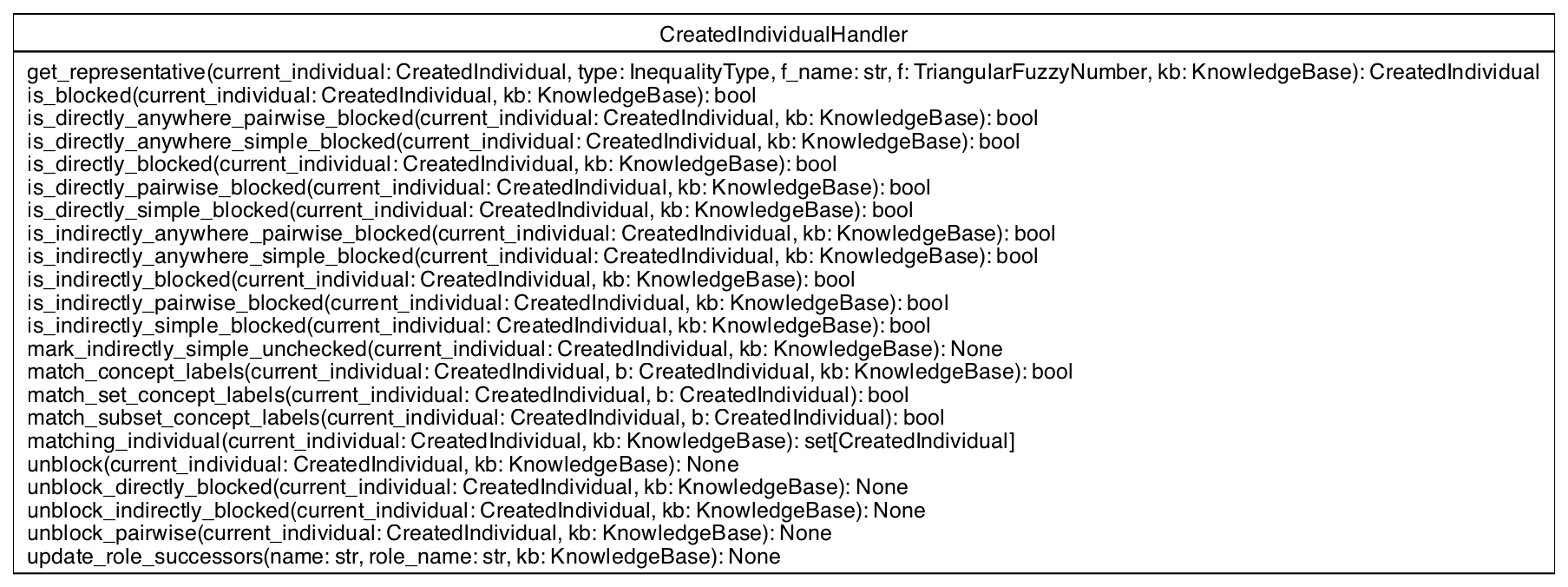}
    \end{figure}

    Specifically, we created the classes \texttt{IndividualHandler} (see \Cref{fig:individual-handler-class}) and \linebreak \texttt{CreatedIndividualHandler} (see \Cref{fig:created-individual-handler-class}) by spinning off the methods that referenced the  \texttt{KnowledgeBase} class into the classes \texttt{Individual} and \texttt{CreatedIndividual} provided in the original implementation. 
    In order to avoid circular loops for the importation of the various classes strongly coupled with the \texttt{KnowledgeBase}, we have enclosed the following classes in the same file: \texttt{KnowledgeBase}, \texttt{DatatypeReasoner} (i.e., the class that handles datatypes and the set of MILP constraints for each of them), \texttt{ClassicalSolver} (i.e., the solver used in case classical crisp logic is used), \texttt{ZadehSolver} (i.e., the solver used in case Zadeh logic is used), \texttt{LukasiewiczSolver} (i.e., the solver used in case \L{}ukasiewicz logic is used), \texttt{IndividualHandler}, and \texttt{CreatedIndividualHandler}.

    A final major difference with the original implementation consists of the addition of new MILP solvers. In particular, the original implementation only allows the use of \texttt{Gurobi Optimizer}~\cite{GurobiOptimizer,GurobiReferenceManual} V8.1.0. In our implementation, additional MILP solvers have been added: \texttt{CPLEX}~\cite{ibm_CPLEX_optimizer}, \texttt{GNU Linear Programming Kit (GLPK)}~\cite{glpk_gnu}, \texttt{COIN Branch and Cut (CBC)}~\cite{cbc_user_guide}, \texttt{HiGH performance Software for linear optimisation (HiGHS)}~\cite{highs_solver}. In particular, three Python libraries were used for the MILP solvers: \texttt{gurobipy v12.0.0}~\cite{gurobipy_12_0_0}, \texttt{mip v1.16rc0}~\cite{mip_1_16rc0}, and \texttt{PuLP v3.2.1}~\cite{pulp_3_2_1}.

    \begin{adjustwidth}{-1cm}{-1cm}
        \begin{figure}
            \centering
            \caption{Representation of the UML Class Diagram of all \texttt{Concept} type classes. For a better graphical representation, all class attributes and methods have been omitted.}
            \label{fig:concept-class-diagram}
            \includegraphics[width=1\linewidth,angle=180]{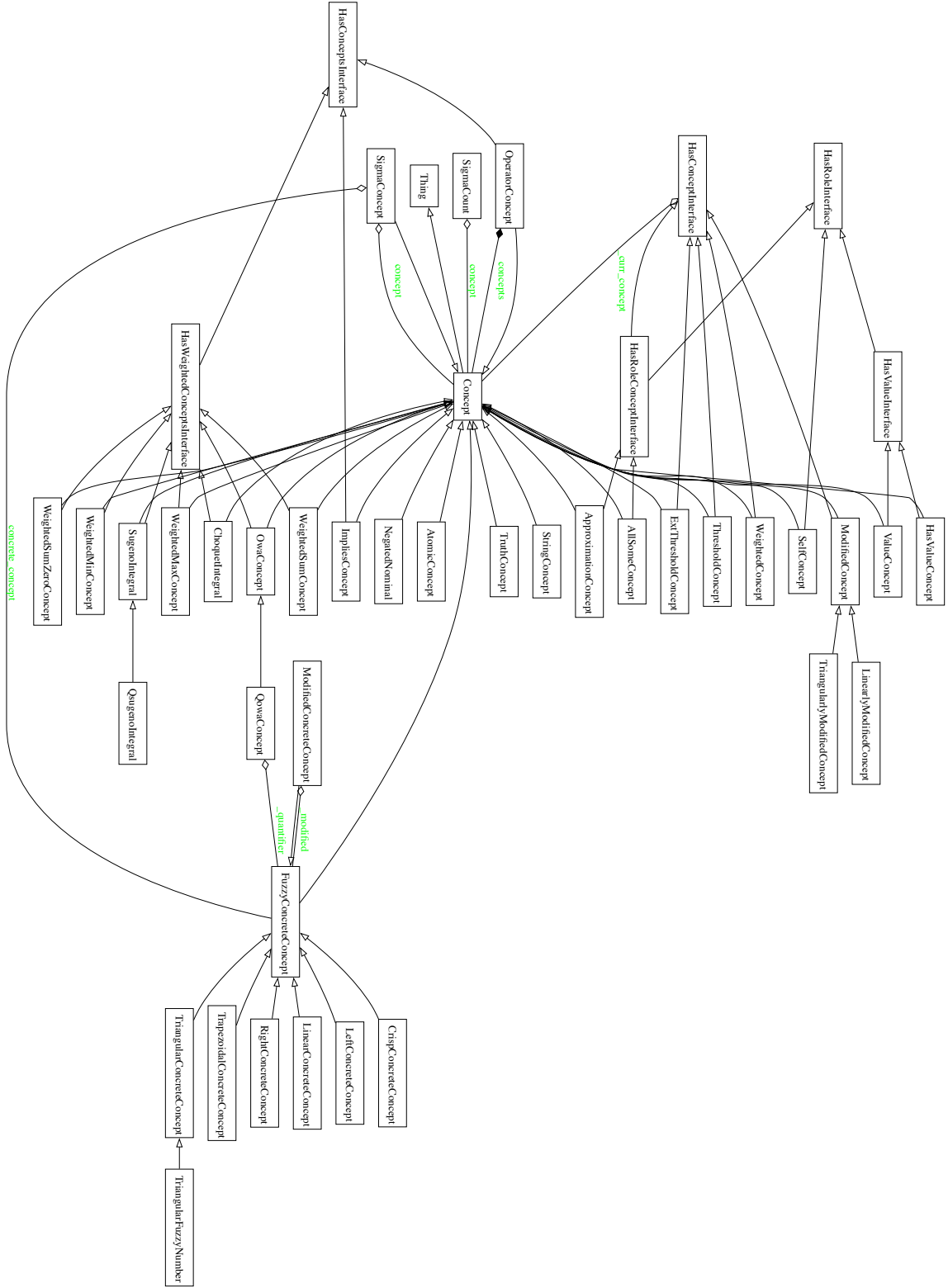}
        \end{figure}
    \end{adjustwidth}
    
    \subsection{Fuzzy OWL~2}

    \nd  In this section, we discuss the changes and improvements made to the implementation of the \texttt{Fuzzy OWL~2} component.

    \nd Specifically, the original implementation involves the use of XML annotations of the various fuzzy elements foreseen in an ontology (classes, properties, individuals, and so on). Currently, however, there is no library in Python that faithfully implements the OWL~2 standard provided by the W3C (see \Cref{sec:related} for details). Therefore, given the current limitations, we have implemented the Python \texttt{pyowl2} library \cite{pyowl2}, whose source code is available on \texttt{PyPI} \cite{pyowl2} and \texttt{GitHub} \cite{pyowl2_github}, and whose documentation is available on \texttt{Read the Docs} \cite{pyowl2_docs}.     
    At present, annotations of annotations have not been implemented, but they are not required by \texttt{Fuzzy OWL~2}.
    
    The biggest differences with the original implementation consist of resolving some semantic inconsistencies and violations of the standard. In particular, the original implementation does not consider the case where the name of a class could coincide with the name of an individual. For example, in the \texttt{fuzzyDL} ontology called \texttt{FuzzyWine.1.0} available at the download directory of the original \texttt{fuzzyDL} software \cite{fuzzydl_web}, there are individuals and classes with the same names. When generating the associated OWL ontology, two logically distinct objects (individuals and classes) with the same IRI (which, by definition, should be unique for each resource) are then generated.

    Furthermore, in the original implementation, the disjointness between \emph{Data Properties} and \emph{Object Properties} is not handled correctly. Specifically, a data property is sometimes converted to both an object property and a data property, which share the same IRI. In this case, the object property is used to capture domain information, while the data property encapsulates range information. This redundancy introduces ambiguity in the structure of the ontology.
    
    Instead, in our implementation, we addressed both issues. Specifically, specific sub-paths of the main ontology namespace were created for each element. That is, the sub-paths \texttt{object-property}, \texttt{data-property}, \texttt{individual}, \texttt{class}, \texttt{annotation-property}, and \texttt{datatype} were created for the insertion of object properties, data properties, individuals, classes, annotation properties, and datatypes, respectively. In this way, if two entities share the same name but have different semantic meanings (e.g., an individual and a class), then they will have distinct IRIs. 
    As for data properties, we made sure that they are created only as data properties with their respective domains and ranges, without creating object properties with the same name.

    A final important difference with the original implementation is the use of Python's standard \texttt{defusedxml} library for parsing string XML annotations in the translation phase from \texttt{fuzzyDL}  to \texttt{Fuzzy OWL~2}. In particular, an ad-hoc parser was created in the original implementation to define the grammar of the XML annotations used to define fuzzy concepts and properties in the ontology. This, however, can lead to several errors in defining the grammar due to the various special symbols contained in an XML (such as angle brackets and double quotes for strings). Python's \texttt{defusedxml} library, on the other hand, already contains within it an XML parser that has been tried and tested over the years and which simplifies the reading of XMLs. In particular, we used this library to transform the string XML annotation into a structured XML object that is much more machine-readable and parsable. We then use this dictionary to check whether the various XML tags adhere to the grammar and return the respective element in the form of a specific class instance.
    For instance, the fuzzy weighted sum concept of the form 
    \begin{center}
        \texttt{0.2 A + 0.3 B + 0.5 C}
    \end{center}
    
    \nd has the XML annotation described in \Cref{listing:fuzzy-owl-2}.

    \begin{listing}[!ht]
        \caption{An XML annotation for \texttt{Fuzzy OWL~2} ontology representing a fuzzy weighted sum concept of the form \texttt{0.2 A + 0.3 B + 0.5 C}.}
        \label{listing:fuzzy-owl-2}
        
\begin{lstlisting}[
            language=XML,
            frame=lines,
            framesep=2mm,
            basicstyle=\ttfamily\scriptsize\linespread{1.2}\selectfont,
            numbers=left
        ]
<fuzzyOwl2 fuzzyType="concept">
    <Concept type="weightedSum">
        <Concept type="weighted" value="0.2" base="A" />
        <Concept type="weighted" value="0.3" base="B" />
        <Concept type="weighted" value="0.5" base="C" />
    </Concept>
</fuzzyOwl2>
\end{lstlisting}        
    \end{listing}

\section{System Architecture}\label{sec:architecture}
    \nd This section describes the architecture of \texttt{fuzzy\_dl\_owl2}.
    The system supports a fuzzy extension of the description logic $\mathcal{SHIF}(\mathbf{D})$ under \L{}ukasiewicz, Zadeh, G\"odel, and classical semantics, and accepts ontologies in both the native \texttt{fuzzyDL} syntax and the annotation-based \texttt{Fuzzy OWL~2} encoding~\cite{BOBILLO20111073}.  Every reasoning task is reduced to a bounded MILP, following the tableau algorithm of Bobillo and Straccia~\cite{bobillo2008fuzzydl}, and dispatched to an external optimisation solver.
    
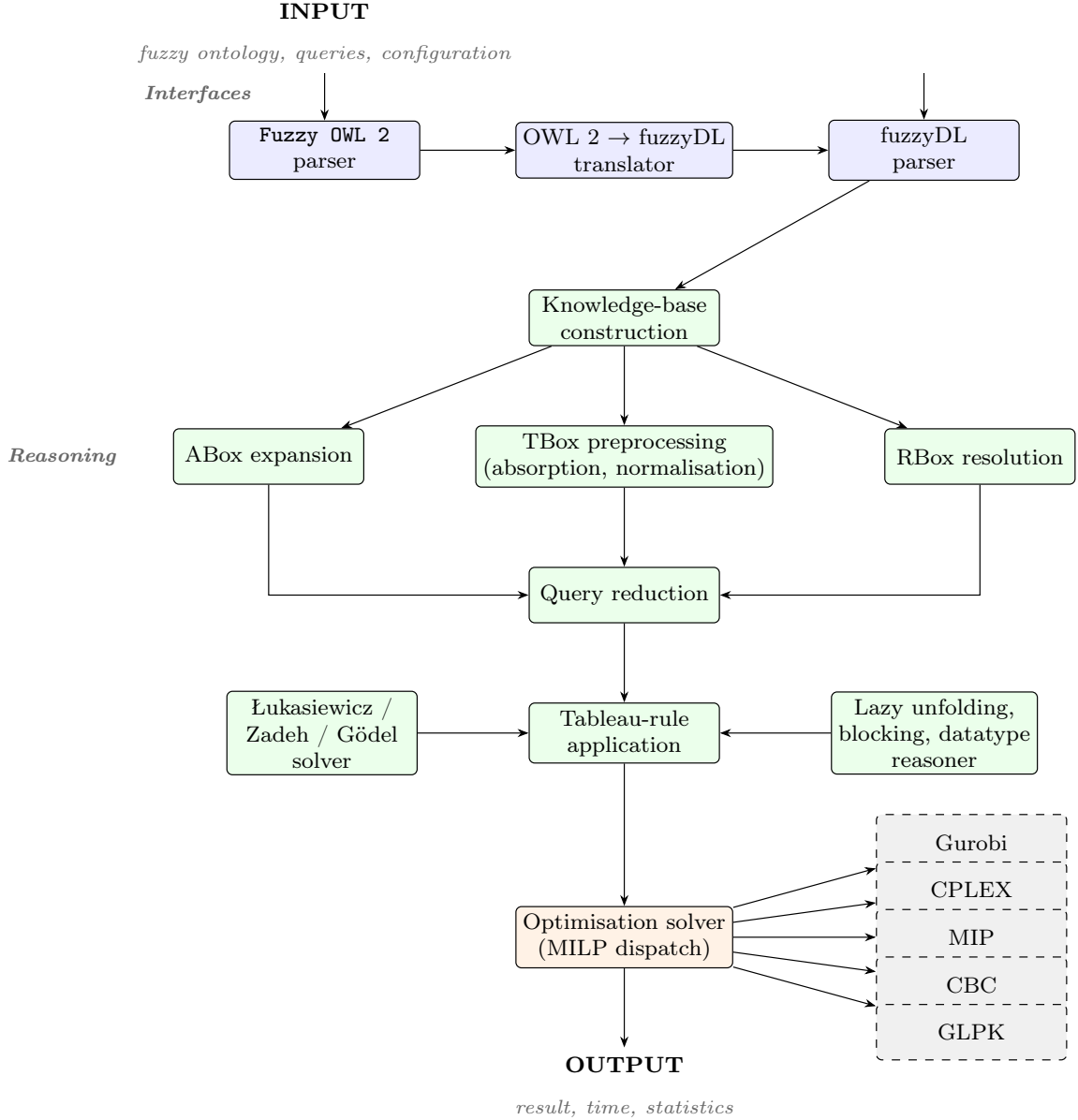
\begin{figure}[tbh!]
        \centering
        \caption{Architecture of the \texttt{fuzzy\_dl\_owl2} reasoner. Solid, rounded boxes are internal components; dashed boxes are external solver back-ends. Arrows denote data flow.}
        \label{fig:architecture}
        
        \resizebox{.95\textwidth}{!}{%
        \begin{tikzpicture}[
        >={Stealth[length=4pt]},
        node distance=6mm and 10mm,
        box/.style  = {draw, rounded corners=2pt, minimum height=7mm,
                        minimum width=24mm, font=\footnotesize, align=center,
                        inner sep=2.5pt},
        iface/.style  = {box, fill=blue!8},
        reason/.style = {box, fill=green!8},
        solver/.style = {box, fill=orange!10},
        ext/.style    = {box, fill=gray!12, dashed},
        lbl/.style    = {font=\scriptsize\itshape, text=black!60},
        io/.style     = {font=\footnotesize\bfseries},
        ]
    
        \node[io] (input) {INPUT};
        \node[lbl, below=1mm of input] (inputdesc)
            {fuzzy ontology, queries, configuration};
    
        \node[iface, below=6mm of inputdesc] (owlparser)
            {\texttt{Fuzzy OWL~2}\\parser};
        \node[iface, right=12mm of owlparser] (owl2dl)
            {OWL~2 $\to$ fuzzyDL\\translator};
        \node[iface, right=12mm of owl2dl] (dlparser)
            {fuzzyDL\\parser};
    
        \node[lbl, above left=1mm and -4mm of owlparser] {\textbf{Interfaces}};
    
        \node[reason, below=14mm of owl2dl] (kbcreate)
            {Knowledge-base\\construction};
    
        \node[reason, below=10mm of kbcreate] (preproc)
            {TBox preprocessing\\(absorption, normalisation)};
    
        \node[reason, left=14mm of preproc] (aboxexp)
            {ABox expansion};
        \node[reason, right=14mm of preproc] (rboxres)
            {RBox resolution};
    
        \node[reason, below=10mm of preproc] (qreduce)
            {Query reduction};
    
        \node[reason, below=10mm of qreduce] (tableau)
            {Tableau-rule\\application};
    
        \node[reason, left=14mm of tableau] (fuzzysolv)
            {\L{}ukasiewicz /\\Zadeh / G\"odel\\solver};
        \node[reason, right=14mm of tableau] (techniques)
            {Lazy unfolding,\\blocking, datatype\\reasoner};
    
        \node[lbl, left=6mm of aboxexp] {\textbf{Reasoning}};
    
        \node[solver, below=18mm of tableau] (milp)
            {Optimisation solver\\(MILP dispatch)};
    
        \node[ext, right=18mm of milp, yshift= 12mm] (gurobi){Gurobi};
        \node[ext, right=18mm of milp, yshift= 6mm] (CPLEX) {CPLEX};
        \node[ext, right=18mm of milp, yshift= 0mm] (mip) {MIP};
        \node[ext, right=18mm of milp, yshift=-6mm] (cbc) {CBC};
        \node[ext, right=18mm of milp, yshift=-12mm] (glpk) {GLPK};
    
        \node[io, below=10mm of milp] (output) {OUTPUT};
        \node[lbl, below=1mm of output] (outputdesc)
            {result, time, statistics};
    
        \draw[->] (inputdesc.south -| owlparser) -- (owlparser);
        \draw[->] (inputdesc.south -| dlparser)  -- (dlparser);
        \draw[->] (owlparser) -- (owl2dl);
        \draw[->] (owl2dl)    -- (dlparser);
        \draw[->] (dlparser)  -- (kbcreate);
        \draw[->] (kbcreate)  -- (preproc);
        \draw[->] (kbcreate)  -- (aboxexp);
        \draw[->] (kbcreate)  -- (rboxres);
        \draw[->] (preproc)   -- (qreduce);
        \draw[->] (aboxexp)   |- (qreduce);
        \draw[->] (rboxres)   |- (qreduce);
        \draw[->] (qreduce)   -- (tableau);
        \draw[->] (fuzzysolv) -- (tableau);
        \draw[->] (techniques)-- (tableau);
        \draw[->] (tableau)   -- (milp);
        \draw[->] (milp)      -- (gurobi);
        \draw[->] (milp)      -- (CPLEX);
        \draw[->] (milp)      -- (mip);
        \draw[->] (milp)      -- (cbc);
        \draw[->] (milp)      -- (glpk);
        \draw[->] (milp)      -- (output);
    
        \end{tikzpicture}%
        }
    \end{figure}

    In \Cref{fig:architecture}, we illustrate the overall pipeline. The upper part manages input parsing and KB construction; the lower part implements the tableau-based reasoning engine.
    
    \subsection{Input Interfaces}
    \label{sec:arch:interfaces}

\nd The system accepts fuzzy ontologies through two input paths (cf.\ the top of \Cref{fig:architecture}).

    \begin{description}
        \item[Native fuzzyDL syntax:]  A hand-written recursive-descent parser tokenises the input into an S-expression stream and dispatches each form to a semantic callback that directly populates the KB.  The parser avoids backtracking and packrat caching and is designed for optional Cython \cite{Cython} compilation, making it competitive with the JavaCC-generated parser of the original Java implementation.
        \item[Fuzzy OWL~2:] The \texttt{Fuzzy OWL~2} front-end reads an OWL~2 ontology via the \texttt{pyowl2} library and extracts fuzzy annotations---membership functions, modifiers, aggregation-concept definitions, and fuzzy property axioms---into a typed intermediate representation. Both RDF/XML and the native \texttt{Fuzzy OWL~2} XML serialisation~\cite{BOBILLO20111073} are supported.\footnote{RDF/XML is the only OWL~2 syntax supported by \texttt{pyowl2}.}. A translator converts this intermediate form into \texttt{fuzzyDL} syntax, which is then parsed by the same recursive-descent front-end.  A reverse translator enables round-trip conversion from \texttt{fuzzyDL} to OWL~2.
    \end{description}
    
    \subsection{Knowledge Base Representation}\label{sec:arch:kb}
    
    \nd The parser populates a central \texttt{KnowledgeBase} object that maintains:
    \begin{itemize}
      \item a Terminological Box (\emph{TBox}) composed of graded concept definitions, primitive inclusions, and general concept inclusions (GCIs) (a class is a subclass of another class to some degree);
      \item an Assertional Box (\emph{ABox}) composed of fuzzy concept assertions (an individual is an instance of class some degree), role assertions, and concrete-domain value assertions;
      \item a Role Box (\emph{RBox}) composed of graded role hierarchy (a role is a subrole of another role to some degree), transitivity, symmetry, reflexivity, functionality, inverse-role declarations, and domain/range restrictions.
    \end{itemize}
    
    \nd Concepts are modelled by a class hierarchy rooted at an abstract \texttt{Concept} class.  The hierarchy covers the standard $\mathcal{SHIF}(\mathbf{D})$ constructors, i.e., conjunction, disjunction, negation, existential and universal restrictions, threshold, implication, as well as fuzzy-specific constructors such as:
    
    \begin{itemize}
      \item \emph{Concrete-domain membership functions:} trapezoidal, triangular, left-shoulder, right-shoulder, linear, crisp, modified, and fuzzy-number shapes;
      \item \emph{Modifiers:} linear and triangular families, supporting \emph{hedges} such as \emph{very} and \emph{slightly};
      \item \emph{Aggregation operators:} weighted sum, strict weighted sum, weighted min, weighted max, ordered weighted averaging, quantified-guided ordered weighted averaging, Choquet integral, Sugeno integral, quasi-Sugeno integral, and $\sigma$-count.
    \end{itemize}
    
    \nd Membership degrees are represented uniformly by a \texttt{Degree} hierarchy with three variants: \texttt{DegreeNumeric} (a constant in $[0,1]$), \texttt{DegreeVariable} (a symbolic degree linked to a MILP decision variable), and \texttt{DegreeExpression} (a linear combination of degree variables).  This allows axioms such as an individual $a$ is an instance of a concept to a degree greater than or equal to $d$, where $d$ is a number, a variable, or a linear expression, to be optimised.
    
    \subsection{Reasoning Pipeline}\label{sec:arch:reasoning}
    
    \nd The reasoning procedure mirrors the algorithm of Bobillo and Straccia~\cite{BobilloStraccia2016} and proceeds in four stages.
    
    \subsubsection{Preprocessing}\label{sec:arch:preproc}
    
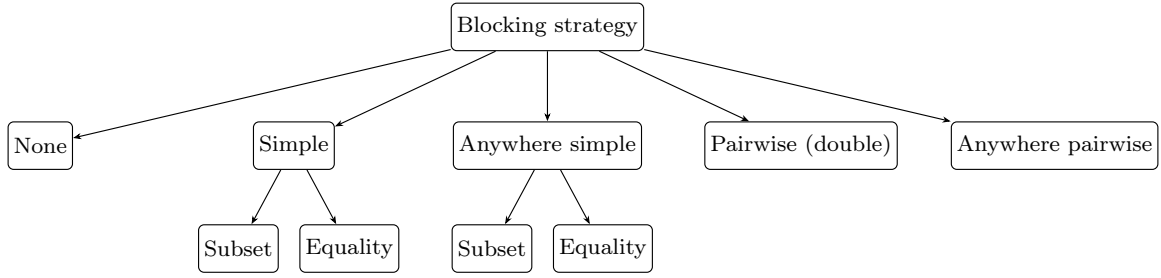
\begin{figure}[tbh!]
        \centering
        \label{fig:blocking}
        \resizebox{0.95\textwidth}{!}{%
        \begin{tikzpicture}[
            >={Stealth[length=3pt]},
            edge from parent/.style={draw,->},
            level 1/.style={sibling distance=32mm, level distance=15mm},
            level 2/.style={sibling distance=14mm, level distance=13mm},
            every node/.style={draw, rounded corners=2pt, minimum height=6mm,font=\footnotesize, align=center, inner sep=2pt},
            ]
            \node {Blocking strategy}
            child { node {None} }
            child { node {Simple}
                child { node {Subset} }
                child { node {Equality} }
            }
            child { node {Anywhere simple}
                child { node {Subset} }
                child { node {Equality} }
            }
            child { node {Pairwise (double)} }
            child { node {Anywhere pairwise} }
            ;
        \end{tikzpicture}}
        
        \caption{Blocking strategies implemented in the reasoner. Each strategy may operate in static or dynamic mode.  The choice is made automatically based on ontology expressivity.}
    \end{figure}

    \nd Before tableau expansion, a KB undergoes several normalisation and optimisation passes:
    
    \begin{enumerate}
      \item \textbf{String encoding.} String-typed concrete features are mapped to consecutive integers so that the MILP solver operates exclusively over numeric domains;
      \item \textbf{RBox resolution.} Inverse roles, role inclusion axioms, reflexive roles, and functional roles are materialised;
      \item \textbf{TBox absorption.} If the TBox is \emph{lazy unfoldable} (acyclic, with no GCIs beyond definitions and primitive inclusions), axioms are segregated into definitions and inclusions for efficient lazy unfolding.  Otherwise, a multi-phase absorption loop iteratively rewrites GCIs 
      until a fixed point is reached;
      \item \textbf{Blocking-strategy selection.} Based on the presence of inverse, transitive, and functional roles and the acyclicity of the TBox, the system automatically selects the least expensive blocking strategy that guarantees termination (\Cref{fig:blocking});
      \item \textbf{Query reduction.}  Each user query is normalised to a single MILP optimisation target.
    \end{enumerate}
    
    \subsubsection{Tableau Expansion and Fuzzy-Logic Encoding}
    \label{sec:arch:tableau}
    
    \nd The tableau algorithm is deterministic: every applicable rule fires exactly once, producing a single completion whose constraints are collected into a unique MILP instance.  For each fuzzy assertion in the ABox, the tableau decomposes the involved concept according to its top-level constructor, generating appropriate MILP variables and inequalities.  The semantics of conjunction ($\otimes$), disjunction ($\oplus$), implication ($\Rightarrow$), and negation ($\ominus$) depend on the chosen fuzzy logic. In \Cref{tab:operators}, we summarised the operator definitions; each is reified into linear constraints via Big-M linearisations (which, essentially, consists of restricting to greater-than inequalities) with auxiliary binary variables. 

    \input{logic_operators}
    
    \nd Three dedicated solver classes --- \texttt{LukasiewiczSolver}, \texttt{ZadehSolver}, and \linebreak \texttt{ClassicalSolver} --- encapsulate these encodings. The classical solver restricts all variables to $\{0,1\}$, recovering crisp two-valued reasoning as a special case.
    
    A dedicated \texttt{DatatypeReasoner} translates concrete-domain assertions involving numeric comparisons, arithmetic feature functions (sum, difference, and product), and fuzzy-number membership into additional MILP constraints. \emph{Lazy unfolding} is applied to primitive concept inclusions and domain/range axioms: an axiom is expanded only for individuals known to instantiate the relevant atomic concept, reducing the size of the generated constraint system.
    
    \subsubsection{MILP Formulation and Solver Dispatch}\label{sec:arch:milp}
    
    \nd The constraint system is represented by a four-level compositional algebra:
    \begin{center}
    \scalebox{0.9}{
        \texttt{Variable} $\;\to\;$ \texttt{Term} ($c \cdot x$)
        $\;\to\;$ \texttt{Expression} $\left( 
        \sum_i c_i x_i\right)$
        $\;\to\;$ \texttt{Inequation} $\left(E 
        \geq \mid \leq 0\right)$
    }
    \end{center}
    
    \nd The \texttt{MILPHelper} class provides the sole interface for allocating decision variables and accumulating constraints; it enforces a canonical form in which the right-hand side of every inequation is zero.
    
    Once the tableau has generated the full constraint system, the helper dispatches it to one of four supported back-ends: Gurobi, CBC, the Python-MIP library, or PuLP (with GLPK, HiGHS, or CPLEX as the underlying engine).  The back-end is selected at configuration time.  Each query is answered by minimising or maximising an objective expression, e.g.~$x_{a:C}$, where $x_{a:C}$ is the decision variable representing the membership degree of being individual~$a$ instance of the concept~$C$.

    \subsection{Performance Engineering of the Front-End}\label{sec:arch:perf}

    \nd Parsing throughput is a practical bottleneck for large ontologies, so the front-end is engineered to approach the speed of the JavaCC-generated parser of the original Java reasoner while remaining pure-Python-installable. Two complementary techniques are used: \emph{(i)} ahead-of-time Cython compilation of the recursive-descent parser, and \emph{(ii)} a multi-backend tokeniser that selects the fastest scanner available in the running environment.

    \subsubsection{Cython Compilation}\label{sec:arch:perf:cython}

    \nd The hot path of the front-end --- the recursive-descent parser (\texttt{dl\_parser\_fast}) and its semantic callbacks --- is compiled ahead of time to a native extension through Cython~\cite{Cython} so that the inner scanning and dispatch loops run as plain C. These parser modules are always cythonised when a C toolchain is present; the Cython tuple builder for the tokeniser (see \Cref{sec:arch:perf:tokenizer}) is the only extension that is conditional, as it links the C scanner emitted by \texttt{re2c}\cite{re2c} or \texttt{flex}\cite{flex} and is therefore skipped when neither generator is installed. Compilation is \emph{optional}: when no C toolchain is available, the same \texttt{.py} sources are imported and interpreted, yielding identical behaviour at lower speed, which keeps a source-only \texttt{pip install} functional while rewarding environments that can build the extensions.

    \subsubsection{Multi-Backend Tokeniser}\label{sec:arch:perf:tokenizer}

    \nd Tokenisation is dispatched at run time to the fastest backend that the current process can import. All backends emit the \emph{same} \texttt{(kind, value, lower, offset)} 4-tuple stream terminated by a single \texttt{T\_EOF} sentinel, so the recursive-descent parser is agnostic to which scanner produced its input. Three backends are provided, in decreasing order of speed:

    \begin{enumerate}
      \item \textbf{Compiled \texttt{re2c}/\texttt{flex} scanner.} A single master  (\texttt{generate.py}) keyword table is the sole source of truth from which the lexer specification, the C token header, and the Python token codes are generated, keeping every layer in sync. At build time the specification is compiled by \texttt{re2c}~\cite{re2c} (preferred) or, if absent, by \texttt{flex}, producing a branchless C scanner that fills three parallel \texttt{int32} span arrays (token kind, byte start, byte length) without allocating any Python objects. A thin Cython module (\texttt{tokenize\_tuples}) then converts those spans into the parser's 4-tuples in C under a released Global Interpreter Lock (GIL), so both the scan and the tuple construction avoid the interpreter. This tuple builder is compiled only when one of the two generators is installed; otherwise, it is omitted, and the system falls back to the backends below.
      \item \textbf{CFFI array tokeniser.} When the Cython tuple builder is unavailable but the C scanner is, the same span arrays are exposed through a \textit{C Foreign Function Interface} (CFFI) wrapper and converted to 4-tuples in Python --- slower than the Cython path but still backed by the native scanner.
      \item \textbf{Pure-Python regular-expression tokeniser.} A single master regular expression provides an always-available fallback that requires no build step, guaranteeing the system runs unmodified on any Python installation.
    \end{enumerate}

    \nd For large files, the compiled backends scan the source once and then materialise 4-tuples in paren-balanced, whole-form batches, bounding peak memory to a single batch rather than the entire token stream. A \texttt{TokenizerHandler} probes the backends at import time and binds the first available one, so deployments transparently obtain native-scanner performance when the extensions are built and correct results when they are not.
    
    \subsection{Reasoning Services} \label{sec:arch:services}
    
    \nd In \Cref{tab:services}, we summarise the reasoning tasks exposed by the system, each implemented as a subclass of an abstract \texttt{Query} class with \texttt{preprocess} and \texttt{solve} methods.
    
    \input{reasoning_services}
    
    \subsection{Relationship to the Java Implementation} \label{sec:arch:java}
    
    \nd The architecture closely follows the original Java reasoner~\cite{BobilloStraccia2016}, preserving the same preprocessing phases, tableau rules, and MILP encodings. The main structural differences with respect to the Java version are:
    \begin{enumerate}
      \item the \texttt{Fuzzy OWL~2} bridge depends on \texttt{pyowl2} rather than the OWL~API;
      \item the parser is a hand-written recursive-descent module rather than a JavaCC-generated parser; and
      \item solver dispatch supports various MILP solvers, namely CBC, MIP, HiGHS, and CPLEX, in addition to Gurobi, which is the only one supported by \texttt{fuzzyDL}.
    \end{enumerate}

    \section{Usage Examples}\label{sec:applications}

    \nd In this section, a concise illustration of the library's practical application is presented. In particular, we show the process of parsing a \texttt{fuzzyDL}  file or loading a Fuzzy OWL ontology and querying it.
    
    In all the applications provided, we employ the configuration contained in the file \enquote{CONFIG.ini}, shown in \Cref{listing:fdl-config}. In particular, this file specifies the following parameters:
    \begin{itemize}
        \item \texttt{debugPrint}: debugging output is disabled;
        \item \texttt{epsilon}: the precision of the computed solution is fixed; for example, $\epsilon = 0.001$ means that the result is computed with accuracy up to the third decimal place;
        \item \texttt{maxIndividuals}: the maximum number of individuals to be handled is specified; the value $-1$ indicates that no upper bound is imposed;
        \item \texttt{owlAnnotationLabel}: the annotation property used to generate the \texttt{Fuzzy OWL~2} RDF/XML ontology is defined; this value will also be used in the subsequent example;
        \item \texttt{milpProvider}: the MILP provider adopted by the reasoner is selected; in the present case, the Python library \texttt{mip} is used.
    \end{itemize}
    
    \begin{listing}
        \caption{Configuration file for the parsing procedure, stored as \enquote{CONFIG.ini}.}
        \label{listing:fdl-config}
        
\begin{lstlisting}[
            language=XML,
            frame=lines,
            framesep=2mm,
            basicstyle=\ttfamily\scriptsize\linespread{1.2}\selectfont,
            numbers=left
        ]
[DEFAULT]
debugPrint = False
epsilon = 0.001
maxIndividuals = -1
owlAnnotationLabel = fuzzyLabel
milpProvider = mip
\end{lstlisting}
        
    \end{listing}

\subsection{Converting a Fuzzy OWL~2 Ontology into a fuzzyDL KB} \label{sec:applications-example-3}
    \nd In this example, we illustrate how to parse a \texttt{Fuzzy OWL~2} ontology in order to obtain the corresponding \texttt{fuzzyDL}  KB file. The resulting directory structure is reported in \Cref{listing:fdl-dir-structure-ex3}.
    
    \nd In this case, only the script shown in \Cref{listing:fdl-script-owl2fdl} is required, together with the configuration file reported in \Cref{listing:fdl-config}. More precisely, the ontology produced in \Cref{sec:applications-example-2}, stored in the file \enquote{example.owl}, is translated back into the corresponding \texttt{fuzzyDL}  KB \enquote{example.fdl}, which is saved in the subdirectory \texttt{./results}. The resulting KB is presented in \Cref{listing:owl2fdl-results}. It is worth observing that the generated file also contains the \texttt{fuzzyDL}  declaration
    \begin{equation*}
        \kw{(define-primitive-concept Car *top*)}
    \end{equation*}
    which states that the class \texttt{Car} is a subclass of the OWL class \texttt{Thing}.
    
    \begin{listing}
        \caption{Directory structure adopted in the example presented in \Cref{sec:applications-example-3}.}
        \label{listing:fdl-dir-structure-ex3}
        
\begin{lstlisting}[
            language=XML,
            frame=lines,
            framesep=2mm,
            basicstyle=\ttfamily\scriptsize\linespread{1.2}\selectfont,
            numbers=left
        ]
-- main/
---- CONFIG.ini
---- script.py
---- example.owl
---- results/
-------- example.fdl
---- logs/
-------- reasoner/
------------ CURRENT_YEAR/
---------------- CURRENT_MONTH/
-------------------- CURRENT_DAY/
------------------------ fuzzydl_HH-MM-SS.log
\end{lstlisting}
        
    \end{listing}
    
    \begin{listing}
        \caption{Python script used to parse the OWL~2 ontology file \enquote{example.owl} and translate it into the corresponding \texttt{fuzzyDL}  KB, stored as \enquote{example.fdl} in the subdirectory \texttt{./results}.}
        \label{listing:fdl-script-owl2fdl}
        
\begin{lstlisting}[
            language=XML,
            frame=lines,
            framesep=2mm,
            basicstyle=\ttfamily\scriptsize\linespread{1.2}\selectfont,
            numbers=left
        ]
from fuzzy_dl_owl2.fuzzyowl2.fuzzyowl2_to_fuzzydl import FuzzyOwl2ToFuzzyDL

# The FuzzyOwl2ToFuzzyDL class is used to translate an OWL 2 ontology 
# into a FDL knowledge base. 
# The constructor takes two arguments: 
#   - the path to the input OWL file and 
#   - the path to the output FDL file. 
fdl = FuzzyOwl2ToFuzzyDL("./results/example.owl", "example.fdl")
# save example.fdl in the subdirectory "./results"
fdl.translate_owl2ontology()  
\end{lstlisting}

    \end{listing}
    
    \begin{listing}
        \caption{\texttt{fuzzyDL}  KB  \enquote{example.fdl} obtained from the Fuzzy OWL ontology \enquote{example.owl} by means of \Cref{listing:fdl-script-owl2fdl}.}
        \label{listing:owl2fdl-results}
        
\begin{lstlisting}[
            language=XML,
            frame=lines,
            framesep=2mm,
            basicstyle=\ttfamily\scriptsize\linespread{1.2}\selectfont,
            numbers=left
        ]
(define-fuzzy-logic lukasiewicz)

(define-modifier very linear-modifier(0.8))

(define-fuzzy-concept eq243 crisp(0.0, 400.0, 243.0, 243.0))
(define-fuzzy-concept geq300 crisp(0.0, 400.0, 300.0, 400.0))

(define-fuzzy-concept High right-shoulder(0.0, 400.0, 180.0, 250.0))

(define-fuzzy-concept VeryHigh modified(very, High))

(define-primitive-concept Car *top*)

(define-concept SportCar (and Car (some speed VeryHigh)))

(instance audi (and Car (some speed eq243)) 1.0)
(instance ferrari (and Car (some speed geq300)) 1.0)
\end{lstlisting}
        
    \end{listing}

    \subsection{Querying a fuzzyDL  KB}\label{sec:applications-example-1}

    \nd In this section, we present an example illustrating how to query a \texttt{fuzzyDL}  KB. The final directory structure adopted in this example is reported in \Cref{listing:fdl-dir-structure}.
    
    More specifically, we load the KB stored in the file \enquote{example.fdl} (see \Cref{listing:fdl-kb}), which describes a \texttt{SportCar} as a \texttt{Car} whose speed is \texttt{very} \texttt{High}. Furthermore, in the final line of the listing, the KB is queried in order to determine the membership degrees of the individuals \texttt{audi} and \texttt{ferrari} with respect to the concept \texttt{SportCar}.
    
    As a preliminary step, a configuration file must be specified in order to initialise the MILP environment constants used by the library. 
    
    \nd For this basic example, it is sufficient to use the Python script reported in \Cref{listing:fdl-script} in order to load the file \enquote{example.fdl}. The output is written to a log file generated within a new directory named \texttt{logs}. More precisely, the log file is stored at the path 
    \begin{center}
        \texttt{logs/reasoner/CURRENT\_YEAR/CURRENT\_MONTH/CURRENT\_DAY}
    \end{center} 
    \nd and is named \mbox{\enquote{fuzzydl\_HH-MM-SS.log}}, where \texttt{HH}, \texttt{MM}, and \texttt{SS} denote the starting time of the script execution. An example is provided by the log file \enquote{fuzzydl\_15-19-48.log} shown in \Cref{listing:fdl-log}. Accordingly, the answer to the first query is $0.92$, meaning that \texttt{audi} belongs to the fuzzy set \texttt{SportCar} with membership degree $0.92$. By contrast, the answer to the second query is $1.0$, that is, \texttt{ferrari} belongs to \texttt{SportCar} with membership degree $1.0$.
    
    The KB may also be parsed in an alternative manner, so that the result is displayed directly on screen rather than being inspected through the log file. For instance, \Cref{listing:fdl-script-alt} reports an alternative script achieving this behaviour. In that case, the parser loads the KB and executes the query without writing the result to the log output. Instead, it returns an instance of the class \texttt{Solution}, which stores the query result and can be displayed in a more readable form. The \texttt{Solution} class provides the method \texttt{is\_consistent\_kb}, which checks whether the KB is consistent, and implements the magic method \texttt{\_\_str\_\_}, which returns a string representation of the solution including the membership degree of the individual with respect to the queried concept.

    \begin{listing}
        \caption{Directory structure adopted in the example presented in \Cref{sec:applications-example-1}.}
        \label{listing:fdl-dir-structure}
        
\begin{lstlisting}[
            language=XML,
            frame=lines,
            framesep=2mm,
            basicstyle=\ttfamily\scriptsize\linespread{1.2}\selectfont,
            numbers=left
        ]
-- main/
---- CONFIG.ini
---- script.py
---- example.fdl
---- logs/
-------- reasoner/
------------ CURRENT_YEAR/
---------------- CURRENT_MONTH/
-------------------- CURRENT_DAY/
------------------------ fuzzydl_HH-MM-SS.log
\end{lstlisting}

    \end{listing}
    
    \begin{listing}
        \caption{\texttt{fuzzyDL}  KB stored in the file \enquote{example.fdl}.}
        \label{listing:fdl-kb}
\begin{lstlisting}[
            language=XML,
            frame=lines,
            framesep=2mm,
            basicstyle=\ttfamily\scriptsize\linespread{1.2}\selectfont,
            numbers=left
        ]
# Defining the logic to be used in the knowledge base.
# In this case, we are using the Lukasiewicz fuzzy logic.
(define-fuzzy-logic lukasiewicz)

# Defining a linear modifier "very" with a parameter of 0.8, 
# which will be used to modify the membership degree of concepts.
(define-modifier very linear-modifier(0.8))

# Defining crisp concepts for speed with specific values
(define-fuzzy-concept eq243 crisp(0, 400, 243, 243))
(define-fuzzy-concept geq300 crisp(0, 400, 300, 400))
# Defining a fuzzy concept for speed with a right shoulder function, 
# where the membership degree starts to increase from 180, 
# reaches 1 at 250, and remains 1 afterwards.
(define-fuzzy-concept High right-shoulder(0, 400, 180, 250))
# Defining a modified concept "VeryHigh" by applying "very" to "High".
(define-fuzzy-concept VeryHigh modified(very, High))

# Defining the complex concept "SportCar" as a Car with a speed that
# is "VeryHigh".
(define-concept SportCar (and Car (some speed VeryHigh)))

# Defining two instances of Car with different speeds.
(instance ferrari (and Car (some speed geq300)) 1)
(instance audi (and Car (some speed eq243)) 1)

# What is the membership degree of "audi" with respect to the 
# concept "SportCar"?
(min-instance? audi SportCar)
# What is the membership degree of "ferrari" with respect to the 
# concept "SportCar"?
(min-instance? ferrari SportCar)
\end{lstlisting}
    \end{listing}

    \begin{listing}
        \caption{Python script used to parse and query the KB stored in the file \enquote{example.fdl}.}
        \label{listing:fdl-script}
        
\begin{lstlisting}[
            language=XML,
            frame=lines,
            framesep=2mm,
            basicstyle=\ttfamily\scriptsize\linespread{1.2}\selectfont,
            numbers=left
        ]
from fuzzy_dl_owl2.fuzzydl.parser import DLParser

DLParser.main("./example.fdl")
# Is audi instance of SportCar ? >= 0.92
# Is ferrari instance of SportCar ? >= 1.0
\end{lstlisting}
        
    \end{listing}

    \begin{listing}
        \caption{Alternative Python script used to parse and query the KB stored in the file \enquote{example.fdl}.}
        \label{listing:fdl-script-alt}
        
\begin{lstlisting}[
            language=XML,
            frame=lines,
            framesep=2mm,
            basicstyle=\ttfamily\scriptsize\linespread{1.2}\selectfont,
            numbers=left
        ]
from fuzzy_dl_owl2.fuzzydl.parser import DLParser
from fuzzy_dl_owl2.fuzzydl.milp.solution import Solution

# parse the knowledge base and get the queries
kb, queries = DLParser.get_kb("./example.fdl") 
# solve the knowledge base to check for consistency and prepare 
# for answering queries
kb.solve_kb() 

# Execute the queries and print the results
for query in queries:
    # Execute the query and get the result as a Solution object
    result: Solution = query.solve(kb)
    # Check if the knowledge base is consistent before printing the result
    if not result.is_consistent_kb():
        break
    # Print the query and the result in a readable format
    print(f"{query}{result}")
    # Is audi instance of SportCar ? >= 0.92
    # Is ferrari instance of SportCar ? >= 1.0
\end{lstlisting}
        
    \end{listing}
    
    \begin{listing}
        \caption{Example log file produced by the parsing and query execution performed with the script in \Cref{listing:fdl-script}.}
        \label{listing:fdl-log}
        
\begin{lstlisting}[
            language=XML,
            frame=lines,
            framesep=2mm,
            basicstyle=\ttfamily\scriptsize\linespread{1.2}\selectfont,
            numbers=left
        ]
DATETIME - INFO -- Knowledge Base parsed in 0.015450833s
DATETIME - INFO -- Using Python-MIP package version 1.16rc0
DATETIME - INFO -- Is audi instance of SportCar? >= 0.92
DATETIME - INFO -- Time (s): 0.171397584
DATETIME - INFO -- Is ferrari instance of SportCar? >= 1.0
DATETIME - INFO -- Time (s): 0.004942084
\end{lstlisting}

    \end{listing}

    \subsection{Converting a fuzzyDL  KB into a Fuzzy OWL~2 Ontology}\label{sec:applications-example-2}
    \nd In this example, we illustrate the conversion of a KB specified in a \texttt{fuzzyDL} file into a \texttt{Fuzzy OWL~2} ontology. To this end, we use the same configuration and KB files as those reported in \Cref{listing:fdl-config,listing:fdl-kb}, respectively. The resulting directory structure for this example is shown in \Cref{listing:fdl-dir-structure-ex2}.
    
\begin{figure}[!htb]
        \centering
        \caption{\texttt{Fuzzy OWL~2} annotation \texttt{fuzzyLabel} and data property \texttt{speed}, obtained from \enquote{example.fdl} by means of \Cref{listing:fdl-script-fdl2owl}.}
        \label{fig:owl-speed-fuzzylabel}
        \begin{minipage}[tbh!]{0.45\textwidth}
            \vspace{0pt}
            \small \textbf{(a)} \texttt{Fuzzy OWL~2} annotation \texttt{fuzzyLabel} produced by \Cref{listing:fdl-script-fdl2owl}.
            \label{fig:owl-fuzzy-label}
            \includegraphics[width=\linewidth]{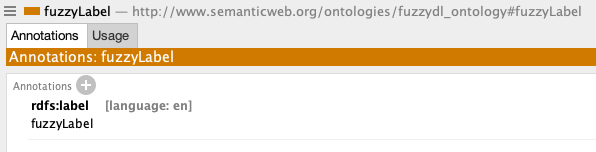}
        \end{minipage}
        \hfill
        \begin{minipage}[tbh!]{0.45\textwidth}
            \vspace{0pt}
            \small \textbf{(b)} OWL~2 data property \texttt{speed} derived from \enquote{example.fdl} by means of \Cref{listing:fdl-script-fdl2owl}.
            \label{fig:owl-speed}
            \includegraphics[width=\linewidth]{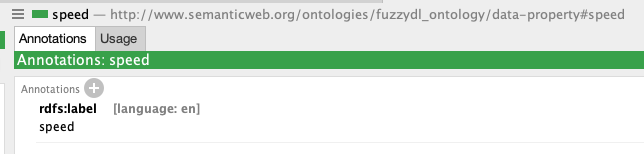}
        \end{minipage}
    \end{figure}

\begin{figure}[!htb]
        \caption{OWL~2 classes, including the class \texttt{SportCar}, obtained from \enquote{example.fdl} by means of \Cref{listing:fdl-script-fdl2owl}.}
        \label{fig:owl-classes-sportcar}
        \centering
        \begin{minipage}[tbh!]{0.45\textwidth}
            \vspace{0pt}
            \small \textbf{(a)} OWL~2 classes derived from \enquote{example.fdl} by means of \Cref{listing:fdl-script-fdl2owl}.
            \label{fig:owl-classes}
            \includegraphics[width=\linewidth]{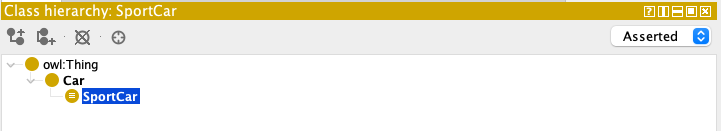}
        \end{minipage}
        \hfill
        \begin{minipage}[tbh!]{0.45\textwidth}
            \vspace{0pt}
            \small \textbf{(b)} OWL~2 class \texttt{SportCar} derived from \enquote{example.fdl} by means of \Cref{listing:fdl-script-fdl2owl}.
            \label{fig:owl-class-sportcar}
            \includegraphics[width=\linewidth]{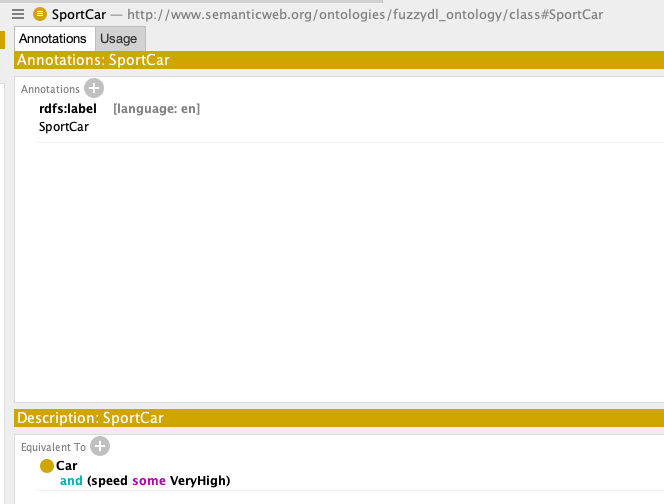}
        \end{minipage}
    \end{figure}

\begin{figure}[!htb]
        \centering
        \caption{OWL~2 individuals \texttt{audi} and \texttt{ferrari}, obtained from \enquote{example.fdl} by means of \Cref{listing:fdl-script-fdl2owl}.}
        \label{fig:owl-individuals}
        \begin{minipage}[tbh!]{0.45\textwidth}
            \vspace{0pt}
            \small \textbf{(a)} OWL~2 individual \texttt{audi} derived from \enquote{example.fdl} by means of \Cref{listing:fdl-script-fdl2owl}.
            \label{fig:owl-individual-audi}
            \includegraphics[width=\linewidth]{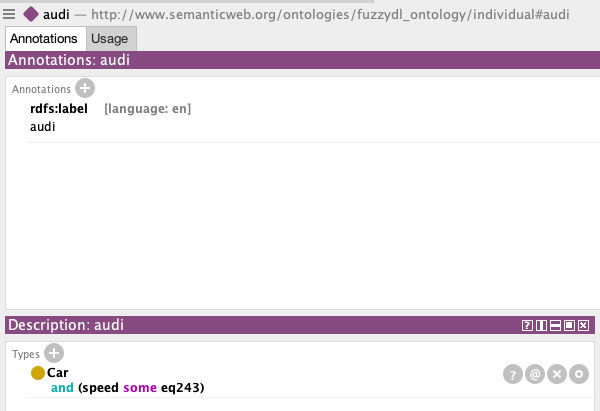}
        \end{minipage}
        \hfill
        \begin{minipage}[tbh!]{0.45\textwidth}
            \vspace{0pt}
            \small \textbf{(b)} OWL~2 individual \texttt{ferrari} derived from \enquote{example.fdl} by means of \Cref{listing:fdl-script-fdl2owl}.
            \label{fig:owl-individual-ferrari}
            \includegraphics[width=\linewidth]{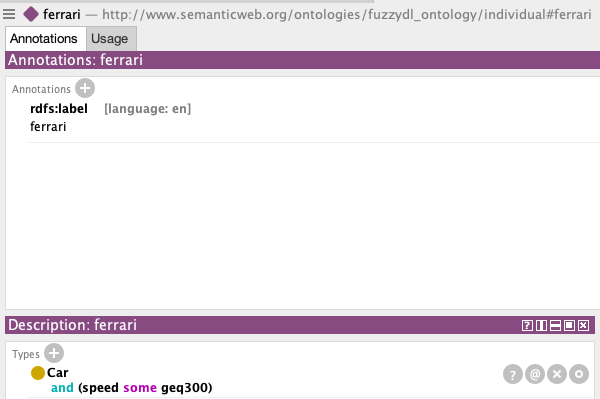}
        \end{minipage}
    \end{figure}

\begin{figure}[!htb]
        \centering
        \caption{\texttt{Fuzzy OWL~2} datatypes \texttt{very}, \texttt{High}, and \texttt{VeryHigh}, obtained from \enquote{example.fdl} by means of \Cref{listing:fdl-script-fdl2owl}.}
        \label{fig:owl-fuzzy-datatypes}
        \begin{minipage}[tbh!]{0.3\textwidth}
            \vspace{0pt}
            \small \textbf{(a)} \texttt{Fuzzy OWL~2} datatype \texttt{very} produced by \Cref{listing:fdl-script-fdl2owl}.
            \label{fig:owl-very}
            \includegraphics[width=\linewidth]{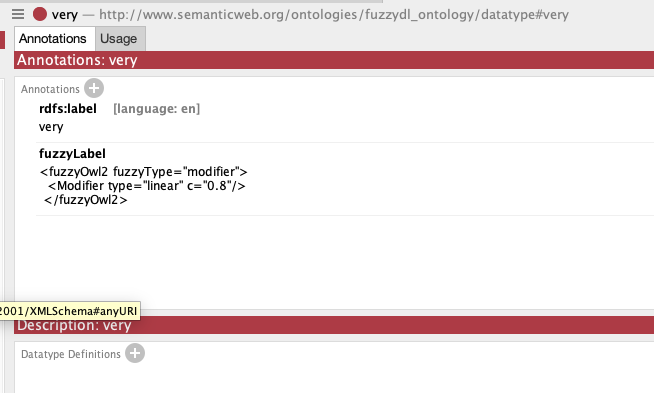}
        \end{minipage}
        \hfill
        \begin{minipage}[tbh!]{0.3\textwidth}
            \vspace{0pt}
            \small \textbf{(b)} \texttt{Fuzzy OWL~2} datatype \texttt{High} derived from \enquote{example.fdl} by means of \Cref{listing:fdl-script-fdl2owl}.
            \label{fig:owl-high}
            \includegraphics[width=\linewidth]{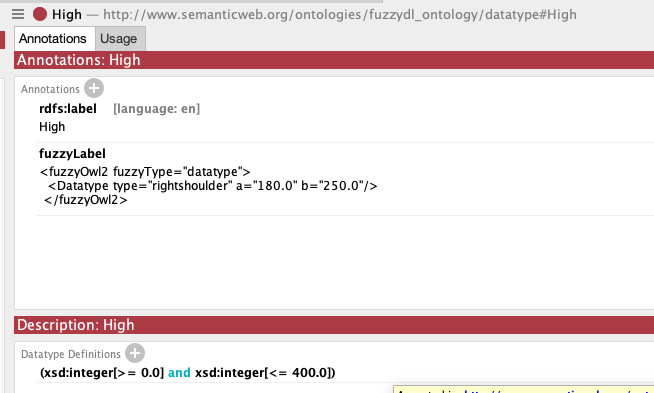}
        \end{minipage}
        \hfill
        \begin{minipage}[tbh!]{0.3\textwidth}
            \vspace{0pt}
            \small \textbf{(c)} \texttt{Fuzzy OWL~2} datatype \texttt{VeryHigh} derived from \enquote{example.fdl} by means of \Cref{listing:fdl-script-fdl2owl}.
            \label{fig:owl-very-high}
            \includegraphics[width=\linewidth]{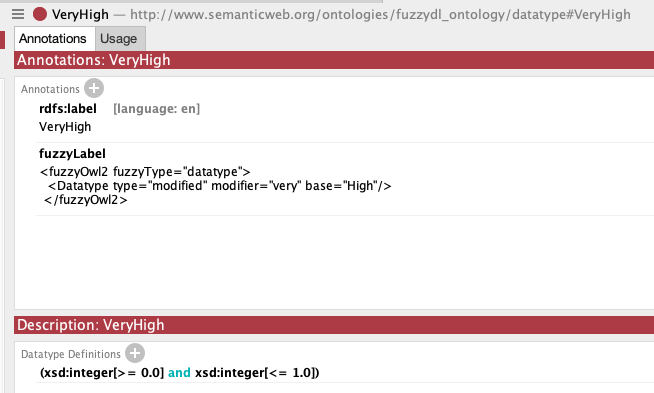}
        \end{minipage}
    \end{figure}

\begin{figure}[!htb]
        \caption{OWL~2 crisp datatypes \texttt{eq243} and \texttt{geq300}, obtained from \enquote{example.fdl} by means of \Cref{listing:fdl-script-fdl2owl}.}
        \label{fig:owl-crisp-datatypes}
        \centering
        \begin{minipage}[tbh!]{0.45\textwidth}
            \vspace{0pt}
            \small \textbf{(a)} OWL~2 crisp datatype \texttt{eq243} produced by \Cref{listing:fdl-script-fdl2owl}.
            \label{fig:owl-eq243}
            \includegraphics[width=\linewidth]{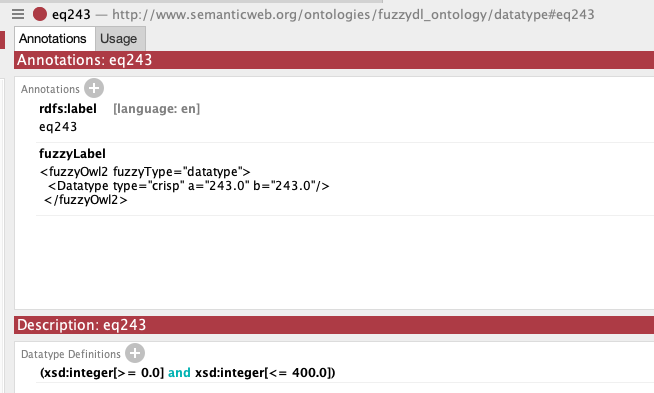}
        \end{minipage}
        \hfill
        \begin{minipage}[tbh!]{0.45\textwidth}
            \vspace{0pt}
            \small \textbf{(b)} OWL~2 crisp datatype \texttt{geq300} derived from \enquote{example.fdl} by means of \Cref{listing:fdl-script-fdl2owl}.
            \label{fig:owl-geq300}
            \includegraphics[width=\linewidth]{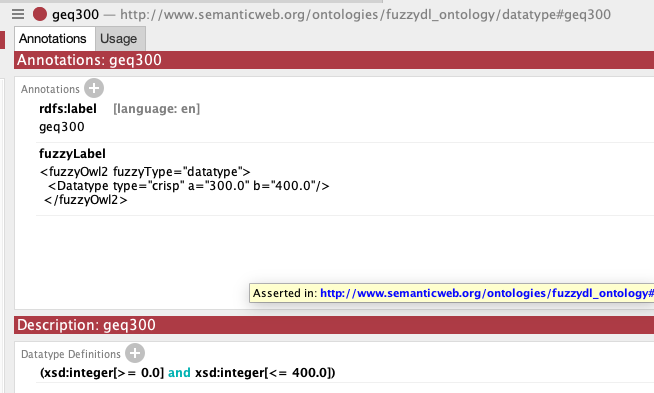}
        \end{minipage}
    \end{figure}

    \nd In this setting, the script reported in \Cref{listing:fdl-script-fdl2owl} is executed in order to translate the \texttt{fuzzyDL}  KB given in \Cref{listing:fdl-kb} into the corresponding OWL~2 ontology. The generated OWL output is graphically presented in \Cref{fig:owl-speed-fuzzylabel,fig:owl-classes-sportcar,fig:owl-individuals,fig:owl-fuzzy-datatypes,fig:owl-crisp-datatypes}.
    
    First, the library generates the OWL annotation property \texttt{fuzzyLabel}, which is used to annotate fuzzy entities, including classes and datatypes. A graphical representation of this annotation is provided in \Cref{fig:owl-speed-fuzzylabel}(b).
    
    Moreover, \Cref{fig:owl-speed-fuzzylabel}(a) displays the generated OWL data property \texttt{speed}, whereas \Cref{fig:owl-classes-sportcar} presents the generated OWL classes, with particular emphasis on the class \texttt{SportCar}.
    
    In \Cref{fig:owl-individuals}, we report the generated individuals \texttt{audi} and \texttt{ferrari}. Finally, \Cref{fig:owl-fuzzy-datatypes,fig:owl-crisp-datatypes} show the generated fuzzy datatypes (\texttt{very}, \texttt{High}, and \texttt{VeryHigh}) and crisp datatypes (\texttt{eq243} and \texttt{geq300}), respectively. We note that the modified datatype \texttt{VeryHigh} is represented directly through the XML fragment reported in \Cref{listing:owl-very-high-xml}.
  
    \begin{listing}
        \caption{Directory structure adopted in the example presented in \Cref{sec:applications-example-2}.}
        \label{listing:fdl-dir-structure-ex2}
        
\begin{lstlisting}[
            language=XML,
            frame=lines,
            framesep=2mm,
            basicstyle=\ttfamily\scriptsize\linespread{1.2}\selectfont,
            numbers=left
        ]
-- main/
---- CONFIG.ini
---- script.py
---- example.fdl
---- results/
-------- example.owl
---- logs/
-------- reasoner/
------------ CURRENT_YEAR/
---------------- CURRENT_MONTH/
-------------------- CURRENT_DAY/
------------------------ fuzzydl_HH-MM-SS.log
\end{lstlisting}

    \end{listing}
    
    \begin{listing}
        \caption{Python script used to parse the \texttt{fuzzyDL}  KB in \enquote{example.fdl} and translate it into the corresponding OWL~2 ontology, stored as \enquote{example.owl} in the subdirectory \texttt{./results}.}
        \label{listing:fdl-script-fdl2owl}
		
\begin{lstlisting}[
            language=XML,
            frame=lines,
            framesep=2mm,
            basicstyle=\ttfamily\scriptsize\linespread{1.2}\selectfont,
            numbers=left
        ]

from fuzzy_dl_owl2.fuzzydl.fuzzydl_to_owl2 import FuzzydlToOwl2

# The FuzzydlToOwl2 class is used to translate an fuzzyDL KB 
# into an OWL 2 ontology. 
# The constructor takes two arguments: 
#   - the path to the input fuzzyDL file and 
#   - the path to the output OWL file. 
fdl = FuzzydlToOwl2("./example.fdl", "example.owl")
# save "example.owl" in the subdirectory "./results"
fdl.run()  
\end{lstlisting}

    \end{listing}
    
    \begin{listing}
        \caption{XML annotation used to represent the modified fuzzy datatype \texttt{VeryHigh} from the \texttt{fuzzyDL}  file \enquote{example.fdl} in the corresponding OWL~2 ontology \enquote{example.owl}.}
        \label{listing:owl-very-high-xml}
        
\begin{lstlisting}[
            language=XML,
            frame=lines,
            framesep=2mm,
            basicstyle=\ttfamily\scriptsize\linespread{1.2}\selectfont,
            numbers=left
        ]
<fuzzyOwl2 fuzzyType="datatype">
    <Datatype type="modified" modifier="very" base="High" />
</fuzzyOwl2>
\end{lstlisting}
        
    \end{listing}
    
    \section{Benchmarks}\label{sec:benchmarks}
    \nd In this section, we present a comparison between the Java and Python codebases through a series of benchmark experiments. All benchmarks were executed on an Apple Mac Mini M4 equipped with $24$ GB of RAM and a $500$ GB SSD. For all benchmarks, we used the optimised Cythonized parser \texttt{DLParserFast} with \texttt{re2c} tokeniser. Moreover, every benchmark was executed considering \L{}ukasiewicz logic. 
    
    The ontologies in our dataset are detailed in Table~\ref{tab:ontos}, which have been taken from~\cite{Bobillo16a}. For each ontology, we include its name (column Ontology) and some relevant features, namely the number of classes, object properties, data properties, number of individuals, and logical axioms. The last column indicates whether there were some missing ontology imports (a referenced ontology to be imported has not been found). 
    
    

   \input{ontoinfo}

    The benchmark results obtained for each KB and for each codebase are represented in \Cref{tab:java-python-benchmarks}. 
    In these experiments, we evaluated the query \kw{(sat?)}, i.e., check satisfiability, over the KBs, excluding those whose execution time exceeded $20$ minutes for both implementations. All runs were made under 
    \L{}ukasiewicz semantics. For the Python codebase, the benchmarks were performed using the following three MILP solvers: Gurobi, CBC, and HiGHS. Each KB was benchmarked ten times to mitigate the effects of CPU idle fluctuations, since the processor was not fully dedicated to the benchmark computations. For each benchmark, both the average execution time and the unbiased standard deviation across the runs were computed. Whenever a KB exceeded the imposed timeout limit and raised a timeout exception, the remaining runs for that KB were skipped. The benchmark measured the overall wall-clock time required by both codebases to parse each KB and solve the query \kw{(sat?)}. 

    \nd In our benchmark comparisons between the Python and Java implementations of the fuzzyDL reasoner, we observe performance characteristics primarily attributable to language-specific runtime overhead rather than to algorithmic differences. Since Python is an interpreted language, it incurs an additional runtime overhead associated with interpreter startup and execution. By contrast, Java executes precompiled bytecode on the Java Virtual Machine, thereby avoiding much of this initial overhead. As a consequence, for smaller KBs, the Java implementation tends to exhibit slightly lower execution times than the Python implementation. However, the Python implementation demonstrates superior performance on computationally intensive tasks while exhibiting overhead on lightweight operations due to interpreter startup costs.

    For small KBs, the Java implementation generally outperforms the Python implementation because of its lower startup overhead. For example:
    \begin{itemize}
        \item \texttt{atom-common}: Java $0.12$\,s versus Python/Gurobi $0.6$\,s (Python approximately $5\times$ slower);
        \item \texttt{SIGKDD-EKAW}: Java $0.11$\,s versus Python/Gurobi $0.63$\,s (Python approximately $5.7\times$ slower).
    \end{itemize}
    
    \nd For large KBs requiring substantial reasoning, the Python implementation outperforms the Java implementation owing to its more efficient MILP generation and tighter solver integration. For instance,
    \begin{itemize}
        \item \texttt{FMA}: Java $27.45$\,s versus Python/Gurobi $3.38$\,s (Python approximately $8.13\times$ faster);
        \item \texttt{FuzzyWine}: Java $62.87$\,s versus Python/Gurobi $14.97$\,s (Python approximately $4.199\times$ faster);
        \item \texttt{chebi}: no Java result available versus Python/Gurobi $9.95$\,s (the Java implementation failed to solve the Knowledge Base within the imposed $20$-minute timeout).
    \end{itemize}
    
    \nd The performance crossover point occurs at approximately $1$--$2$ seconds of solver time, where Python's startup overhead (roughly $0.5$s) becomes negligible with respect to the overall computation time. As pointed out before, this behaviour is typical of Python architectures, in which fixed startup costs are progressively amortised as the computational workload increases. In production settings, using persistent Python processes would largely eliminate this discrepancy while preserving the algorithmic advantages of the Python implementation.

    An interesting case is represented by the KB \path{mosquito_insecticide_resistance.obo}. The Java implementation solved this benchmark in $3.29$ seconds, whereas the fastest Python configuration required $12.21$ seconds, making it approximately four times slower. Our investigation of this anomaly led to the same conclusion as for the smaller KBs: the performance gap is largely attributable to the overhead associated with runtime object creation and initialisation in the Python implementation.

    A further interesting case is represented by the KB \path{propreo.TBox}, for which the HiGHS solver outperformed Gurobi. This observation suggests that supporting multiple MILP solvers is valuable not only for providing access to open-source alternatives but also because a given KB may be solved more efficiently by a solver other than the highly optimised Gurobi. In other words, solver performance appears to depend not only on the overall quality of the optimisation engine, but also on the structural characteristics of the specific MILP instance generated from the KB.

    Finally, we note that four KBs (\path{propreo}, \path{organic-compound-complex}, \path{AirSystem}, and \path{chebi}) could not be solved by the Java implementation within the imposed $20$-minute timeout, whereas the Python implementation completed all four benchmarks in less than one minute.

    \input{java_python_benchamarks}
       
    \section{Conclusions}\label{sec:conclusions}


    \nd While classical ontology languages struggle to handle vague or imprecise knowledge, fuzzy ontologies enrich them with fuzzy logic to model partial truths. Our contribution consists of a complete Python re-engineering of the \texttt{fuzzyDL}  reasoner and the \texttt{Fuzzy OWL~2} framework~\cite{BobilloStraccia2016,BOBILLO20111073}.
    Specifically, we have released \texttt{fuzzy-dl-owl2}, porting the original Java-based software into Python to address several shortcomings, including semantic inconsistencies, architectural rigidity, and limited solver options. Key Improvements in \texttt{fuzzy-dl-owl2} are the following:
    
    \begin{description}    
        \item[Architectural Redesign of \texttt{fuzzyDL} :] The original \texttt{Concept} class was overly complex, containing attributes that were not always necessary (e.g., assigning a scalar weight to every concept even if unused). Our Python port simplifies this into a modular class hierarchy where Concept only holds basic attributes, relying on specific sub-classes for complex logic.
        \item[Decoupling:] To avoid circular imports—a major issue in Python—the developers decoupled various classes from the \texttt{KnowledgeBase} class, centralising its management and spinning off specific methods into new handler classes;
        \item[Expanded MILP Solver Support:] The original implementation relied only on the Gurobi solver for MILP. The new framework supports a broader range of solvers, including CPLEX and open-source alternatives such as GLPK, CBC, and HiGHS, and uses Python libraries such as gurobipy, MIP, and PuLP.
        \item[Semantic Fixes in \texttt{Fuzzy OWL~2}:] The original software struggled with IRI ambiguities, such as generating identical IRIs if an individual and a class shared the same name. The new implementation assigns specific sub-paths 
        (e.g., individual/, class/) to ensure distinct IRIs. It also fixes redundancy issues where data properties were incorrectly converted into both object and data properties.
        \item[Improved XML Parsing and OWL~2 Compliance:] We have developed a library, called \texttt{pyowl2} \cite{pyowl2}, with the aim of supporting complex annotations, which are part of the OWL~2 standard and are a fundamental ingredient of \texttt{Fuzzy OWL~2}. Additionally, it replaces the original custom, error-prone XML parser with Python's standard XML library, which safely translates string XML annotations into more easily machine-readable structured XML objects.
    \end{description}
    
    \nd In summary, the resulting Python framework provides an extensible, portable, and theoretically grounded platform for reasoning with fuzzy ontologies. The source code and documentation are publicly available on \texttt{GitHub} and \texttt{Read the Docs} to encourage community adoption.
    
    For future work, we plan to employ this library in ontological applications along the lines discussed in \cite{CastronovoFilipponeGaliciLaRosaPavoneTabacchi:DSA-ISC:2025,FilipponeLaRosaTabacchi:SUM:2024}. Furthermore, we look forward to porting also other reasoners such as \texttt{DeLorean} and implementing a reasoning algorithm dedicated to some \texttt{Fuzzy OWL~2} Profiles,~\cite{owl2_profiles_w3c} such as \texttt{Fuzzy OWL~2} EL~\cite{Bobillo18}. 

    \section*{CRediT authorship contribution statement}
    
    \begin{description}
        \item[Fernando Bobillo:] Supervision, Formal analysis, Software, Writing - Review \& Editing.
        \item[Giuseppe Filippone:] Conceptualization, Methodology, Formal analysis, Software, Investigation, Visualization, Writing - Original Draft, Writing - Review \& Editing. 
        \item[Gianmarco La Rosa:] Methodology, Formal analysis, Visualization, Writing - Original Draft, Writing - Review \& Editing.
        \item[Umberto Straccia:] Supervision, Formal analysis, Software, Writing - Review \& Editing. 
        \item[Marco Elio Tabacchi:] Supervision, Project administration, Funding acquisition, Writing - Review \& Editing.
    \end{description}
 
    \section*{Funding}
    \nd Giuseppe Filippone, Gianmarco La Rosa, and Marco Elio Tabacchi acknowledge financial support from the \textbf{Sustainability Decision Framework} (\textbf{SDF}) Research Project --  CUP \textbf{B79J23000540005} -- Grant Assignment Decree No. \textbf{5486} adopted on \textbf{2023-08-04} by the Italian Ministry of University and Research (MUR). 
    Gianmarco La Rosa also acknowledges support from the \textbf{Gruppo Nazionale per le Strutture Algebriche, Geometriche e le loro Applicazioni} (\textbf{GNSAGA}) of the \textbf{Istituto Nazionale di Alta Matematica} (\textbf{INdAM}) \enquote{Francesco Severi}.
    Marco Elio Tabacchi also acknowledges financial support under the National Recovery and Resilience Plan (NRRP), Mission 4, Component 2, Investment 1.1, Call for tender No. 1409 published on 14.9.2022 by the Italian Ministry of University and Research (MUR), funded by the European Union -- NextGenerationEU -- Project Title \textbf{Quantum Models for Logic, Computation and Natural Processes} (\textbf{QM4NP}) -- CUP \textbf{B53D23030160001} -- Grant Assignment Decree No. \textbf{1371} adopted on \textbf{2023-09-01} by the Italian Ministry of University and Research (MUR).
    Fernando Bobillo and Umberto Straccia acknowledge support from the project \textbf{PID2024-159530OB-I00} funded by MICIU/AEI/10.13039/501100011033/FEDER, UE. 
    Fernando Bobillo also acknowledges support from project \textbf{UZ2024-IyA-02} (funded by the University of Zaragoza).

    \section*{Data availability}
    \nd The Python libraries \texttt{fuzzy-dl-owl2} and \texttt{pyowl2} are freely accessible from \texttt{PyPi} and \texttt{GitHub} repositories. The former is available at \url{https://pypi.org/project/fuzzy-dl-owl2/} and \url{https://github.com/SDF-Unipa/fuzzy_dl_owl2}, while the latter is available at \url{https://pypi.org/project/pyowl2/} and \url{https://github.com/SDF-Unipa/pyowl2}.
    
    \section*{Conflict of Interest}
    \nd The authors declare no conflicts of interest.

    \section*{Acknowledgements}
    \nd We would like to express our gratitude to Vincenzo Reina, a student on the bachelor's degree course in Computer Science at the Department of Mathematics and Computer Science at the Università degli Studi di Palermo, for his contribution to the testing phase on the GitHub repository at the following address: \url{https://github.com/renvins/fuzzy_dl_owl2_testing}. In particular, the student tested the process of translating from FDL code to the respective fuzzy ontology in \texttt{Fuzzy OWL~2}, and vice versa.

    \section*{Abbreviations}
    \nd The following abbreviations are used in this manuscript:	\\

    \begin{tabular}{@{}ll}
        CBC & COIN Branch and Cut\\
        CPLEX & IBM ILOG CPLEX Optimization Studio\\
        DL & Description Logic\\
        FDL & Fuzzy Description Logic\\
        GLPK & GNU Linear Programming Kit\\
        HiGHS & High-performance parallel linear optimization software \\
        IRI & Internationalized Resource Identifier\\
        MILP & Mixed Integer Linear Programming\\
        GCI & General Concept Inclusion\\
        PCD & Primitive Concept Definition\\
        OWL & Web Ontology Language\\
        RDF & Resource Description Framework\\
        UML & Unified Modeling Language\\
        W3C & World Wide Web Consortium\\
        XML & Extensible Markup Language\\
        OWA & Ordered Weighted Average\\
        QOWA & Quantified-Guided Ordered Weighted Average\\
        EBNF & Extended Backus--Naur Form\\
    \end{tabular}

    \appendix

    \section{Grammar and supported syntax}\label{sec:grammar}
    \nd In this section, we present the full grammar handled by our implementation and the syntax of the XML annotations of \texttt{Fuzzy OWL~2}. We refer the reader to \cite{fuzzydlowl2_docs} for more details.
    
    The \texttt{fuzzyDL} reasoner extends classical description logics with fuzzy concepts and roles \cite{BobilloStraccia2016}. The \texttt{fuzzyDL}  concrete syntax is a parenthesised prefix language.  The grammar is reported in Extended Backus--Naur Form (EBNF).  Terminals are written in typewriter font, nonterminals are written as $\nt{nonterminal}$, alternatives are separated by the vertical line $(|)$, optional fragments are denoted as $[\cdot]$, and repeated fragments are enclosed in $\{\cdot\}$.  The symbols $\{\cdot\}^{+}$ and $\{\cdot\}^{2+}$ denote one-or-more and two-or-more repetitions, respectively, whereas $\{\cdot\}^{n}$ denote exactly $n$ repetitions.

    In \Cref{tab:ebnf-lexical-top-level}, we provide the lexical and the top-level grammar of \texttt{fuzzyDL}  in EBNF. The start symbol is given by \nt{knowledge\_base}.
    
    In \Cref{tab:ebnf-datatypes}, we provide the grammar for concrete fuzzy concepts, fuzzy numbers, features, and datatype restrictions in EBNF. In \Cref{tab:ebnf-concepts}, we provide the grammar for the expression of the concepts in EBNF. Lastly, in \Cref{tab:ebnf-axioms-queries}, we provide the grammar for axioms, constraints, statements, and queries in EBNF.
    
    The XML annotation syntax represents \texttt{Fuzzy OWL~2} annotations using a root element \kw{fuzzyOwl2} and a \kw{fuzzyType} attribute.

    In the following parts of this section, we describe in detail the grammars for \texttt{fuzzyDL}  and \texttt{Fuzzy OWL~2}.

    \input{grammar_ebnf_1}
    \input{grammar_ebnf_2}
    \input{grammar_ebnf_3}
    \input{grammar_ebnf_4}

    \subsection{Logic semantics and truth constants}
    \nd At the beginning of a KB, the user specifies the underlying logical semantics. 
    We support three distinct semantic settings, which determine the interpretation of negation, conjunction, disjunction, and implication \cite{fuzzydl_syntax}. In particular, we have 
    \begin{itemize}
        \item \textbf{Zadeh semantics}: this setting employs {\L}ukasiewicz negation, the G\"odel $t$-norm, the G\"odel $t$-conorm, and Kleene--Dienes implication. An exception arises in the case of general concept inclusions, where the degree of membership in the subsumed concept must be less than or equal to the degree of membership in the subsumer concept. This semantics is retained primarily for compatibility with earlier work on fuzzy DLs.
    
        \item \textbf{{\L}ukasiewicz semantics}: this setting adopts {\L}ukasiewicz negation, $t$-norm, $t$-conorm, and implication.
    
        \item \textbf{Classical semantics}: this setting recovers the standard crisp interpretation of conjunction, disjunction, negation, and implication.
    \end{itemize}
    
    \nd In addition, a KB may introduce rational-valued \nt{truth\_constants} by means of declarations of the form \kw{(define-truth-constant $c\ n$)}. These named constants may subsequently be used, for instance, as lower bounds in axioms and queries.

    \subsection{Concrete fuzzy concepts and modifiers}
    
    \nd Fuzzy membership functions may be declared explicitly by means of the \linebreak\kw{define-fuzzy-concept} directive. The language provides several standard families of membership functions, including crisp intervals, left-shoulder and right-shoulder functions, triangular and trapezoidal functions, as well as linear and modified constructions. More precisely, the 
    non-terminal \nt{concept\_type} may take one of the following forms:
    \begin{itemize}
        \item \kw{crisp}$(k_1,k_2,a,b)$, denoting a membership function equal to $1$ on the interval $[a,b]$ and equal to $0$ outside the interval $[k_1,k_2]$;
        
        \item \kw{left-shoulder}$(k_1,k_2,a,b)$, defining a $Z$-shaped membership function over the interval $[a,b]$, with support contained in $[k_1,k_2]$;
        
        \item \kw{right-shoulder}$(k_1,k_2,a,b)$, defining an $S$-shaped membership function over the interval $[a,b]$, with support contained in $[k_1,k_2]$;
        
        \item \kw{triangular}$(k_1,k_2,a,b,c)$ and \kw{trapezoidal}$(k_1,k_2,a,b,c,d)$, corresponding to the usual triangular and trapezoidal membership functions, with support contained in $[k_1,k_2]$;
        
        \item \kw{linear}$(k_1,k_2,a,b)$, representing a linear hedge, with support contained in $[k_1,k_2]$; and
        \item \kw{modified}$(M,F)$, which applies a previously declared modifier $M$ to a datatype restriction $F$.
    \end{itemize}
    
    \nd Modifiers are used to alter the shape of an existing fuzzy concept. The language supports two modifier types, namely \kw{linear-modifier}$(c)$, which defines a scaling hedge with $c>0$, and \kw{triangular-modifier}$(a,b,c)$. A modifier is declared by means of
    \begin{equation*}
        \kw{(define-modifier modifier\_name modifier\_type(parameters))}
    \end{equation*}
    \nd and may subsequently be referenced within a \kw{modified} concept type.

    \subsection{Fuzzy numbers and arithmetic}
    \nd Reasoning with concrete data often requires approximate numerical values. To this end, \texttt{fuzzyDL} supports triangular fuzzy numbers of the form $(a,b,c)$, as well as rational numbers $n$, which are implicitly interpreted as the degenerate triangular fuzzy number $(n,n,n)$. Before fuzzy numbers can be used, the admissible domain range $[k_1,k_2]$ must be declared by means of
    \begin{equation*}
        \kw{(define-fuzzy-number-range $k_1\ k_2$)}.
    \end{equation*}
    \nd If the directive specifying the admissible domain range $[k_1,k_2]$ is omitted, the default range is taken to be $[-\infty,\infty]$.
    
    Subsequently, a named fuzzy number $N$ may be introduced by a declaration of the form
    \begin{equation*}
        \kw{(define-fuzzy-number $N\ \mathrm{expression}$)},
    \end{equation*}
    \nd where \kw{expression} may consist of a constant, a triangular triple, or an arithmetic combination of fuzzy numbers, as specified in \Cref{tab:fuzzy-number-definitions}.

    \input{grammar_fuzzy_numbers}

    \subsection{Features and datatype restrictions}

    \nd \emph{Features} (or functional datatypes) associate individuals with concrete values via data property. A feature $F$ is declared by \kw{(functional $F$)}, and its range is specified by a clause of the form \kw{(range $F$ type parameters)} (see \Cref{tab:feature-definitions}). The admissible types are \kw{*integer* $k_1\,k_2$} and \kw{*real* $k_1\,k_2$}, which denote intervals of integers and reals (viz. rationals), respectively, together with \kw{*string*} and \kw{*boolean*}. 
    
    Datatype restrictions (see \Cref{tab:restriction-definitions}) are then used to constrain individuals on the basis of feature values. Their semantics depend on the choice of an integer $b$, which is required in order to determine the correct solution by means of a MILP solver. The value $b$ must belong to the interval $[k_1,k_2] \subseteq [-k_{\infty}, k_{\infty}]$ associated with the feature $F$. Likewise, the values assigned to \texttt{variable}, $\mathrm{function}(F_1, \ldots, F_n)$, and the range of \texttt{fuzzy\_number} must lie within $[-k_{\infty}, k_{\infty}]$, where $k_{\infty}$ denotes the maximal representable number and depends on the selected MILP solver (see \Cref{tab:milp-kinfty}). We note that the value of $k_{\infty}$ may differ across MILP solvers for numerical reasons. In particular, excessively large values may lead to the accumulation of numerical errors and, consequently, to distorted results. The values currently reported were chosen because they produce consistent results on the test files provided.
    
    An example of a datatype restriction is given by \kw{(>= var $F$)}, which denotes the set of all individuals whose value for $F$ is at least the variable \texttt{var}; similarly, \kw{(<= var $F$)} and \kw{(= var $F$)} express upper-bound and exact-value restrictions, respectively. The variable \texttt{var} may be replaced either by a real number or by a fuzzy-number expression. Function-based restrictions, such as \kw{(>= $F$ (f+ $F_1$ $F_2$))}, compare the value of a feature with a fuzzy number obtained by applying a function; in this example, the function corresponds to the sum of the two fuzzy numbers $F_1$ and $F_2$. Such restrictions are interpreted by means of supremum or infimum operators over the datatype domain, according to the t-norm and implication determined by the chosen semantics. 
    
    Within datatype restrictions, the variable \texttt{var} must be declared beforehand by means of \texttt{(free var)}, in accordance with the constraints introduced in \ref{subsec:crisp-declarations-and-constraints}. Alternatively, \texttt{var} may be replaced directly by a value, namely an integer, a real number, a string, or a Boolean constant (\texttt{true}, \texttt{false}), depending on the range of the feature $F$.

    \input{grammar_feature_definitions}
    \input{grammar_restriction_definitions}
    \input{grammar_milp_kinfty}

    \subsection{Concepts}
    \nd In this section, we describe all concept expressions introduced by the grammar. Their formal mathematical definitions are reported in \Cref{tab:concept-expression-definitions}.

    \input{grammar_concept_definitions}

    \subsubsection{Fuzzy relations and similarity}
    \nd Beyond crisp concepts and roles, one can introduce fuzzy similarity or equivalence relations between individuals. A fuzzy similarity relation $s$ is declared by 
    \begin{equation*}
    \kw{(define-fuzzy-similarity $s$)} \ .    
    \end{equation*}
    \nd Analogously, 
    \begin{equation*}
    \kw{(define-fuzzy-equivalence $s$)}
    \end{equation*}
    \nd introduces a fuzzy equivalence relation. These relations can be employed in rough set approximations. Upper and lower approximations of a concept $C$ with respect to a relation $s$ are constructed using \kw{(ua $s\ C$)} for the upper approximation and \kw{(la $s\ C$)} for the lower approximation; tight and loose variants \kw{tua}, \kw{lua}, \kw{tla}, \kw{lla} refine the approximation and are defined using nested suprema and infima. A local reflexivity construct \kw{(self $S$)} enforces that an individual is related to itself by role $S$~\cite{bobillo12}.
    
    We remark that the fuzzy relation $s$ must be declared in advance either as a fuzzy similarity relation or as a fuzzy equivalence relation, namely by means of \texttt{(define-fuzzy-similarity $s$)} or \texttt{(define-fuzzy-equivalence $s$)}, respectively.

    \subsubsection{Concept expressions}
    \nd Every complex concept is defined recursively from atomic concepts, role restrictions, and datatypes by means of prefix operators. The distinguished symbols \kw{*top*} and \kw{*bottom*} denote, respectively, the universal concept, corresponding to the truth constant $1$, and the inconsistent concept, corresponding to the truth constant $0$. Given an interpretation 
    $\I = (\highi{\Delta},\highi{\cdot})$, an atomic concept $A$ is evaluated as the function $A^{\mathcal{I}}(x)\colon \highi{\Delta} \to [0,1]$. Standard conjunctions and disjunctions are expressed by \kw{(and $C_1\ C_2$)} and \kw{(or $C_1\ C_2$)}, and are interpreted by means of the selected t-norm and t-conorm. The grammar also includes specialised connectives corresponding to specific fuzzy semantics. For instance, \kw{(g-and $C_1\ C_2$)} uses the G\"odel t-norm, whereas \kw{(l-and $C_1\ C_2$)} uses the {\L}ukasiewicz t-norm. Analogous variants are available for disjunction and implication, namely \kw{g-or}, \kw{l-or}, \kw{g-implies}, \kw{l-implies}, \kw{kd-implies}, and \kw{z-implies}. Negation is written as \kw{(not $C$)} and is interpreted according to the negation function associated with the chosen semantics.
    
    Role restrictions make it possible to quantify over related individuals. The existential restriction \kw{(some $R\ C$)} denotes the set of individuals having an $R$-successor in $C$ \cite{straccia2014foundations}; semantically, it is evaluated as the supremum, over all fillers, of the t-norm between the interpretation of the role and that of the concept. Interpretations have to be \emph{witnessed} in the sense that the supremum has to be attained at some individual of the interpretation domain~\cite{HajekP06}. Dually, the universal restriction \kw{(all $R\ C$)} requires all $R$-successors to belong to $C$ and is interpreted through an infimum involving the corresponding fuzzy implication.  Interpretations have to be \emph{witnessed} in the sense that the supremum and infimum have to be attained at some individual of the interpretation domain~\cite{HajekP06}.
    Moreover, the restriction \kw{(some $R$ $a$)} denotes the individual value restriction associated with the role $R$ and the individual $a$, that is, the concept whose interpretation is given by $R^\mathcal{I}(x,a^\mathcal{I})$.
    
    Datatype restrictions may also occur as concept expressions. They are written as \kw{(>= var $F$)}, \kw{(<= var $F$)}, and \kw{(= var $F$)}.    In such expressions, \texttt{var} may also be replaced by a fuzzy number. Their interpretation follows the same principles adopted for datatype restrictions in feature definitions.

    \subsubsection{Special constructs and aggregations}
    \nd The language also provides several higher-level concept constructors.
    Threshold concepts are written as \kw{([>= var] $C$)} and \kw{([<= var] $C$)}.
    These return the value $C^{\mathcal{I}}(x)$ whenever the membership degree of the concept is at least, or at most, the threshold $\mathrm{var}$, respectively, and $0$ otherwise.
    Weighted concepts combine concept membership degrees with numerical coefficients. A basic weighted concept is written as \kw{($n\ C$)}, where the membership degree of $C$ is multiplied by the scalar $n \in \mathbb{R}$.    
    Weighted sums are written as
    \begin{equation*}
        \kw{(w-sum $(n_1\ C_1) \ldots (n_k\ C_k)$)}.
    \end{equation*}
    \nd This constructor returns the sum of the weighted membership degrees of the concepts involved.    
    The weighted maximum operator is written as
    \begin{equation*}
        \kw{(w-max $(v_1\ C_1) \ldots (v_k\ C_k)$)},
    \end{equation*}
    \nd whereas the weighted minimum operator is written as
    \begin{equation*}
        \kw{(w-min $(v_1\ C_1) \ldots (v_k\ C_k)$)}.
    \end{equation*}
    \nd These operators return, respectively, the maximum or minimum obtained by combining the weights with the corresponding concept degrees through the appropriate min--max scheme.    
    Lastly, a zero-constrained weighted sum is written as
    \begin{equation*}
        \kw{(w-sum-zero $(n_1\ C_1) \ldots (n_k\ C_k)$)}.
    \end{equation*}
    \nd This constructor behaves as a weighted sum, except that it returns $0$ whenever at least one component has value $0$.

    For complex decision-making tasks, \texttt{fuzzyDL} provides several aggregation operators~\cite{ASOC2013}.
    The Ordered Weighted Averaging (OWA) operator is written as
    \begin{equation*}
        \kw{(owa $(w_1, \ldots, w_n)$ $(C_1, \ldots, C_n)$)}.
    \end{equation*}
    \nd It combines the concepts $(C_1, \ldots, C_n)$ by means of the weights $(w_1, \ldots, w_n)$, yielding the value $\sum_i w_i y_i$, where $y_i$ denotes the $i$th largest value among $C_i^I(x)$.
    A quantifier-guided OWA is expressed as
    \begin{equation*}
        \kw{(q-owa $Q$ $(C_1, \ldots, C_n)$)}.
    \end{equation*}
    \nd In this case, the weights are induced by a fuzzy quantifier $Q$, usually specified by a right-shoulder or linear function.
    The Choquet integral is written as
    \begin{equation*}
        \kw{(choquet $(v_1, \ldots, v_n)$ $(C_1, \ldots, C_n)$)}.
    \end{equation*}
    \nd Its value is computed as $y_1 v_1 + \sum_{i=2}^n (y_i - y_{i-1}) v_i$.    
    The Sugeno integral is given by
    \begin{equation*}
        \kw{(sugeno $(v_1, \ldots, v_n)$ $(C_1, \ldots, C_n)$)}.
    \end{equation*}
    \nd It returns
    \begin{equation*}
        \max_{i=1}^n \min(y_i, \mu_i),
    \end{equation*}
    \nd where $y_i$ is the $i-th$ largest of the $C_i^I(x)$, and $\mu_i$ is recursively defined by
    \begin{align*}
        \mu_1 &= ow_1,\\
        \mu_i &= ow_i \oplus \mu_{i-1}, \qquad i = 2,\ldots,n,
    \end{align*}
    \nd where $ow_i$ denotes the weight $v_i$ associated with the $i$th largest value among the $C_i^I(x)$.
    The quasi-Sugeno integral is expressed as
    \begin{equation*}
        \kw{(q-sugeno $(v_1, \ldots, v_n)$ $(C_1, \ldots, C_n)$)}.
    \end{equation*}
    \nd It returns
    \begin{equation*}
        \max_{i=1}^n\ (y_i \otimes_{L} \mu_i),
    \end{equation*}
    \nd that is, it employs the {\L}ukasiewicz t-norm in place of the minimum operator used in the Sugeno integral.  
    Finally, the sigma-count construct is written as
    \begin{equation*}
        \kw{(sigma-count $R\ C\ [ a_1\ \ldots\ a_k ] \ F_C$)}.
    \end{equation*}
    \nd This construct counts how many individuals $a_i$ satisfy the concept $C$ and the fuzzy concrete concept $F_C$ along the role $R$. If the list of individuals is empty, all individuals in the ontology are considered.  
    Moreover, the fuzzy concrete concept $F_C$ occurring in a \texttt{sigma-count} concept must be defined beforehand by means of \texttt{define-fuzzy-concept} as a \texttt{left-shoulder}, \texttt{right-shoulder}, or \texttt{triangular} concept.
    
    Important constraints on these operators require all weights $n_i$, $w_i$, and $v_i$ to belong to the interval $[0,1]$, and to satisfy
    \begin{equation*}
        \begin{array}{cc}
        \sum n_i \leq 1, & \sum w_i = 1.
        \end{array}
    \end{equation*}
    \nd Furthermore, Choquet \cite{grabisch2011aggregation} integral is defined with respect to a normalised capacity (or fuzzy measure); accordingly, in the discrete nested-set representation, the relevant normalisation condition is
    \begin{equation*}
        0 \leq v_1 \leq \ldots \leq v_n = 1.
    \end{equation*}
    
    \nd Finally, fuzzy numbers may only occur within existential restrictions, universal restrictions, and datatype restrictions. In threshold concepts, \texttt{var} may be replaced by a value $w \in [0,1]$. 

    \subsection{Crisp declarations and constraints}\label{subsec:crisp-declarations-and-constraints}
    \nd Although \texttt{fuzzyDL} is primarily intended for fuzzy reasoning, it also provides support for crisp concepts and crisp roles. A crisp concept $C$ is declared by means of 
    \begin{equation*}
        \kw{(crisp-concept $C$)} \ ,    
    \end{equation*}
    \nd  whereas a crisp role $R$ is declared by 
    \begin{equation*}
        \kw{(crisp-role $R$)} \ .    
    \end{equation*}
    \nd Such declarations enforce \texttt{classical} semantics, in the sense that the corresponding concept or role are mapped into  
    $\{0,1\}$ rather than into $[0,1]$, for a given individual. 
        
    In addition, the language allows the specification of linear constraints on variables. These are written in the form
    \begin{equation*}
        \kw{(constraint $a_1 \ast \mathrm{var}_1 + \ldots + a_k \ast \mathrm{var}_k\ \mathrm{OP}\ v$)},
    \end{equation*}
    \nd where $\texttt{OP} \in \{\kw{>=}, \kw{<=}, \kw{=}\}$. The language also supports declarations of binary and free variables. In particular, 
    \begin{equation*}
        \kw{(constraint binary var)} 
    \end{equation*}
    \nd restricts the variable \texttt{var} to the set $\{0,1\}$, whereas 
    \begin{equation*}
        \kw{(constraint free var)} 
    \end{equation*}
    \nd allows \texttt{var} to range over all real numbers. These constraints are handled by the underlying MILP solver and are especially useful in optimisation settings.

    \subsection{Axioms}
    \nd A KB is formed by a set of axioms involving assertions on individuals, concept inclusions and definitions, disjointness conditions, role restrictions, and structural role properties. Their formal semantics are given in \Cref{tab:axiom-definitions}. For convenience, we summarise below the admissible axiom forms provided by the grammar. The degree $d$ is optional and, if omitted, is assumed to be equal to $1.0$.

    Before introducing the individual axiom types, two general remarks are in order. First, transitive roles cannot be declared functional. Secondly, if the selected \texttt{fuzzyDL} logic is \texttt{zadeh}, then the implication $(\Longrightarrow)$ is interpreted as Zadeh's set-inclusion operator.
    
    \begin{description}
        \item[Concept assertions] The  \kw{(instance $a\ C\ [d]$)} axiom states that the individual $a$ belongs to the concept $C$ to degree at least $d$.
    
        \item[Role assertions] The  \kw{(related $a\ b\ R\ [d]$)}  axiom states that the pair $(a,b)$ belongs to the role $R$ to degree at least $d$.
    
        \item[Standard general concept inclusions] The  \kw{(implies $C_1\ C_2\ [d]$)} axiom expresses that $C_1$ is subsumed by $C_2$ to degree at least $d$, according to the default implication adopted by the language.
    
        \item[G\"odel general concept inclusions] The  \kw{(g-implies $C_1\ C_2\ [d]$)} axiom specifies subsumption from $C_1$ to $C_2$ to degree at least $d$ by means of G\"odel implication.
    
        \item[Kleene--Dienes general concept inclusions] The  \kw{(kd-implies $C_1\ C_2\ [d]$)} axiom specifies subsumption from $C_1$ to $C_2$ to degree at least $d$ by means of Kleene--Dienes implication.
    
        \item[{\L}ukasiewicz general concept inclusions] The  \kw{(l-implies $C_1\ C_2\ [d]$)} axiom specifies subsumption from $C_1$ to $C_2$ to degree at least $d$ by means of {\L}ukasiewicz implication.
    
        \item[Zadeh general concept inclusions] The  \kw{(z-implies $C_1\ C_2\ [d]$)} axiom specifies subsumption from $C_1$ to $C_2$ to degree at least $d$ by means of Zadeh implication.
    
        \item[Concept definitions] The  \kw{(define-concept $A\ C$)} axiom gives a complete definition of the atomic concept $A$ by identifying it with the concept expression $C$. By contrast, \kw{(define-primitive-concept $A\ C$)} only states that $A$ is subsumed by $C$. The axiom \kw{(equivalent-concepts $C_1\ C_2$)} declares that $C_1$ and $C_2$ are extensionally equivalent.
    
        \item[Disjointness axioms] The  \kw{(disjoint $C_1\ \ldots\ C_k$)} axiom imposes pairwise disjointness on the listed concepts.
    
        \item[Disjoint union axioms] The  \kw{(disjoint-union $C_1\ \ldots\ C_k$)} axiom states that the first concept is the disjoint union of the remaining concepts, that is, it coincides with their union and they are pairwise disjoint.
    
        \item[Range restrictions] The  \kw{(range $R\ C$)} axiom constrains the range of the role $R$ to the concept $C$.
    
        \item[Domain restrictions] The  \kw{(domain $R\ C$)} axiom constrains the domain of the role $R$ to the concept $C$.
    
        \item[Functional roles] The  \kw{(functional $R$)} axiom declares that the role $R$ is functional.
    
        \item[Inverse-functional roles] The \kw{(inverse-functional $R$)} axiom declares that the role $R$ is inverse-functional.
    
        \item[Reflexive roles] The  \kw{(reflexive $R$)} axiom  declares that every individual is related to itself by means of $R$.
    
        \item[Symmetric roles] The  \kw{(symmetric $R$)} axiom declares that the role $R$ is symmetric.
    
        \item[Transitive roles] The \kw{(transitive $R$)} axiom declares that the role $R$ is transitive.
    
        \item[Role inclusion axioms] The \kw{(implies-role $R_1\ R_2\ [d]$)} axiom states that $R_1$ is included in $R_2$ to degree at least $d$.
    
        \item[Inverse role axioms] The  \kw{(inverse $R_1\ R_2$)} axiom states that $R_1$ and $R_2$ are mutually inverse roles.
    \end{description}

    \input{grammar_axioms}

    \subsection{Show statements and queries}
    \nd For debugging and inspection purposes, \texttt{fuzzyDL} provides a family of \kw{show} statements (see \Cref{tab:show-statements}). For instance,
    \begin{equation*}
        \kw{(show-concrete-fillers $F_1\ \ldots\ F_n$)}
    \end{equation*}
    displays the values of the fillers of the features $F_i$, for $i=1,\ldots,n$, whereas
    \begin{equation*}
        \kw{(show-concrete-fillers-for ind $F_1\ \ldots\ F_n$)}
    \end{equation*}
    shows the values of the fillers of the same features for the individual \texttt{ind}. Moreover,
    \begin{equation*}
        \kw{(show-concrete-instance-for ind $F\ C_1\ \ldots\ C_n$)}
    \end{equation*}
    returns the degrees to which the filler of the feature $F$ for the individual \texttt{ind} is an instance of each concept $C_i$, for $i=1,\ldots,n$.
    
    The language also supports inspection of abstract role fillers. In particular,
    \begin{equation*}
        \kw{(show-abstract-fillers $R_1\ \ldots\ R_n$)}
    \end{equation*}
    shows the fillers of the roles $R_i$, together with their membership to any concept, while
    \begin{equation*}
        \kw{(show-abstract-fillers-for ind $R_1\ \ldots\ R_n$)}
    \end{equation*}
    provides the corresponding information restricted to the individual \texttt{ind}. In addition,
    \begin{equation*}
        \kw{(show-concepts $a_1\ \ldots\ a_n$)}
    \end{equation*}
    shows the membership of the individuals $a_i$ to any concept, for $i=1,\ldots,n$,
    \begin{equation*}
        \kw{(show-instances $C_1\ \ldots\ C_n$)}
    \end{equation*}
    shows the values of the instances of the concepts $C_i$, and
    \begin{equation*}
        \kw{(show-variables $x_1\ \ldots\ x_n$)}
    \end{equation*}
    shows the values of the variables $x_i$.
    
    Queries (see \Cref{tab:query-definitions}) are instead used to test properties of a KB or to optimise variables. Basic consistency is checked by means of 
    \begin{equation*}
     \kw{(sat?)} \ ,
    \end{equation*}
    \nd which determines whether $\mathcal{K}$ is consistent.
    Membership queries are expressed by
    \begin{equation*}
    \kw{(max-instance?~a C)}
    \end{equation*}
    \nd and 
    \begin{equation*}
    \kw{(min-instance?~a C)} \ ,
    \end{equation*}
    \nd which return, respectively,
    \begin{equation*}
        \sup \{n \mid \mathcal{K} \models \text{(instance a C n)}\}
    \end{equation*}
    and
    \begin{equation*}
        \inf \{n \mid \mathcal{K} \models \text{(instance a C n)}\},
    \end{equation*}
    that is, the maximal and minimal degrees to which the individual $a$ is an instance of the concept $C$ w.r.t.~a KB $\mathcal{K}$. Likewise,
    \begin{equation*}
        \kw{(all-instances?~C)}
    \end{equation*}
    computes \kw{(min-instance?~a C)} for every individual of a KB $\mathcal{K}$.
    
    Analogous queries are available for role membership:
    \begin{equation*}
        \kw{(max-related?~a b R)}
    \end{equation*}
    and
    \begin{equation*}
        \kw{(min-related?~a b R)},
    \end{equation*}
    which compute the supremum and infimum of all degrees $n$ such that $\mathcal{K} \models \text{(related a b R n)}$.
    
    Concept subsumption can also be queried. The commands 
    \begin{equation*}
    \kw{(max-subs?~C D)}     
    \end{equation*}
    \nd  and 
    \begin{equation*}
    \kw{(min-subs?~C D)} 
    \end{equation*}
    \nd  return, respectively,
    \begin{equation*}
        \sup \{n \mid \mathcal{K} \models \text{(implies D C n)}\}
    \end{equation*}
    and
    \begin{equation*}
        \inf \{n \mid \mathcal{K} \models \text{(implies D C n)}\}.
    \end{equation*}
    Corresponding variants are provided for specific implication semantics, namely
    \kw{max-g-subs?}, \kw{min-g-subs?}, \kw{max-l-subs?}, \kw{min-l-subs?}, \kw{max-kd-subs?}, and \kw{min-kd-subs?}, based on G\"odel, {\L}ukasiewicz, and Kleene--Dienes implication, respectively.
    
    Satisfiability queries are written as 
    \begin{equation*}
    \kw{(max-sat?~C [a])} 
    \end{equation*}
    \nd and 
    \begin{equation*}
    \kw{(min-sat?~C [a])} \ ,
    \end{equation*}
    \nd which compute, respectively,
    \begin{equation*}
        \sup_{\mathcal{I}} \sup_{a \in \Delta^\mathcal{I}} C^\mathcal{I}(a)
    \end{equation*}
    and
    \begin{equation*}
        \inf_{\mathcal{I}} \inf_{a \in \Delta^\mathcal{I}} C^\mathcal{I}(a).
    \end{equation*}
    If the optional individual $a$ is specified, the query returns the corresponding maximum or minimum satisfiability degree for that individual; otherwise, it returns the solution for an arbitrary individual.
    
    Variable optimisation is supported by 
    \begin{equation*}
    \kw{(max-var?~var)}
    \end{equation*}
    \nd and 
    \begin{equation*}
    \kw{(min-var?~var)} \ ,
    \end{equation*}
    \nd which determine, respectively, the supremum and infimum of the values of \texttt{var} compatible with the consistency of $\mathcal{K}$.
    
    Finally, \texttt{fuzzyDL} provides defuzzification queries for transforming fuzzy outputs associated with concrete features into crisp values. In particular,
    \begin{equation*}
        \kw{(defuzzify-lom?~C a F)}
    \end{equation*}
    defuzzifies the value of $F$ by selecting the largest of the maxima method,
    \begin{equation*}
        \kw{(defuzzify-mom?~C a F)}
    \end{equation*}
    uses the middle of the maxima method, and
    \begin{equation*}
        \kw{(defuzzify-som?~C a F)}
    \end{equation*}
    uses the smallest of the maxima method. Moreover, \kw{(bnp?~f)} computes the Best Non-Fuzzy Performance (BNP) of a fuzzy number $f$ \cite{straccia2014foundations}. 
    
    We remark that, in the defuzzification queries, the concept $C$ represents a collection of Mamdani-style IF--THEN fuzzy rules specifying how the value of the concrete feature $F$ is to be determined.

    \input{grammar_show_filters}
    \input{grammar_queries}

    \subsection{Fuzzy OWL~2 XML annotation syntax} \label{subsec:fuzzyowl2-xml-syntax}
    \nd \texttt{Fuzzy OWL~2} represents fuzzy information by means of string XML-valued OWL~2 annotations attached to an otherwise ordinary OWL~2 ontology. This choice preserves compatibility with standard OWL~2 tools, which may process the underlying crisp ontology without modification, while fuzzy-aware tools can inspect the XML payload and reconstruct the intended fuzzy semantics. The general annotation wrapper is shown in \Cref{listing:fuzzyowl2-wrapper}. In the original \texttt{Fuzzy OWL~2} proposal, fuzzy annotations are enclosed within a \texttt{<FuzzyOwl2>} element whose \texttt{fuzzyType} attribute identifies the nature of the annotated object, such as an ontology-level logic specification, a concept, a datatype, a role, a modifier, or an axiom \cite{bobillo2010representing,BOBILLO20111073}. In the present implementation, the wrapper is likewise denoted by \texttt{fuzzyOwl2}, but it contains a single child element whose tag directly determines the relevant syntactic category.
    The grammar presented in this work should therefore be understood as an implementation-oriented schema fragment for the XML annotation payload, rather than as a full OWL~2 serialisation. More precisely, it specifies the XML content used as the value of the corresponding OWL annotation property.
    
    Ontology-level annotations determine the fuzzy semantics according to which the ontology is to be interpreted (see \Cref{listing:fuzzyowl2-ontology}). The current implementation supports the logics \texttt{lukasiewicz}, \texttt{zadeh}, and \texttt{classical}. If no ontology-level annotation is provided, the default semantics is \texttt{lukasiewicz}, as stated in the documentation \cite{fuzzydlowl2_docs}.
    
    Concept annotations (see \Cref{listing:fuzzyowl2-concepts}) cover a broad class of fuzzy concept constructors, including modified concepts, weighted concepts, weighted aggregations, integral-based aggregations, quantified OWA concepts, and fuzzy nominals. A modified concept applies a previously defined modifier to a base concept. Weighted aggregation concepts contain nested weighted sub-concepts. OWA, Choquet, Sugeno, and quasi-Sugeno concepts include both a sequence of weights and a corresponding sequence of concept names. In the quantified OWA case, the explicit weight vector is replaced by a quantifier. Lastly, nominal concepts associate a degree with a named individual.
    
    Datatype annotations (see \Cref{listing:fuzzyowl2-datatypes}) are used to encode concrete fuzzy membership functions. In the implementation, these are represented by means of a \texttt{Datatype} element equipped with a \texttt{type} attribute specifying the shape of the membership function, together with numerical attributes identifying the relevant breakpoints. The supported forms are crisp, left-shoulder, right-shoulder, triangular, trapezoidal, and linear datatypes, as well as modified datatypes. The latter are obtained by applying a modifier, either linear or triangular, to a previously defined datatype.
    
    Role and modifier annotations (see \Cref{listing:fuzzyowl2-roles-modifiers}) follow the same wrapper convention. A role may be modified by a modifier that has previously been declared, while modifiers themselves may be either linear or triangular. These XML encodings correspond directly to the parenthesised syntax 
    \begin{equation*}
        \kw{(define-modifier $M$ linear-modifier($c$))}
    \end{equation*}
    and 
    \begin{equation*}
        \kw{(define-modifier $M$ triangular-modifier($a,b,c$))}
    \end{equation*}
    introduced earlier.
    
    Finally, axiom annotations (see \Cref{listing:fuzzyowl2-axiom}) associate a fuzzy truth degree with an ordinary OWL~2 axiom. This mechanism is used for fuzzy concept assertions, role assertions, inclusions, and other axiom types whose crisp OWL~2 counterpart remains explicitly present in the ontology.
    
    \begin{listing}
        \caption{Top-level XML annotation wrapper for \texttt{Fuzzy OWL~2} syntax.}
        \label{listing:fuzzyowl2-wrapper}

        
\begin{lstlisting}[
            language=XML,
            frame=lines,
            framesep=2mm,
            basicstyle=\ttfamily\scriptsize\linespread{1.2}\selectfont,
            numbers=left
        ]
<fuzzyOwl2 fuzzyType="ontology|concept|datatype|role|modifier|axiom">
    <!-- Top XML annotation -->
</fuzzyOwl2>
\end{lstlisting}  
        
    \end{listing}
    
    \begin{listing}
        \caption{Ontology-level fuzzy logic annotation.}
        \label{listing:fuzzyowl2-ontology}
        
\begin{lstlisting}[
            language=XML,
            frame=lines,
            framesep=2mm,
            basicstyle=\ttfamily\scriptsize\linespread{1.2}\selectfont,
            numbers=left
        ]
<fuzzyOwl2 fuzzyType="ontology">
    <FuzzyLogic logic="lukasiewicz|goedel|zadeh" />
</fuzzyOwl2>
\end{lstlisting} 
        
    \end{listing}
    
    \begin{listing}
        \caption{Concept-level \texttt{Fuzzy OWL~2} XML annotation syntax.}
        \label{listing:fuzzyowl2-concepts}
        
\begin{lstlisting}[
            language=XML,
            frame=lines,
            framesep=2mm,
            basicstyle=\ttfamily\scriptsize\linespread{1.2}\selectfont,
            numbers=left
        ]
<!-- Modified concept -->
<fuzzyOwl2 fuzzyType="concept">
    <Concept type="modified" modifier="ModifierIRI" base="BaseConceptIRI" />
</fuzzyOwl2>

<!-- Weighted atomic concept -->
<fuzzyOwl2 fuzzyType="concept">
    <Concept type="weighted" value="0.5" base="BaseConceptIRI" />
</fuzzyOwl2>

<!-- Weighted aggregation: type is weightedMinimum, weightedMaximum,
     weightedSum, or weightedSumZero -->
<fuzzyOwl2 fuzzyType="concept">
    <Concept type="weightedMaximum">
        <Concept type="weighted" value="0.5" base="SubConcept1IRI" />
        <Concept type="weighted" value="0.5" base="SubConcept2IRI" />
    </Concept>
</fuzzyOwl2>

<!-- OWA / Choquet / Sugeno / quasi-Sugeno aggregation -->
<fuzzyOwl2 fuzzyType="concept">
    <Concept type="owa|choquet|sugeno|quasisugeno">
        <Weights>
            <Weight>0.6</Weight>
            <Weight>0.4</Weight>
        </Weights>
        <Names>
            <Name>Concept1IRI</Name>
            <Name>Concept2IRI</Name>
        </Names>
    </Concept>
</fuzzyOwl2>

<!-- Quantifier-guided OWA -->
<fuzzyOwl2 fuzzyType="concept">
    <Concept type="qowa" quantifier="QuantifierIRI">
        <Names>
            <Name>Concept1IRI</Name>
            <Name>Concept2IRI</Name>
        </Names>
    </Concept>
</fuzzyOwl2>

<!-- Fuzzy nominal -->
<fuzzyOwl2 fuzzyType="concept">
    <Concept type="nominal" value="0.8" individual="IndividualIRI" />
</fuzzyOwl2>
\end{lstlisting}

    \end{listing}
    
    \begin{listing}
        \caption{Datatype-level \texttt{Fuzzy OWL~2} XML annotation syntax.}
        \label{listing:fuzzyowl2-datatypes}
        
\begin{lstlisting}[
            language=XML,
            frame=lines,
            framesep=2mm,
            basicstyle=\ttfamily\scriptsize\linespread{1.2}\selectfont,
            numbers=left
        ]
<!-- Crisp interval -->
<fuzzyOwl2 fuzzyType="datatype">
    <Datatype type="crisp" a="0" b="1" />
</fuzzyOwl2>

<!-- Left-shoulder membership function -->
<fuzzyOwl2 fuzzyType="datatype">
    <Datatype type="leftshoulder" a="0" b="1" />
</fuzzyOwl2>

<!-- Right-shoulder membership function -->
<fuzzyOwl2 fuzzyType="datatype">
    <Datatype type="rightshoulder" a="0" b="1" />
</fuzzyOwl2>

<!-- Triangular membership function -->
<fuzzyOwl2 fuzzyType="datatype">
    <Datatype type="triangular" a="0" b="0.5" c="1" />
</fuzzyOwl2>

<!-- Trapezoidal membership function -->
<fuzzyOwl2 fuzzyType="datatype">
    <Datatype type="trapezoidal" a="0" b="0.25" c="0.75" d="1" />
</fuzzyOwl2>

<!-- Linear membership function -->
<fuzzyOwl2 fuzzyType="datatype">
    <Datatype type="linear" a="0" b="1" />
</fuzzyOwl2>

<!-- Modified membership function using modifier 'very' and datatype 'High' -->
<fuzzyOwl2 fuzzyType="datatype">
    <Datatype type="modified" modifier="very" base="High" />
</fuzzyOwl2>
\end{lstlisting}

    \end{listing}
    
    \begin{listing}
        \caption{Role and modifier \texttt{Fuzzy OWL~2} XML annotation syntax.}
        \label{listing:fuzzyowl2-roles-modifiers}
        
\begin{lstlisting}[
            language=XML,
            frame=lines,
            framesep=2mm,
            basicstyle=\ttfamily\scriptsize\linespread{1.2}\selectfont,
            numbers=left
        ]
<!-- Modified role -->
<fuzzyOwl2 fuzzyType="role">
    <Role type="modified" modifier="ModifierIRI" base="BaseRoleIRI" />
</fuzzyOwl2>

<!-- Linear modifier -->
<fuzzyOwl2 fuzzyType="modifier">
    <Modifier type="linear" c="0.5" />
</fuzzyOwl2>

<!-- Triangular modifier -->
<fuzzyOwl2 fuzzyType="modifier">
    <Modifier type="triangular" a="0" b="0.5" c="1" />
</fuzzyOwl2>
\end{lstlisting}

    \end{listing}
    
    \begin{listing}
        \caption{Axiom-level degree annotation.}
        \label{listing:fuzzyowl2-axiom}
        
\begin{lstlisting}[
            language=XML,
            frame=lines,
            framesep=2mm,
            basicstyle=\ttfamily\scriptsize\linespread{1.2}\selectfont,
            numbers=left
        ]
<fuzzyOwl2 fuzzyType="axiom">
    <Degree value="0.5" />
</fuzzyOwl2>
\end{lstlisting}

    \end{listing}
    
    \section{\texttt{fuzzyDL} Parser}\label{sec:appendix-b}
    \nd The \texttt{DLParser} class provides a \texttt{pyparsing}-based \cite{pyparsing} front end for the \texttt{fuzzyDL} language. Its behaviour is semantic as well as syntactic: once a directive has been recognised, the associated parse action either creates domain objects (for instance, concepts, modifiers, fuzzy numbers, features, MILP expressions, or query objects), or updates the current \texttt{KnowledgeBase}. The parser is therefore not merely a lexer--parser layer, but the component that initialises the internal semantic representation consumed by the reasoner.
    
    The entry point of the grammar is the method \texttt{get\_grammatics()}, which defines all admissible top-level forms. These include declarations, concept expressions, datatype restrictions, axioms, constraints, inspection commands, and queries. The methods \texttt{parse\_string()} and \texttt{parse\_string\_opt()} execute the grammar on an input string, whereas \texttt{get\_kb()} initialises the parser state, loads the configuration, creates a fresh \texttt{KnowledgeBase}, parses the input file, and returns both the populated KB and the list of parsed queries.
    
    A central feature of the implementation is that concept syntax is interpreted compositionally. Atomic names, truth constants, datatype restrictions, role restrictions, modifiers, thresholds, weighted constructors, OWA / integral operators, and approximation operators are translated directly into semantic objects such as \texttt{OperatorConcept}, \texttt{AllSomeConcept}, \\ \texttt{HasValueConcept}, \texttt{ThresholdConcept}, \texttt{WeightedConcept}, \texttt{OwaConcept}, \\ \texttt{ChoquetIntegral}, \texttt{SugenoIntegral}, \texttt{QsugenoIntegral}, and \texttt{SigmaConcept}. Similarly, numeric and symbolic degrees are normalised into \texttt{DegreeNumeric}, \texttt{DegreeExpression}, or \texttt{DegreeVariable} objects, whilst linear constraints are mapped to \texttt{Variable}, \texttt{Term}, \texttt{Expression}, and \texttt{Inequation} objects in the MILP layer.
    
    At the top level, declarations and axioms immediately modify the KB. Examples include the declaration of fuzzy logic, modifiers, fuzzy concrete concepts, fuzzy numbers, concrete features, crisp concepts or roles, similarity and equivalence relations, and structural axioms such as role properties or concept inclusions. By contrast, query directives do not change the ontology itself; instead, they instantiate specialised query objects, such as \texttt{KbSatisfiableQuery}, \texttt{MaxInstanceQuery}, \texttt{MinSubsumesQuery}, \\ 
    \texttt{MaxRelatedQuery}, \texttt{LomDefuzzifyQuery}, or \texttt{BnpQuery}, which are appended to \texttt{DLParser.queries\_list} for subsequent execution.
    
    \Cref{tab:dl-parser-directives} summarises the directives recognised by the grammar and the corresponding actions performed by the parser. In this context, \texttt{kb} denotes the instance of the class \texttt{KnowledgeBase} stored as an attribute of \texttt{DLParser}. The description follows the current implementation in \texttt{dl\_parser.py}, where the code exhibits an evident implementation caveat, which is stated explicitly.

    \input{appendix_fuzzydl_parser}

    %
        
    {
    \footnotesize
    \bibliographystyle{plain} 
    \bibliography{references}
    }    

\end{document}

%% file: logic_operators.tex
\begin{table}[tbh!]
    \centering
    \renewcommand{\arraystretch}{1.5}
    \rowcolors{2}{gray!10}{white}
    \caption{Fuzzy operators supported by the reasoner.}
    \label{tab:operators}
    \begin{tabular}{@{}lccc@{}}
        \toprule
        \hiderowcolors
        \textbf{Operator} & \textbf{\L{}ukasiewicz} & \textbf{G\"odel} & \textbf{Zadeh} \\
        \showrowcolors
        \midrule
        $\alpha \otimes \beta$      & $\max(\alpha{+}\beta{-}1, 0)$
            & $\min(\alpha,\beta)$ & $\min(\alpha,\beta)$ \\
        $\alpha \oplus \beta$       & $\min(\alpha{+}\beta, 1)$
            & $\max(\alpha,\beta)$ & $\max(\alpha,\beta)$ \\
        $\ominus\,\alpha$              & $1 - \alpha$
            & $\left\{\begin{array}{ll}
                1, & \alpha = 0 \\
                0, & \text{otherwise}
            \end{array}\right.$
            & $1-\alpha$ \\[5pt]
        $\alpha \Rightarrow \beta$  & $\min(1, 1{-}\alpha{+}\beta)$
            & $\left\{\begin{array}{ll}
                1, & \alpha \le \beta \\
                \beta, & \text{otherwise}
            \end{array}\right.$
            & $\max(1{-}\alpha,\beta)$ \\[5pt]
        \bottomrule
    \end{tabular}
\end{table}

%% file: reasoning_services.tex
\begin{table}[tbh!]
    \centering
    \renewcommand{\arraystretch}{1.5}
    \rowcolors{2}{gray!10}{white}
    \caption{Supported reasoning services.  BSD\,=\,best satisfiability degree; MSD\,=\,minimal satisfiability degree; BED\,=\,best entailment degree; MED\,=\,maximal entailment degree; LOM\,=\,largest of maxima; SOM\,=\,smallest of maxima; MOM\,=\,middle of maxima.}
    \label{tab:services}
    \footnotesize
    \begin{tabularx}{\textwidth}{lX}
        \toprule
        \hiderowcolors
        \textbf{Service} & \textbf{Description} \\
        \showrowcolors
        \midrule
        KB consistency        & Existence of a model satisfying every axiom \\
        Concept satisfiability & BSD / MSD of concept w.r.t.~the KB \\
        Concept subsumption   & Degree of a GCI (General Concept Inclusion) \\
        Instance checking     & Min / max degree of a concept assertion \\
        Entailment            & BED / MED of a non-graded axiom \\
        All instances         & All instances of a concept and their BED \\
        Variable optimisation & Min / max of a variable \\
        Defuzzification       & LOM, SOM, or MOM of a concrete role \\
        Classification        & Compute subsumption lattice \\
        BNP                   & Best non-fuzzy performance of a triangular
                                fuzzy number \\
        Related               & Degree of a role assertion \\
        \bottomrule
    \end{tabularx}
\end{table}

%% file: ontoinfo.tex
\begingroup

\scriptsize
\rowcolors{2}{gray!10}{white}
\begin{longtable}{@{}
    >{\raggedright\arraybackslash}p{0.2\textwidth}
    >{\centering\arraybackslash}p{0.1\textwidth}
    >{\centering\arraybackslash}p{0.1\textwidth}
    >{\centering\arraybackslash}p{0.1\textwidth}
    >{\centering\arraybackslash}p{0.1\textwidth}
    >{\centering\arraybackslash}p{0.1\textwidth}
    >{\centering\arraybackslash}p{0.1\textwidth}
@{}}
    \caption{Ontology dataset and statistics.}
    \label{tab:ontos}\\
    
    \toprule
    \hiderowcolors
    \textbf{Ontology} & \textbf{Classes} & \textbf{Object Properties} & \textbf{Data Properties} & \textbf{Individuals} & \textbf{Logical Axioms} & \textbf{Missing Imports} \\
    \showrowcolors
    \midrule
    \endfirsthead

    \toprule
    \hiderowcolors
    \textbf{Ontology} & \textbf{Classes} & \textbf{Object Properties} & \textbf{Data Properties} & \textbf{Individuals} & \textbf{Logical Axioms} & \textbf{Missing Imports} \\
    \showrowcolors
    \midrule
    \endhead

    \bottomrule
    \hiderowcolors
    \multicolumn{2}{r}{\small Continue on next page}\\
    \showrowcolors
    \endfoot
    
    \bottomrule 
    \endlastfoot
            
    \path{AirSystem} & 114 & 111 & 10 & 0 & 1407 &  \\
    \path{amino-acid} & 46 & 5 & 1 & 0 & 563 &  \\
    \path{atom-common} & 14 & 5 & 0 & 0 & 89 &  \\
    \path{cancer_my} & 90 & 13 & 3 & 4 & 200 &  \\
    \path{cancer_ra} & 89 & 13 & 3 & 0 & 184 &  \\
    \path{chebi} & 20979 & 10 & 0 & 0 & 282349 &  \\
    \path{chemical} & 48 & 9 & 11 & 0 & 114 &  \\
    \path{cton} & 17033 & 43 & 0 & 0 & 33203 &  \\
    \path{earthrealm} & 1709 & 183 & 33 & 113 & 3082 & $\star$ \\
    \path{economy} & 339 & 45 & 8 & 482 & 1625 &  \\
    \path{EMAP.obo} & 13731 & 1 & 0 & 0 & 13730 &  \\
    \path{FBbt_XP} & 7225 & 21 & 0 & 12 & 25741 &  \\
    \path{FMA} & 78983 & 8 & 0 & 0 & 168826 &  \\
    \path{fmaOwlDlComponent_1_4_0} & 6488 & 165 & 0 & 98 & 18775 &  \\
    \path{FuzzyWine} & 178 & 15 & 7 & 138 & 678 &  \\
    \path{galen-ians-full-doctored} & 2748 & 413 & 0 & 0 & 4736 &  \\
    \path{gene_ontology_edit.obo} & 26225 & 4 & 0 & 0 & 42656 &  \\
    \path{goslim} & 161 & 0 & 0 & 79 & 237 &  \\
    \path{GRO} & 420 & 18 & 7 & 1 & 827 &  \\
    \path{human_activities} & 162 & 7 & 0 & 0 & 163 & $\star$ \\
    \path{legal-role} & 9 & 1 & 0 & 0 & 9 &  \\
    \path{lubm} & 43 & 25 & 7 & 1555 & 8612 &  \\
    \path{matchmaking} & 108 & 12 & 34 & 0 & 235 &  \\
    \path{mosquito_insecticide_resistance.obo} & 4286 & 10 & 0 & 0 & 4282 &  \\
    \path{mygrid-moby-service} & 36 & 12 & 23 & 0 & 95 & $\star$ \\
    \path{NCI} & 27652 & 70 & 0 & 0 & 46940 &  \\
    \path{numerics} & 5370 & 510 & 65 & 1334 & 10348 & $\star$ \\
    \path{ontology} & 475 & 8 & 0 & 0 & 1959 &  \\
    \path{organic-compound-complex} & 49 & 2 & 0 & 0 & 25 &  \\
    \path{pathway.obo} & 49 & 2 & 0 & 0 & 25 &  \\
    \path{people.fd} & 60 & 14 & 0 & 22 & 107 &  \\
    \path{periodic-table-complex} & 165 & 0 & 0 & 0 & 58 & $\star$ \\
    \path{photography} & 189 & 29 & 16 & 0 & 364 &  \\
    \path{pizza} & 100 & 8 & 0 & 5 & 712 &  \\
    \path{po} & 50 & 44 & 13 & 4 & 183 &  \\
    \path{process} & 1342 & 122 & 26 & 113 & 2313 &  \\
    \path{property} & 1193 & 120 & 26 & 113 & 2153 & $\star$ \\
    \path{propreo} & 401 & 34 & 5 & 46 & 732 &  \\
    \path{propreo.Tbox} & 401 & 34 & 5 & 0 & 732 &  \\
    \path{relative-places.bug} & 7 & 16 & 0 & 0 & 79 &  \\
    \path{relative-places} & 17 & 29 & 0 & 0 & 123 &  \\
    \path{SIGKDD-EKAW} & 123 & 50 & 11 & 0 & 367 &  \\
    \path{so-xp.obo} & 1660 & 22 & 0 & 0 & 1943 &  \\
    \path{space} & 1193 & 120 & 26 & 113 & 2153 & $\star$ \\
    \path{spatial.obo} & 106 & 13 & 0 & 0 & 190 &  \\
    \path{subatomic-particle-complex} & 52 & 1 & 0 & 0 & 75 & $\star$ \\
    \path{tambis-patched} & 395 & 100 & 2 & 0 & 597 &  \\
    \path{teleost_taxonomy.obo} & 36076 & 1 & 0 & 0 & 36069 &  \\
    \path{thesaurus} & 65231 & 125 & 71 & 0 & 94288 &  \\
    \path{Transportation} & 445 & 89 & 4 & 183 & 1157 &  \\
    \path{units} & 12 & 3 & 5 & 103 & 358 &  \\
    \path{worm_phenotype_xp.obo} & 1841 & 31 & 0 & 0 & 1173 &  \\
\end{longtable}
\endgroup

%% file: java_python_benchamarks.tex
\begingroup
\renewcommand{\arraystretch}{2}
\rowcolors{2}{gray!10}{white}
{\scriptsize
\begin{longtable}{@{}%
    >{\arraybackslash}p{0.15\textwidth}%
    >{\centering\arraybackslash}p{0.20\textwidth}%
    >{\centering\arraybackslash}p{0.10\textwidth}%
    >{\centering\arraybackslash}p{0.10\textwidth}%
    >{\centering\arraybackslash}p{0.10\textwidth}%
    >{\centering\arraybackslash}p{0.20\textwidth}%
@{}}
    \caption{Average wall-clock execution time for each KB, together with the ratio between the best Python execution time and the Java execution time. Wall-clock times are reported as mean $\pm$ unbiased standard deviation, expressed in seconds. The best execution time for each Knowledge Base is highlighted in bold. If a Knowledge Base required more than $20$ minutes to complete, its execution time is reported as not available (N/A), and the corresponding ratio is omitted.}
    \label{tab:java-python-benchmarks}\\
    
    \toprule
    \hiderowcolors
    \multirow{2}{=}{\centering\arraybackslash\textbf{KB}} &
    \multirow{2}{=}{\centering\arraybackslash\textbf{Java wall time (s), mean $\pm$ SD}} &
    \multicolumn{3}{c}{\centering\arraybackslash\textbf{Python wall time (s), mean $\pm$ SD}} & \multirow{2}{=}{\centering\arraybackslash\textbf{Ratio (Python / Java wall-time)}}\\
    \cmidrule(lr){3-5}
    {} & {} & \textbf{Gurobi} & \textbf{CBC} & \textbf{HiGHS} & {}\\
    \showrowcolors
    \midrule
    \endfirsthead

    \toprule
    \hiderowcolors
    \multirow{2}{=}{\centering\arraybackslash\textbf{Knowledge Base}} &
    \multirow{2}{=}{\centering\arraybackslash\textbf{Java wall time (s), mean $\pm$ SD}} &
    \multicolumn{3}{c}{\textbf{Python wall time (s), mean $\pm$ SD}} & \multirow{2}{=}{\centering\arraybackslash\textbf{Ratio (Python / Java wall-time)}}\\
    \cmidrule(lr){3-5}
    {} & {} & \textbf{Gurobi} & \textbf{CBC} & \textbf{HiGHS} & {}\\
    \showrowcolors
    \midrule
    \endhead

    \bottomrule
    \hiderowcolors
    \multicolumn{2}{r}{\small Continue on next page}\\
    \showrowcolors
    \endfoot

    \bottomrule
    \endlastfoot
    
    \path{propreo} & \textbf{N/A} & $\mathbf{42.42 \pm 0.18}$ & $156.85 \pm 183.98$ & $48.40 \pm 0.13$ & \textbf{N/A} \\
    \path{organic-compound-complex} & \textbf{N/A} & $\mathbf{35.86 \pm 0.25}$ & $172.09 \pm 49.80$ & $45.95 \pm 1.87$ & \textbf{N/A} \\
    \path{AirSystem} & \textbf{N/A} & $\mathbf{35.52 \pm 2.05}$ & $613.60 \pm 55.56$ & $135.43 \pm 8.65$ & \textbf{N/A} \\
    \path{chebi} & \textbf{N/A} & $\mathbf{9.95 \pm 0.31}$ & $11.31 \pm 0.11$ & $10.64 \pm 0.07$ & \textbf{N/A} \\
    \path{FMA} & $27.45 \pm 2.18$ & $\mathbf{3.38 \pm 0.05}$ & $3.70 \pm 0.04$ & $3.60 \pm 0.04$ & $8.128$ \\
    \path{FuzzyWine} & $62.87 \pm 5.34$ & $\mathbf{14.97 \pm 0.08}$ & $21.14 \pm 0.09$ & $16.77 \pm 0.36$ & $4.199$ \\
    \path{chemical} & $3.38 \pm 0.10$ & $\mathbf{1.13 \pm 0.02}$ & $1.35 \pm 0.05$ & $1.15 \pm 0.04$ & $2.983$ \\
    \path{lubm} & $60.85 \pm 4.09$ & $\mathbf{21.29 \pm 0.16}$ & $26.26 \pm 0.12$ & $21.49 \pm 0.09$ & $2.858$ \\
    \path{propreo.TBox} & $75.66 \pm 3.59$ & $39.52 \pm 0.56$ & $91.08 \pm 61.24$ & $\mathbf{27.94 \pm 18.44}$ & $2.708$ \\
    \path{FBbt_XP} & $1.59 \pm 0.04$ & $\mathbf{1.25 \pm 0.05}$ & $1.31 \pm 0.05$ & $1.27 \pm 0.05$ & $1.269$ \\
    \path{cancer_my} & $\mathbf{0.87 \pm 0.04}$ & $1.02 \pm 0.03$ & $1.38 \pm 0.08$ & $1.25 \pm 0.05$ & $0.851$ \\
    \path{NCI} & $\mathbf{0.56 \pm 0.08}$ & $1.09 \pm 0.04$ & $1.14 \pm 0.05$ & $1.10 \pm 0.04$ & $0.510$ \\
    \path{gene_ontology_edit.obo} & $\mathbf{0.47 \pm 0.04}$ & $1.00 \pm 0.05$ & $1.03 \pm 0.05$ & $1.02 \pm 0.05$ & $0.466$ \\
    \path{earthrealm} & $\mathbf{0.33 \pm 0.02}$ & $0.72 \pm 0.04$ & $0.82 \pm 0.04$ & $0.79 \pm 0.04$ & $0.462$ \\
    \path{numerics} & $\mathbf{0.33 \pm 0.03}$ & $0.73 \pm 0.03$ & $0.80 \pm 0.03$ & $0.76 \pm 0.04$ & $0.459$ \\
    \path{process} & $\mathbf{0.34 \pm 0.02}$ & $0.75 \pm 0.05$ & $0.87 \pm 0.03$ & $0.82 \pm 0.05$ & $0.456$ \\
    \path{property} & $\mathbf{0.33 \pm 0.02}$ & $0.74 \pm 0.03$ & $0.84 \pm 0.05$ & $0.84 \pm 0.05$ & $0.452$ \\
    \path{teleost_taxonomy.obo} & $\mathbf{0.45 \pm 0.03}$ & $1.01 \pm 0.04$ & $1.05 \pm 0.04$ & $1.02 \pm 0.05$ & $0.445$ \\
    \path{thesaurus} & $\mathbf{0.99 \pm 0.05}$ & $2.27 \pm 0.06$ & $2.29 \pm 0.08$ & $2.30 \pm 0.11$ & $0.439$ \\
    \path{cton} & $\mathbf{0.40 \pm 0.02}$ & $0.93 \pm 0.05$ & $0.98 \pm 0.05$ & $0.99 \pm 0.04$ & $0.431$ \\
    \path{amino-acid} & $\mathbf{0.31 \pm 0.05}$ & $0.74 \pm 0.02$ & $0.95 \pm 0.05$ & $0.94 \pm 0.06$ & $0.425$ \\
    \path{space} & $\mathbf{0.33 \pm 0.02}$ & $0.79 \pm 0.04$ & $0.88 \pm 0.06$ & $0.80 \pm 0.04$ & $0.422$ \\
    \path{fmaOwlDlComponent_1_4_0} & $\mathbf{0.31 \pm 0.02}$ & $0.79 \pm 0.04$ & $0.83 \pm 0.03$ & $0.81 \pm 0.02$ & $0.393$ \\
    \path{people.fd} & $\mathbf{0.30 \pm 0.01}$ & $0.79 \pm 0.04$ & $0.93 \pm 0.07$ & $0.87 \pm 0.04$ & $0.383$ \\
    \path{EMAP.obo} & $\mathbf{0.28 \pm 0.03}$ & $0.79 \pm 0.05$ & $0.84 \pm 0.04$ & $0.84 \pm 0.07$ & $0.349$ \\
    \path{economy} & $\mathbf{0.21 \pm 0.01}$ & $0.68 \pm 0.03$ & $0.75 \pm 0.04$ & $0.71 \pm 0.03$ & $0.302$ \\
    \path{po} & $\mathbf{0.24 \pm 0.01}$ & $0.82 \pm 0.07$ & $1.90 \pm 0.26$ & $1.22 \pm 0.23$ & $0.288$ \\
    \path{subatomic-particle-complex} & $\mathbf{0.19 \pm 0.02}$ & $0.71 \pm 0.07$ & $0.76 \pm 0.06$ & $0.72 \pm 0.02$ & $0.271$ \\
    \path{mosquito_insecticide_resistance.obo} & $\mathbf{3.29 \pm 0.09}$ & $12.21 \pm 0.16$ & $12.30 \pm 0.09$ & $12.23 \pm 0.06$ & $0.270$ \\
    \path{galen-ians-full-doctored} & $\mathbf{0.23 \pm 0.02}$ & $0.86 \pm 0.03$ & $0.89 \pm 0.04$ & $0.90 \pm 0.04$ & $0.270$ \\
    \path{so-xp.obo} & $\mathbf{0.17 \pm 0.02}$ & $0.66 \pm 0.04$ & $0.74 \pm 0.04$ & $0.70 \pm 0.04$ & $0.263$ \\
    \path{Transportation} & $\mathbf{0.17 \pm 0.02}$ & $0.64 \pm 0.04$ & $0.71 \pm 0.03$ & $0.68 \pm 0.03$ & $0.262$ \\
    \path{relative-places} & $\mathbf{0.16 \pm 0.03}$ & $0.66 \pm 0.03$ & $0.82 \pm 0.05$ & $0.74 \pm 0.05$ & $0.244$ \\
    \path{relative-places.bug} & $\mathbf{0.16 \pm 0.01}$ & $0.65 \pm 0.02$ & $0.80 \pm 0.03$ & $0.71 \pm 0.05$ & $0.241$ \\
    \path{units} & $\mathbf{0.16 \pm 0.02}$ & $0.66 \pm 0.03$ & $0.70 \pm 0.02$ & $0.69 \pm 0.04$ & $0.238$ \\
    \path{GRO} & $\mathbf{0.15 \pm 0.00}$ & $0.62 \pm 0.03$ & $0.66 \pm 0.04$ & $0.65 \pm 0.03$ & $0.234$ \\
    \path{ontology} & $\mathbf{0.14 \pm 0.01}$ & $0.62 \pm 0.03$ & $0.64 \pm 0.03$ & $0.66 \pm 0.02$ & $0.233$ \\
    \path{worm_phenotype_xp.obo} & $\mathbf{0.16 \pm 0.01}$ & $0.70 \pm 0.04$ & $0.74 \pm 0.04$ & $0.73 \pm 0.03$ & $0.223$ \\
    \path{mygrid-moby-service} & $\mathbf{0.14 \pm 0.01}$ & $0.61 \pm 0.02$ & $0.65 \pm 0.03$ & $0.65 \pm 0.03$ & $0.222$ \\
    \path{cancer_ra} & $\mathbf{0.13 \pm 0.02}$ & $0.64 \pm 0.05$ & $0.69 \pm 0.05$ & $0.64 \pm 0.04$ & $0.212$ \\
    \path{periodic-table-complex} & $\mathbf{0.15 \pm 0.01}$ & $0.71 \pm 0.03$ & $0.76 \pm 0.07$ & $0.72 \pm 0.06$ & $0.211$ \\
    \path{pizza} & $\mathbf{0.12 \pm 0.01}$ & $0.62 \pm 0.03$ & $0.72 \pm 0.04$ & $0.69 \pm 0.04$ & $0.200$ \\
    \path{tambis-patched} & $\mathbf{0.13 \pm 0.01}$ & $0.66 \pm 0.04$ & $0.71 \pm 0.03$ & $0.68 \pm 0.04$ & $0.197$ \\
    \path{atom-common} & $\mathbf{0.12 \pm 0.01}$ & $0.60 \pm 0.03$ & $0.66 \pm 0.03$ & $0.64 \pm 0.03$ & $0.196$ \\
    \path{pathway.obo} & $\mathbf{0.12 \pm 0.02}$ & $0.64 \pm 0.04$ & $0.69 \pm 0.07$ & $0.64 \pm 0.04$ & $0.195$ \\
    \path{goslim} & $\mathbf{0.12 \pm 0.01}$ & $0.61 \pm 0.02$ & $0.66 \pm 0.02$ & $0.64 \pm 0.04$ & $0.194$ \\
    \path{human_activities} & $\mathbf{0.11 \pm 0.01}$ & $0.61 \pm 0.03$ & $0.66 \pm 0.04$ & $0.64 \pm 0.03$ & $0.178$ \\
    \path{SIGKDD-EKAW} & $\mathbf{0.11 \pm 0.01}$ & $0.63 \pm 0.03$ & $0.66 \pm 0.02$ & $0.64 \pm 0.04$ & $0.177$ \\
    \path{photography} & $\mathbf{0.11 \pm 0.01}$ & $0.65 \pm 0.05$ & $0.67 \pm 0.04$ & $0.66 \pm 0.03$ & $0.176$ \\
    \path{matchmaking} & $\mathbf{0.10 \pm 0.01}$ & $0.58 \pm 0.02$ & $0.64 \pm 0.03$ & $0.62 \pm 0.02$ & $0.175$ \\
    \path{spatial.obo} & $\mathbf{0.11 \pm 0.01}$ & $0.66 \pm 0.06$ & $0.67 \pm 0.04$ & $0.66 \pm 0.04$ & $0.167$ \\
    \path{legal-role} & $\mathbf{0.09 \pm 0.01}$ & $0.61 \pm 0.03$ & $0.65 \pm 0.03$ & $0.63 \pm 0.03$ & $0.147$ \\

\end{longtable}
}
\endgroup

%% file: grammar_ebnf_1.tex
\begingroup

\scriptsize
\rowcolors{2}{gray!10}{white}
\begin{longtable}{@{}
    >{\raggedright\arraybackslash}p{0.3\textwidth}
    >{\raggedright\arraybackslash}p{0.6\textwidth}
@{}}
    \caption{Lexical and top-level fuzzyDL grammar in EBNF.} \label{tab:ebnf-lexical-top-level}\\
    
    \toprule
    \hiderowcolors
    \textbf{Nonterminal} & \textbf{EBNF production} \\
    \showrowcolors
    \midrule
    \endfirsthead

    \toprule
    \hiderowcolors
    \textbf{Nonterminal} & \textbf{EBNF production} \\
    \showrowcolors
    \midrule
    \endhead

    \bottomrule
    \hiderowcolors
    \multicolumn{2}{r}{\small Continue on next page}\\
    \showrowcolors
    \endfoot
    
    \bottomrule 
    \endlastfoot
    
    $\nt{name}$
    &
    \texttt{["] [a-zA-Z\_] [a-zA-Z0-9\_]* ["]} \\
    
    $\nt{number}$
    &
    \texttt{[+-]? [0-9]+ (\textbackslash.[0-9]+)} \\
    
    $\nt{logic}$
    &
    \kw{lukasiewicz} $|$ \kw{zadeh} $|$ \kw{classical} \\
    
    $\nt{knowledge\_base}$
    &
    $\{\,\nt{declaration} \mid \nt{constraint} \mid \nt{statement} \mid
    \nt{axiom} \mid \nt{query}\,\}^{*}$ \\
    
    $\nt{declaration}$
    &
    \ensuremath{\begin{array}{l}
        \nt{define\_logic} \\
        \mid\ \nt{truth\_constant}\\
        \mid\ \nt{modifier} \\
        \mid\ \nt{fuzzy\_concept}\\
        \mid\ \nt{fuzzy\_number\_range}\\
        \mid\ \nt{fuzzy\_number}\\
        \mid\ \nt{feature}\\
        \mid\ \nt{feature\_range}\\
        \mid\ \nt{crisp\_declaration}\\
        \mid\ \nt{fuzzy\_relation} \\
    \end{array}}\\
    
    $\nt{define\_logic}$
    &
    \kw{(define-fuzzy-logic} $\nt{logic}$\kw{)} \\
    
    $\nt{truth\_constant}$
    &
    \kw{(define-truth-constant} $\nt{name}$ $\nt{number}$\kw{)} \\
    
    $\nt{modifier}$
    &
    \kw{(define-modifier} $\nt{name}$ $\nt{modifier\_type}$\kw{)} \\
    
    $\nt{modifier\_type}$
    &
    \ensuremath{\begin{array}{l}
        \kw{linear-modifier(}\nt{number}\kw{)}\\
        \mid\ \kw{triangular-modifier(} \nt{number}\ {\left\{\kw{,}\ \nt{number}\right\}}^2\kw{)}\\
    \end{array}} \\
    
    $\nt{crisp\_declaration}$
    &
    \ensuremath{\begin{array}{l}
        \kw{(crisp-concept} \{\nt{name}\}^{+}\kw{)}\\
        \mid\ \kw{(crisp-role} \{\nt{name}\}^{+}\kw{)}\\
    \end{array}} \\
    
    $\nt{fuzzy\_relation}$
    &
    \ensuremath{\begin{array}{l}
        \kw{(define-fuzzy-similarity} \nt{name}\kw{)}\\
        \mid\ \kw{(define-fuzzy-equivalence} \nt{name}\kw{)}\\
    \end{array}} \\
\end{longtable}
\endgroup

%% file: grammar_ebnf_2.tex
\begingroup
\scriptsize
\rowcolors{2}{gray!10}{white}
\begin{longtable}{@{}
    >{\raggedright\arraybackslash}p{0.3\textwidth}
    >{\raggedright\arraybackslash}p{0.6\textwidth}
@{}}
    \caption{Concrete fuzzy concepts, fuzzy numbers, features, and datatype restrictions in EBNF.}
    \label{tab:ebnf-datatypes}\\
    
    \toprule
    \hiderowcolors
    \textbf{Nonterminal} & \textbf{EBNF production} \\
    \showrowcolors
    \midrule
    \endfirsthead

    \toprule
    \hiderowcolors
    \textbf{Nonterminal} & \textbf{EBNF production} \\
    \showrowcolors
    \midrule
    \endhead

    \bottomrule
    \hiderowcolors
    \multicolumn{2}{r}{\small Continue on next page}\\
    \showrowcolors
    \endfoot

    \bottomrule
    \endlastfoot
    
    $\nt{fuzzy\_concept}$
    &
    \kw{(define-fuzzy-concept} $\nt{name}$ $\nt{concept\_type}$\kw{)} \\
    
    $\nt{concept\_type}$
    &
    \ensuremath{\begin{array}{l}
        \kw{crisp(}\nt{number}\ {\{\kw{,}\ \nt{number}\}}^{3}\kw{)} \\
        \mid\ \kw{left-shoulder(}\nt{number}\ {\{\kw{,}\ \nt{number}\}}^{3}\kw{)} \\
        \mid\ \kw{right-shoulder(}\nt{number}\ {\{\kw{,}\ \nt{number}\}}^{3}\kw{)} \\
        \mid\ \kw{triangular(}\nt{number}\ {\{\kw{,}\ \nt{number}\}}^{3}\kw{)} \\
        \mid\ \kw{trapezoidal(}\nt{number}\ {\{\kw{,}\ \nt{number}\}}^{4}\kw{)} \\
        \mid\ \kw{linear(}\nt{number}\ {\{\kw{,}\ \nt{number}\}}^{3}\kw{)} \\
        \mid\ \kw{modified(}\nt{name}\kw{,}\ \nt{name}\kw{)}
    \end{array}} \\
    
    $\nt{fuzzy\_number\_range}$
    &
    \kw{(define-fuzzy-number-range} $\nt{number}$ $\nt{number}$\kw{)} \\
    
    $\nt{fuzzy\_number}$
    &
    \kw{(define-fuzzy-number} $\nt{name}$ $\nt{fuzzy\_number\_expression}$\kw{)} \\
    
    $\nt{fuzzy\_number\_expression}$
    &\ensuremath{\begin{array}{l}
        \nt{name} \\
        \mid\ \nt{number} \\
        \mid\ \kw{(}\nt{number}\kw{,}\ \nt{number}\kw{,}\ \nt{number}\kw{)} \\
        \mid\ \kw{(f+}\ \{\nt{fuzzy\_number\_expression}\}^{+}\kw{)} \\
        \mid\ \kw{(f-}\ \nt{fuzzy\_number\_expression}\ \nt{fuzzy\_number\_expression}\kw{)} \\
        \mid\ \kw{(f*}\ \{\nt{fuzzy\_number\_expression}\}^{+}\kw{)} \\
        \mid\ \kw{(f/}\ \nt{fuzzy\_number\_expression}\ \nt{fuzzy\_number\_expression}\kw{)}
    \end{array}}\\
    
    $\nt{feature}$
    &
    \kw{(functional} $\nt{name}$\kw{)} \\
    
    $\nt{feature\_range}$
    &
    \ensuremath{\begin{array}{l}
        \kw{(range}\ \nt{name}\ \kw{*integer*}\ \nt{number}\ \nt{number}\kw{)} \\
        \mid\ \kw{(range}\ \nt{name}\ \kw{*real*}\ \nt{number}\ \nt{number}\kw{)} \\
        \mid\ \kw{(range}\ \nt{name}\ \kw{*string*}\kw{)} \\
        \mid\ \kw{(range}\ \nt{name}\ \kw{*boolean*}\kw{)}
    \end{array}}\\
    
    $\nt{restriction}$
    &
    \ensuremath{\begin{array}{l}
        \kw{(}\nt{restriction\_operator}\ \nt{name}\ \nt{name}\kw{)} \\
        \mid\ \kw{(}\nt{restriction\_operator}\ \nt{name}\ \nt{restriction\_function}\kw{)} \\
        \mid\ \kw{(}\nt{restriction\_operator}\ \nt{name}\ \nt{fuzzy\_number}\kw{)}
    \end{array}}\\
    
    $\nt{restriction\_operator}$
    &
    \kw{>=} $|$ \kw{<=} $|$ \kw{=} \\
    
    $\nt{restriction\_function}$
    &\ensuremath{\begin{array}{l}
        \nt{number} \\
        \mid\ \nt{name} \\
        \mid\ \nt{number}\ [\kw{*}]?\ \nt{restriction\_function} \\
        \mid\ \nt{restriction\_function}\ \kw{-}\ \nt{restriction\_function} \\
        \mid\ \{\nt{restriction\_function}\ \kw{+}\}^{+}\ \nt{restriction\_function}
    \end{array}}\\
\end{longtable}
\endgroup

%% file: grammar_ebnf_3.tex
\begingroup
\scriptsize
\rowcolors{2}{gray!10}{white}
\begin{longtable}{@{}
    >{\raggedright\arraybackslash}p{0.3\textwidth}
    >{\raggedright\arraybackslash}p{0.6\textwidth}
@{}}
    \caption{Concept-expression grammar in EBNF.}
    \label{tab:ebnf-concepts}\\
    
    \toprule
    \hiderowcolors
    \textbf{Nonterminal} & \textbf{EBNF production} \\
    \showrowcolors
    \midrule
    \endfirsthead

    \toprule
    \hiderowcolors
    \textbf{Nonterminal} & \textbf{EBNF production} \\
    \showrowcolors
    \midrule
    \endhead

    \bottomrule
    \hiderowcolors
    \multicolumn{2}{r}{\small Continue on next page}\\
    \showrowcolors
    \endfoot

    \bottomrule
    \endlastfoot
    
    $\nt{concept}$
    &
    \ensuremath{\begin{array}{l}
        \kw{*top*} \\
        \mid\ \kw{*bottom*} \\
        \mid\ \nt{name} \\
        \mid\ \nt{restriction} \\
        \mid\ \kw{(and}\ \nt{concept}\ \nt{concept}\kw{)} \\
        \mid\ \kw{(g-and}\ \nt{concept}\ \nt{concept}\kw{)} \\
        \mid\ \kw{(l-and}\ \nt{concept}\ \nt{concept}\kw{)} \\
        \mid\ \kw{(or}\ \nt{concept}\ \nt{concept}\kw{)} \\
        \mid\ \kw{(g-or}\ \nt{concept}\ \nt{concept}\kw{)} \\
        \mid\ \kw{(l-or}\ \nt{concept}\ \nt{concept}\kw{)} \\
        \mid\ \kw{(not}\ \nt{concept}\kw{)} \\
        \mid\ \kw{(implies}\ \nt{concept}\ \nt{concept}\kw{)} \\
        \mid\ \kw{(g-implies}\ \nt{concept}\ \nt{concept}\kw{)} \\
        \mid\ \kw{(l-implies}\ \nt{concept}\ \nt{concept}\kw{)} \\
        \mid\ \kw{(kd-implies}\ \nt{concept}\ \nt{concept}\kw{)} \\
        \mid\ \kw{(all}\ \nt{name}\ \nt{concept}\kw{)} \\
        \mid\ \kw{(some}\ \nt{name}\ \nt{concept}\kw{)} \\
        \mid\ \kw{(some}\ \nt{name}\ \nt{name}\kw{)} \\
        \mid\ \kw{(ua}\ \nt{name}\ \nt{concept}\kw{)} \\
        \mid\ \kw{(lua}\ \nt{name}\ \nt{concept}\kw{)} \\
        \mid\ \kw{(tua}\ \nt{name}\ \nt{concept}\kw{)} \\
        \mid\ \kw{(la}\ \nt{name}\ \nt{concept}\kw{)} \\
        \mid\ \kw{(lla}\ \nt{name}\ \nt{concept}\kw{)} \\
        \mid\ \kw{(tla}\ \nt{name}\ \nt{concept}\kw{)} \\
        \mid\ \kw{(self}\ \nt{concept}\kw{)} \\
        \mid\ \kw{(}\nt{name}\ \nt{concept}\kw{)} \\
        \mid\ \kw{(}\nt{fuzzy\_number}\kw{)} \\
        \mid\ \kw{([>=}\ \nt{name}\ \kw{]}\ \nt{concept}\kw{)} \\
        \mid\ \kw{([<=}\ \nt{name}\ \kw{]}\ \nt{concept}\kw{)} \\
        \mid\ \kw{(}\nt{number}\ \nt{concept}\kw{)} \\
        \mid\ \kw{(w-sum}\ \{\kw{(}\nt{number}\ \nt{concept}\kw{)}\}^{+}\kw{)} \\
        \mid\ \kw{(w-max}\ \{\kw{(}\nt{number}\ \nt{concept}\kw{)}\}^{+}\kw{)} \\
        \mid\ \kw{(w-min}\ \{\kw{(}\nt{number}\ \nt{concept}\kw{)}\}^{+}\kw{)} \\
        \mid\ \kw{(w-sum-zero}\ \{\kw{(}\nt{number}\ \nt{concept}\kw{)}\}^{+}\kw{)} \\
        \mid\ \kw{(owa}\ \kw{(}\{\nt{number}\}^{+}\kw{)}\ \kw{(}\{\nt{concept}\}^{+}\kw{)} \kw{)} \\
        \mid\ \kw{(q-owa}\ \nt{name}\ \{\nt{concept}\}^{+}\kw{)} \\
        \mid\ \kw{(choquet}\ \kw{(}\{\nt{number}\}^{+}\kw{)}\ \kw{(}\{\nt{concept}\}^{+}\kw{)} \kw{)} \\
        \mid\ \kw{(sugeno}\ \kw{(}\{\nt{number}\}^{+}\kw{)}\ \kw{(}\{\nt{concept}\}^{+}\kw{)}\kw{)} \\
        \mid\ \kw{(q-sugeno}\ \kw{(}\{\nt{number}\}^{+}\kw{)}\ \kw{(}\{\nt{concept}\}^{+}\kw{)}\kw{)} \\
        \mid\ \kw{(sigma-count}\ \nt{name}\ \nt{concept}\ \kw{\{}\{\nt{name}\}^{+}\kw{\}}\ \nt{name}\kw{)}
    \end{array}}
\end{longtable}
\endgroup

%% file: grammar_ebnf_4.tex
\begingroup

\rowcolors{2}{gray!10}{white}
\scriptsize
\begin{longtable}{@{}
    >{\raggedright\arraybackslash}p{0.3\textwidth}
    >{\raggedright\arraybackslash}p{0.6\textwidth}
@{}}
    \caption{Axioms, constraints, statements, and queries in EBNF.}
    \label{tab:ebnf-axioms-queries}\\
    
    \toprule
    \hiderowcolors
    \textbf{Nonterminal} & \textbf{EBNF production} \\
    \showrowcolors
    \midrule
    \endfirsthead

    \toprule
    \hiderowcolors
    \textbf{Nonterminal} & \textbf{EBNF production} \\
    \showrowcolors
    \midrule
    \endhead

    \bottomrule
    \hiderowcolors
    \multicolumn{2}{r}{\small Continue on next page}\\
    \showrowcolors
    \endfoot

    \bottomrule
    \endlastfoot
    
    $\nt{degree}$
    &
    $\nt{number} \mid \nt{expression} \mid \nt{name}$ \\
    
    $\nt{axiom}$
    &
    \ensuremath{\begin{array}{l}
        \kw{(instance}\ \nt{name}\ \nt{concept}\ [\nt{degree}]?\kw{)} \\
        \mid\ \kw{(related}\ \nt{name}\ \nt{name}\ \nt{name}\ [\nt{degree}]?\kw{)} \\
        \mid\ \kw{(implies}\ \nt{concept}\ \nt{concept}\ [\nt{number}]?\kw{)} \\
        \mid\ \kw{(g-implies}\ \nt{concept}\ \nt{concept}\ [\nt{number}]?\kw{)} \\
        \mid\ \kw{(kd-implies}\ \nt{concept}\ \nt{concept}\ [\nt{number}]?\kw{)} \\
        \mid\ \kw{(l-implies}\ \nt{concept}\ \nt{concept}\ [\nt{number}]?\kw{)} \\
        \mid\ \kw{(z-implies}\ \nt{concept}\ \nt{concept}\ [\nt{number}]?\kw{)} \\
        \mid\ \kw{(define-concept}\ \nt{name}\ \nt{concept}\kw{)} \\
        \mid\ \kw{(define-primitive-concept}\ \nt{name}\ \nt{concept}\kw{)} \\
        \mid\ \kw{(equivalent-concepts}\ \nt{concept}\ \nt{concept}\kw{)} \\
        \mid\ \kw{(disjoint}\ \{\nt{concept}\}^{+}\kw{)} \\
        \mid\ \kw{(disjoint-union}\ \{\nt{concept}\}^{+}\kw{)} \\
        \mid\ \kw{(range}\ \nt{name}\ \nt{concept}\kw{)} \\
        \mid\ \kw{(domain}\ \nt{name}\ \nt{concept}\kw{)} \\
        \mid\ \kw{(functional}\ \nt{name}\kw{)} \\
        \mid\ \kw{(inverse-functional}\ \nt{name}\kw{)} \\
        \mid\ \kw{(reflexive}\ \nt{name}\kw{)} \\
        \mid\ \kw{(symmetric}\ \nt{name}\kw{)} \\
        \mid\ \kw{(transitive}\ \nt{name}\kw{)} \\
        \mid\ \kw{(implies-role}\ \nt{name}\ \nt{name}\ [\nt{number}]?\kw{)} \\
        \mid\ \kw{(inverse}\ \nt{name}\ \nt{name}\kw{)}
    \end{array}}\\
    
    $\nt{operator}$
    &
    \kw{>=} $|$ \kw{<=} $|$ \kw{=} \\
    
    $\nt{term}$
    &
    \ensuremath{\begin{array}{l}
        \nt{number} \\
        \mid\ \nt{name} \\
        \mid\ \nt{number}\ [\kw{*}]\ \nt{term} \\
        \mid\ \nt{name}\ [\kw{*}]\ \nt{term}
    \end{array}}\\
    
    $\nt{expression}$
    &
    $\{\nt{term}\ \kw{+}\}^{+}$ $\nt{term}$ \\
    
    $\nt{inequality\_constraint}$
    &
    $\nt{expression}$ $\nt{operator}$ $\nt{number}$ \\
    
    $\nt{constraint}$
    &
    \ensuremath{\begin{array}{l}
        \kw{(constraints}\ \nt{inequality\_constraint}\kw{)} \\
        \mid\ \kw{(constraints}\ \kw{binary}\ \nt{name}\kw{)} \\
        \mid\ \kw{(constraints}\ \kw{free}\ \nt{name}\kw{)}
    \end{array}}\\
    
    $\nt{statement}$
    &
    \ensuremath{\begin{array}{l}
        \kw{(show-concrete-fillers}\ \{\nt{name}\}^{+}\kw{)} \\
        \mid\ \kw{(show-concrete-fillers-for}\ \{\nt{name}\}^{2+}\kw{)} \\
        \mid\ \kw{(show-concrete-instance-for}\ \{\nt{name}\}^{3+}\kw{)} \\
        \mid\ \kw{(show-abstract-fillers}\ \{\nt{name}\}^{+}\kw{)} \\
        \mid\ \kw{(show-abstract-fillers-for}\ \{\nt{name}\}^{2+}\kw{)} \\
        \mid\ \kw{(show-concepts}\ \{\nt{name}\}^{+}\kw{)} \\
        \mid\ \kw{(show-instances}\ \{\nt{name}\}^{+}\kw{)} \\
        \mid\ \kw{(show-variables}\ \{\nt{name}\}^{+}\kw{)} \\
        \mid\ \kw{(show-language)}
    \end{array}}\\
    
    $\nt{query}$
    & \ensuremath{\begin{array}{l}
        \kw{(sat?)} \\
        \mid\ \kw{(max-instance?}\ \nt{name}\ \nt{concept}\kw{)} \\
        \mid\ \kw{(min-instance?}\ \nt{name}\ \nt{concept}\kw{)} \\
        \mid\ \kw{(all-instances?}\ \nt{concept}\kw{)} \\
        \mid\ \kw{(max-related?}\ \nt{name}\ \nt{name}\ \nt{name}\kw{)} \\
        \mid\ \kw{(min-related?}\ \nt{name}\ \nt{name}\ \nt{name}\kw{)} \\
        \mid\ \kw{(max-subs?}\ \nt{concept}\ \nt{concept}\kw{)} \\
        \mid\ \kw{(min-subs?}\ \nt{concept}\ \nt{concept}\kw{)} \\
        \mid\ \kw{(max-g-subs?}\ \nt{concept}\ \nt{concept}\kw{)} \\
        \mid\ \kw{(min-g-subs?}\ \nt{concept}\ \nt{concept}\kw{)} \\
        \mid\ \kw{(max-l-subs?}\ \nt{concept}\ \nt{concept}\kw{)} \\
        \mid\ \kw{(min-l-subs?}\ \nt{concept}\ \nt{concept}\kw{)} \\
        \mid\ \kw{(max-kd-subs?}\ \nt{concept}\ \nt{concept}\kw{)} \\
        \mid\ \kw{(min-kd-subs?}\ \nt{concept}\ \nt{concept}\kw{)} \\
        \mid\ \kw{(max-sat?}\ \nt{concept}\ [\nt{name}]?\kw{)} \\
        \mid\ \kw{(min-sat?}\ \nt{concept}\ [\nt{name}]?\kw{)} \\
        \mid\ \kw{(max-var?}\ \nt{name}\kw{)} \\
        \mid\ \kw{(min-var?}\ \nt{name}\kw{)} \\
        \mid\ \kw{(defuzzify-lom?}\ \nt{concept}\ \nt{name}\ \nt{name}\kw{)} \\
        \mid\ \kw{(defuzzify-mom?}\ \nt{concept}\ \nt{name}\ \nt{name}\kw{)} \\
        \mid\ \kw{(defuzzify-som?}\ \nt{concept}\ \nt{name}\ \nt{name}\kw{)} \\
        \mid\ \kw{(bnp?}\ \nt{name}\kw{)}
    \end{array}}\\
\end{longtable}
\endgroup

%% file: grammar_fuzzy_numbers.tex
\begin{table}[!htb]
    \centering
    \renewcommand{\arraystretch}{1.5}
    \rowcolors{2}{gray!10}{white}
    \caption{Definitions of fuzzy-number expressions.}
    \label{tab:fuzzy-number-definitions}
    {\footnotesize
    \begin{tabular}{lll}
        \toprule
        \hiderowcolors
        \textbf{Example} & {} & \textbf{Definition} \\
        \showrowcolors
        \midrule
        (a, b, c) & triangular fuzzy number & $\displaystyle (a,b,c)$ \\
        n & real number & $\displaystyle (n,n,n)$ \\
        (f+ f1 f2 $\ldots$ fn) & sums fuzzy numbers elementwise & $\displaystyle \left(\sum_{i = 0}^{n} a_i, \sum_{i = 0}^{n} b_i, \sum_{i = 0}^{n} c_i\right)$ \\
        (f- f1 f2) & subtracts two fuzzy numbers & $\displaystyle (a_1-c_2, b_1 - b_2, c_1 - a_2)$ \\
        (f* f1 f2 $\ldots$ fn) & multiplies fuzzy numbers elementwise & $\displaystyle \left(\prod_{i = 0}^{n} a_i, \prod_{i = 0}^{n} b_i, \prod_{i = 0}^{n} c_i\right)$ \\
        (f/ f1 f2) & divides two fuzzy numbers & $\displaystyle \left(\frac{a_1}{c_2}, \frac{b_1}{b_2}, \frac{c_1}{a_2}\right)$ \\
        \bottomrule
    \end{tabular}
    }
\end{table}

%% file: grammar_feature_definitions.tex
\begin{table}[!htb]
    \centering
    \renewcommand{\arraystretch}{1.5}
    \rowcolors{2}{gray!10}{white}
    \caption{Definitions of features and feature ranges.}
    \label{tab:feature-definitions}
    {\footnotesize
    \begin{tabular}{ll}
        \toprule
        \hiderowcolors
        \textbf{Rule} & \textbf{Meaning} \\
        \showrowcolors
        \midrule
        (functional F) & Define the feature F \\
        (range F \texttt{*integer*} $k_1$ $k_2$) & The range of $F$ is an integer number in $[k_1, k_2]$ \\
        (range F \texttt{*real*} $k_1$ $k_2$) & The range of $F$ is a rational number in $[k_1, k_2]$ \\
        (range F \texttt{*string*}) & The range of $F$ is a string \\
        (range F \texttt{*boolean*}) & The range of $F$ are booleans \\
        \bottomrule
    \end{tabular}
    }
\end{table}

%% file: grammar_restriction_definitions.tex
\begin{table}[!htb]
    \centering
    \renewcommand{\arraystretch}{1.5}
    \rowcolors{2}{gray!10}{white}
    \caption{Definitions of datatype and feature restrictions.}
    \label{tab:restriction-definitions}
    {\footnotesize
    \begin{tabular}{ll}
        \toprule
        \hiderowcolors
        \textbf{Restriction} & \textbf{Definition} \\
        \showrowcolors
        \midrule
        $(\mathrm{>=}\ F\ \text{variable})$ & $\displaystyle \sup_{b \in {\Delta}_D} \left[F^\mathcal{I} (x, b) \otimes (b \geq \text{variable})\right]$ \\
        $(\mathrm{<=}\ F\ \text{variable})$ & $\displaystyle \sup_{b \in \Delta_D} \left[F^\mathcal{I} (x, b) \otimes (b \leq \text{variable})\right]$ \\
        $(=\ F\ \text{variable}) $ & $\displaystyle \sup_{b \in \Delta_D} \left[F^\mathcal{I} (x, b) \otimes (b = \text{variable})\right]$ \\
        $(\mathrm{>=}\ F\ \text{fuzzy\_number})$ & $\displaystyle \sup_{b^\prime, b \in \Delta_D} \left[F^\mathcal{I} (x, b) \otimes (b \geq b^\prime) \otimes {\text{fuzzy\_number}(b^\prime)}^\mathcal{I}\right]$ \\
        $(\mathrm{<=}\ F\ \text{fuzzy\_number})$ & $\displaystyle \sup_{b^\prime, b \in \Delta_D} \left[F^\mathcal{I} (x, b) \otimes (b \leq b^\prime) \otimes {\text{fuzzy\_number}(b^\prime)}^\mathcal{I}\right]$ \\
        $(=\ F\ \text{fuzzy\_number})$ & $\displaystyle \sup_{b^\prime, b \in \Delta_D} \left[F^\mathcal{I} (x, b) \otimes (b = b^\prime) \otimes {\text{fuzzy\_number}(b^\prime)}^\mathcal{I}\right]$ \\
        $(\mathrm{>=}\ F\ \mathrm{function}(F_1, \ldots, F_n))$ & $\displaystyle \sup_{b \in \Delta_D} \left[F^\mathcal{I} (x, b) \otimes (b \geq {\mathrm{function}(F_1, \ldots, F_n)}^{\mathcal{I}})\right]$ \\
        $(\mathrm{<=}\ F\ \mathrm{function}(F_1, \ldots, F_n))$ & $\displaystyle \sup_{b \in \Delta_D} \left[F^\mathcal{I} (x, b) \otimes (b \leq {\mathrm{function}(F_1, \ldots, F_n)}^{\mathcal{I}})\right]$ \\
        $(=\ F\ \mathrm{function}(F_1, \ldots, F_n))$ & $\displaystyle \sup_{b \in \Delta_D} \left[F^\mathcal{I} (x, b) \otimes (b = {\mathrm{function}(F_1, \ldots, F_n)}^{\mathcal{I}})\right]$ \\
        \bottomrule
    \end{tabular}
    }
\end{table}

%% file: grammar_milp_kinfty.tex
\begin{table}[!htb]
    \centering
    \renewcommand{\arraystretch}{1.5}
    \rowcolors{2}{gray!10}{white}
    \caption{Solver-dependent values of $k_{\infty}$.}
    \label{tab:milp-kinfty}
    \footnotesize
    \begin{tabular}{lc}
        \toprule
        \hiderowcolors
        \textbf{MILP Solver} & $k_{\infty}$ \\
        \showrowcolors
        \midrule
        Gurobi & $1000 \cdot (2^{31} - 1)$ \\
        PULP CBC & $2^{31} - 1$ \\
        MIP & $2^{31} - 1$ \\
        PULP GLPK & $2^{28} - 1$ \\
        PULP HiGHS & $2^{28} - 1$ \\
        PULP CPLEX & $2^{28} - 1$ \\
        \bottomrule
    \end{tabular}
\end{table}

%% file: grammar_concept_definitions.tex
\begingroup

\scriptsize
\renewcommand{\arraystretch}{1.5}
\rowcolors{2}{gray!10}{white}
\begin{longtable}{%
    @{}>{\raggedright\arraybackslash}p{0.3\textwidth}>{\raggedright\arraybackslash}p{0.2\textwidth}>{\raggedright\arraybackslash}p{0.4\textwidth}@{}%
}
    \caption{Definitions of fuzzyDL concept expressions.}
    \label{tab:concept-expression-definitions}\\
    
    \toprule
    \hiderowcolors
    \textbf{Expression} & {} & \textbf{EBNF production} \\
    \showrowcolors
    \midrule
    \endfirsthead

    \toprule
    \hiderowcolors
    \textbf{Expression} & {} & \textbf{EBNF production} \\
    \showrowcolors
    \midrule
    \endhead

    \bottomrule
    \hiderowcolors
    \multicolumn{3}{r}{\small Continue on next page}\\
    \showrowcolors
    \endfoot

    \bottomrule
    \endlastfoot

    \texttt{*top*} & top concept & $\top\ =\ 1$ \\
    \texttt{*bottom*} & bottom concept & $\perp\ =\ 0$ \\
    A & atomic concept $A$ & $A^\mathcal{I}(x)$ \\
    CFC & concrete fuzzy concept $CFC$ (e.g., crisp, left-shoulder, and so on) & $\mathrm{CFC}^\mathcal{I}(x)$ \\
    DR & datatype restriction $DR$ & $\mathrm{DR}^\mathcal{I}(x)$ \\
    (and $C_1$ $C_2$) & concept conjunction of $C_1$ and $C_2$ & $C_1^\mathcal{I}(x) \otimes C_2^\mathcal{I}(x)$ \\
    (g-and $C_1$ $C_2$) & G\"odel conjunction of $C_1$ and $C_2$ & $C_1^\mathcal{I}(x) \otimes_G C_2^\mathcal{I}(x)$ \\
    (l-and $C_1$ $C_2$) & {\L}ukasiewicz conjunction of $C_1$ and $C_2$ & $C_1^\mathcal{I}(x) \otimes_L C_2^\mathcal{I}(x)$ \\
    (or $C_1$ $C_2$) & concept disjunction of $C_1$ and $C_2$ & $C_1^\mathcal{I}(x) \oplus C_2^\mathcal{I}(x)$ \\
    (g-or $C_1$ $C_2$) & G\"odel disjunction of $C_1$ and $C_2$ & $C_1^\mathcal{I}(x) \oplus_G C_2^\mathcal{I}(x)$ \\
    (l-or $C_1$ $C_2$) & {\L}ukasiewicz disjunction of $C_1$ and $C_2$ & $C_1^\mathcal{I}(x) \oplus_L C_2^\mathcal{I}(x)$ \\
    (not $C$) & concept $C$ negation & $\ominus_L C^\mathcal{I}(x)$ \\
    (implies $C_1$ $C_2$) & concept implication between $C_1$ and $C_2$ & $C_1^\mathcal{I}(x) \Longrightarrow C_2^\mathcal{I}(x)$ \\
    (g-implies $C_1$ $C_2$) & G\"odel implication between $C_1$ and $C_2$ & $C_1^\mathcal{I}(x) \Longrightarrow_G C_2^\mathcal{I}(x)$ \\
    (l-implies $C_1$ $C_2$) & {\L}ukasiewicz implication between $C_1$ and $C_2$ & $C_1^\mathcal{I}(x) \Longrightarrow_L C_2^\mathcal{I}(x)$ \\
    (kd-implies $C_1$ $C_2$) & Kleene-Dienes implication between $C_1$ and $C_2$ & $C_1^\mathcal{I}(x) \Longrightarrow_\mathrm{KD} C_2^\mathcal{I}(x)$ \\
    (all $R$ $C$) & universal role $R$ restriction for concept $C$ & $\displaystyle\inf_{y \in \Delta^\mathcal{I}}\ \{ R^\mathcal{I}(x, y) \Longrightarrow C^\mathcal{I}(y) \} $ \\
    (some $R$ $C$) & existential role $R$ restriction for concept $C$ & $\displaystyle\sup_{y \in \Delta^\mathcal{I}}\ \{ R^\mathcal{I}(x, y) \otimes C^\mathcal{I}(y) \} $ \\
    (some $R$ $a$) & individual value restriction for role $R$ and individual $a$ & $R^\mathcal{I}(x, a) $ \\
    (ua $s$ $C$) & upper approximation for a fuzzy relation $s$ and individual $a$ & $\displaystyle\sup_{y \in \Delta^\mathcal{I}}\ \{ s^\mathcal{I}(x, y) \otimes C^\mathcal{I}(y) \} $ \\
    (lua $s$ $C$) & loose upper approximation for a fuzzy relation $s$ and individual $a$ & $\displaystyle\sup_{z \in X}\ \{ s^\mathcal{I}(x, z) \otimes \sup_{y \in \Delta^\mathcal{I}} s^\mathcal{I}(y, z) \otimes C^\mathcal{I}(x) \}$ \\
    (tua $s$ $C$) & tight upper approximation for a fuzzy relation $s$ and individual $a$ & $\displaystyle\inf_{z \in X}\ \{ s^\mathcal{I}(x, z) \Longrightarrow \sup_{y \in \Delta^\mathcal{I}} s^\mathcal{I}(y, z) \otimes C^\mathcal{I}(x) \}$ \\
    (la $s$ $C$) & lower approximation for a fuzzy relation $s$ and individual $a$ & $\displaystyle\inf_{y \in \Delta^\mathcal{I}}\ s^\mathcal{I}(x, y) \Longrightarrow C^\mathcal{I}(y) $ \\
    (lla $s$ $C$) & loose lower approximation for a fuzzy relation $s$ and individual $a$ & $\displaystyle\sup_{z \in X}\ \{ s^\mathcal{I}(x, z) \otimes \inf_{y \in \Delta^\mathcal{I}} s^\mathcal{I}(y, z) \otimes C^\mathcal{I}(x) \}$ \\
    (tla $s$ $C$) & tight lower approximation for a fuzzy relation $s$ and individual $a$ & $\displaystyle\inf_{z \in X}\ \{ s^\mathcal{I}(x, z) \Longrightarrow \inf_{y \in \Delta^\mathcal{I}} s^\mathcal{I}(y, z) \otimes C^\mathcal{I}(x) \}$ \\
    (self C) & local reflexivity concept & $\displaystyle C^\mathcal{I}(x)(x, x)$ \\
    (MOD C) & modifier MOD applied to concept $C$ & $\displaystyle {f_m}(C^\mathcal{I}(x))$, where $f_m$ is the modifier associated to MOD \\
    (FN) & fuzzy number FN & $\mathrm{FM}^\mathcal{I}(x)$ \\
    ([>= var ] C) & threshold concept &
    $\displaystyle \left\{
    \begin{aligned}
    C^\mathcal{I}(x), &\quad \text{if } C^\mathcal{I}(x) \geq \mathrm{var},\\
    0, &\quad \text{otherwise}
    \end{aligned}
    \right.$ \\
    
    ([<= var ] C) & threshold concept &
    $\displaystyle\left\{
    \begin{aligned}
    C^\mathcal{I}(x), &\quad \text{if } C^\mathcal{I}(x) \leq \mathrm{var},\\
    0, &\quad \text{otherwise}
    \end{aligned}
    \right.$ \\
    (n C) & weighted concept C with weight n & $n C^\mathcal{I}(x)$ \\
    (w-max $(v_1\ C_1) \ldots (v_k\ C_k)$) & weighted max of concepts & $\displaystyle\max_{i=1}^{k} \min \{v_i, x_i\}$ \\
    (w-min $(v_1\ C_1) \ldots (v_k\ C_k)$) & weighted min of concepts & $\displaystyle\min_{i=1}^{k} \max \{1 - v_i, x_i\}$ \\
    (w-sum $(n_1\ C_1) \ldots (n_k\ C_k)$) & weighted sum of concepts & $\displaystyle \sum_{i=1}^{k} n_i C_i^\mathcal{I}(x)$ \\
    (w-sum-zero $(n_1\ C_1) \ldots (n_k\ C_k)$) & weighted sum of concepts with zero handling & if $\displaystyle C_i^\mathcal{I}(x) = 0$ for some $i \in \{1, \ldots, k\}$, then $0$; otherwise $\displaystyle \sum_{i=1}^{k} n_i C_i^\mathcal{I}(x)$ \\
    (owa $(w_1, \ldots, w_n)$ $(C_1, \ldots, C_n)$) & OWA aggregation operator & $\displaystyle\sum_{i=1}^n w_i y_i $ \\
    (q-owa $Q\ (C_1, \ldots, C_n)$) & quantifier-guided OWA with name $Q$, where $Q$ is a right-shoulder or a linear function & $\displaystyle\sum_{i=1}^n w_i y_i $, where $w_i = Q(\frac{i}{n}) - Q(\frac{i - 1}{n})$ \\
    (choquet $(v_1, \ldots, v_n)$ $(C_1, \ldots, C_n)$) & Choquet integral & $\displaystyle y_1 v_1 + \sum_{i=2}^n (y_i - y_{i - 1}) v_i $ \\
    (sugeno $(v_1, \ldots, v_n)$ $(C_1, \ldots, C_n)$) & Sugeno integral & $\displaystyle \max_{i=1}^n \min \{y_i, mu_i\}$ \\
    (q-sugeno $(v_1, \ldots, v_n)$ $(C_1, \ldots, C_n)$) & Quasi-Sugeno integral & $\displaystyle \max_{i=1}^n y_i \otimes_L mu_i $ \\
    (sigma-count $R$ $C$ $[a_1\ \ldots\ a_k]$ $F_C$) & A Sigma-Count concept with role $R$ and associated to the concept $C$, the individuals $a_i$, and the fuzzy concrete concept $F_C$ &  {\tiny $\displaystyle \highi{F_C}\left(\sum_{a_i \in \{a_1\ \ldots\ a_k \}}\left(\highi{R}\left(x, \highi{a_i}\right) \otimes \highi{C}\left(\highi{a_i}\right)\right)\right) $} \\
\end{longtable}
\endgroup

%% file: grammar_axioms.tex
\begin{table}[!htb]
    \centering
    \renewcommand{\arraystretch}{2.5}
    \rowcolors{2}{gray!10}{white}
    \caption{Definitions of the \texttt{fuzzyDL} axioms. The degree value $d$ is optional; if omitted, it is assumed to be equal to $1.0$.}
    \label{tab:axiom-definitions}
    \scriptsize
    \begin{tabularx}{\textwidth}{lX}
        \toprule
        \hiderowcolors
        \textbf{Axiom} & \textbf{Definition} \\
        \showrowcolors
        \midrule
        (instance a C [d]) & $C^\mathcal{I}(a^\mathcal{I}) \geq d $ \\
        (related a b R [d]) & $R^\mathcal{I}(a^\mathcal{I}, b^\mathcal{I}) \geq d $ \\
        (implies $C_1$ $C_2$ [d]) & $\displaystyle\inf_{x \in \Delta^\mathcal{I}}\ C_1^\mathcal{I}(x) \Longrightarrow C_2^\mathcal{I}(x) \geq d$ \\
        (g-implies $C_1$ $C_2$ [d]) & $\displaystyle\inf_{x \in \Delta^\mathcal{I}}\ C_1^\mathcal{I}(x) \Longrightarrow_G C_2^\mathcal{I}(x) \geq d$ \\
        (kd-implies $C_1$ $C_2$ [d]) & $\displaystyle\inf_{x \in \Delta^\mathcal{I}}\ C_1^\mathcal{I}(x) \Longrightarrow_{\mathrm{KD}} C_2^\mathcal{I}(x) \geq d$ \\
        (l-implies $C_1$ $C_2$ [d]) & $\displaystyle\inf_{x \in \Delta^\mathcal{I}}\ C_1^\mathcal{I}(x) \Longrightarrow_{L} C_2^\mathcal{I}(x) \geq d $ \\
        (z-implies $C_1$ $C_2$ [d]) & $\displaystyle\inf_{x \in \Delta^\mathcal{I}}\ C_1^\mathcal{I}(x) \Longrightarrow_Z C_2^\mathcal{I}(x) \geq d$ \\
        (define-concept A C) & $\displaystyle\forall_{x \in \Delta^\mathcal{I}}\ A^\mathcal{I}(x) = C^\mathcal{I}(x) $ \\
        (define-primitive-concept A C) & $\displaystyle\inf_{x \in \Delta^\mathcal{I}}\ A^\mathcal{I}(x) \leq C^\mathcal{I}(x) $ \\
        (equivalent-concepts $C_1$ $C_2$) & $\displaystyle\forall_{x \in \Delta^\mathcal{I}}\ C_1^\mathcal{I}(x) = C_2^\mathcal{I}(x) $ \\
        (disjoint $C_1\ \ldots\ C_k$) & (implies (g-and $C_i$ $C_j$) \texttt{*bottom*}), i.e., $\displaystyle\forall_{i, j \in \{1, \ldots, k\}, i < j}\ (C_i^\mathcal{I}(x) \otimes_G C_j^\mathcal{I}(x)) \Longrightarrow \perp$ \\
        (disjoint-union $C_1\ \ldots\ C_k$) & $\displaystyle C_1 = \bigoplus_{i=2}^k C_i$ and $\displaystyle\forall_{i, j \in \{2, \ldots, k\}, i < j}\ (C_i^\mathcal{I}(x) \otimes_G C_j^\mathcal{I}(x)) \Longrightarrow \perp$ \\
        (range R $C$) & (implies \texttt{*top*} (all R C)), i.e., $\displaystyle \top \Longrightarrow \inf_{y \in \Delta^\mathcal{I}}\ \{ R^\mathcal{I}(x, y) \Longrightarrow C^\mathcal{I}(y) \}$ \\
        (domain R $C$) & (implies (some R \texttt{*top*}) C), i.e., $\displaystyle\sup_{y \in \Delta^\mathcal{I}}\ \{ R^\mathcal{I}(x, y) \otimes \top \} \Longrightarrow C^\mathcal{I}(x)$ \\
        (functional R) & $\displaystyle\min(R^\mathcal{I}(a, b), R^\mathcal{I}(a, c))>0 \Longrightarrow b = c$ \\
        (inverse-functional R) & $\displaystyle\min(R^\mathcal{I}(b, a), R^\mathcal{I}(c, a))>0 \Longrightarrow b = c$ \\
        (reflexive R) & $\displaystyle\forall_{a \in \Delta^\mathcal{I}}\ R^\mathcal{I}(a, a) = 1$ \\
        (symmetric R) & $\displaystyle\forall_{a, b \in \Delta^\mathcal{I}}\ R^\mathcal{I}(a, b) = R^\mathcal{I}(b, a)$ \\
        (transitive R) & $\displaystyle\forall_{a, b \in \Delta^\mathcal{I}}\ R^\mathcal{I}(a, b) \geq \sup_{c \in \Delta^\mathcal{I}} R^\mathcal{I}(a, c) \otimes R^\mathcal{I}(c, b)$ \\
        (implies-role $R_1$ $R_2$ [d]) & $\displaystyle\inf_{x, y \in \Delta^\mathcal{I}}\ R_1^\mathcal{I}(x, y) \Longrightarrow_L R_2^\mathcal{I}(x, y) \geq d$ \\
        (inverse $R_1$ $R_2$) & $\displaystyle R_1^\mathcal{I} \equiv {(R_2^\mathcal{I})}^{-1}$ \\
        \bottomrule
    \end{tabularx}
\end{table}

%% file: grammar_show_filters.tex
\begin{table}[!thb]
    \centering
    \renewcommand{\arraystretch}{1.5}
    \rowcolors{2}{gray!10}{white}
    \caption{Usage of \texttt{fuzzyDL} shows statements.}
    \label{tab:show-statements}
    \scriptsize
    \begin{tabularx}{\textwidth}{lX}
        \toprule
        \hiderowcolors
        \textbf{Statement} & \textbf{Meaning} \\
        \showrowcolors
        \midrule
        (show-concrete-fillers $F_1\ \ldots\ F_n$) & show value of the fillers of the features $F_i$, for $i = 1, \ldots, n$ \\
        (show-concrete-fillers-for \texttt{ind} $F_1\ \ldots\ F_n$) & show value of the fillers of $F_i$, with $i = 1, \ldots, n$, for the individual \texttt{ind} \\
        (show-concrete-instance-for \texttt{ind} $F\ C_1\ \ldots\ C_n$) & show degrees of being the $F$ filler of the individual \texttt{ind} instance of $C_i$, for $i = 1, \ldots, n$ \\
        (show-abstract-fillers $R_1\ \ldots\ R_n$) & show fillers of $R_i$ and membership to any concept \\
        (show-abstract-fillers-for \texttt{ind} $R_1\ \ldots\ R_n$) & show fillers of $R_i$, with $i = 1, \ldots, n$, for the individual \texttt{ind} and membership to any concept \\
        (show-concepts $a_1\ \ldots\ a_n$) & show membership of the individuals $a_i$ to any concept, for $i = 1, \ldots, n$ \\
        (show-instances $C_1\ \ldots\ C_n$) & show value of the instances of the concepts $C_i$, for $i = 1, \ldots, n$ \\
        (show-variables $x_1\ \ldots\ x_n$) & show value of the variables $x_i$, for $i = 1, \ldots, n$ \\
        \bottomrule
    \end{tabularx}
\end{table}

%% file: grammar_queries.tex
\begin{table}[!htb]
    \centering
    \renewcommand{\arraystretch}{2}
    \rowcolors{2}{gray!10}{white}
    \caption{Definitions of the \texttt{fuzzyDL} queries. In the queries \texttt{max-sat?} and \texttt{min-sat?}, the individual $a$ is optional; if omitted, the query returns the solution for a generic individual $a$.}
    \label{tab:query-definitions}
    
    \scriptsize  
    \begin{tabularx}{\textwidth}{lX}
        \toprule
        \hiderowcolors
        \textbf{Query} & \textbf{Definition} \\
        \showrowcolors
        \midrule
        (sat?) & Check if $\mathcal{K}$ is consistent \\
        (max-instance? a C) & $\displaystyle\sup\ \{n \mid \mathcal{K} \models \text{(instance a C n)}\}$ \\
        (min-instance? a C) & $\displaystyle\inf\ \{n \mid \mathcal{K} \models \text{(instance a C n)}\}$ \\
        (all-instances? C) & (min-instance? a C) for every individual of $\mathcal{K}$ \\
        (max-related? a b R) & $\displaystyle\sup\ \{n \mid \mathcal{K} \models \text{(related a b R n)} \}$ \\
        (min-related? a b R) & $\displaystyle\inf\ \{n \mid \mathcal{K} \models \text{(related a b R n)} \}$ \\
        (max-subs? C D) & $\displaystyle\sup\ \{n \mid \mathcal{K} \models \text{(implies D C n)} \}$ \\
        (min-subs? C D) & $\displaystyle\inf\ \{n \mid \mathcal{K} \models \text{(implies D C n)} \}$ \\
        (max-g-subs? C D) & $\displaystyle\sup\ \{n \mid \mathcal{K} \models \text{(g-implies D C n)} \}$ \\
        (min-g-subs? C D) & $\displaystyle\inf\ \{n \mid \mathcal{K} \models \text{(g-implies D C n)} \}$ \\
        (max-l-subs? C D) & $\displaystyle\sup\ \{n \mid \mathcal{K} \models \text{(l-implies D C n)} \}$ \\
        (min-l-subs? C D) & $\displaystyle\inf\ \{n \mid \mathcal{K} \models \text{(l-implies D C n)} \}$ \\
        (max-kd-subs? C D) & $\displaystyle\sup\ \{n \mid \mathcal{K} \models \text{(kd-implies D C n)} \}$ \\
        (min-kd-subs? C D) & $\displaystyle\inf\ \{n \mid \mathcal{K} \models \text{(kd-implies D C n)} \}$ \\
        (max-sat? C [a]) & $\displaystyle\sup_{\mathcal{I}}\ \sup_{a \in \Delta^\mathcal{I}}\ C^\mathcal{I}(a)$ \\
        (min-sat? C [a]) & $\displaystyle\inf_{\mathcal{I}}\ \inf_{a \in \Delta^\mathcal{I}}\ C^\mathcal{I}(a)$ \\
        (max-var? var) & $\displaystyle\sup\ \{\text{var} \mid \mathcal{K} \text{ is consistent}\}$ \\
        (min-var? var) & $\displaystyle\inf\ \{\text{var} \mid \mathcal{K} \text{ is consistent}\}$ \\
        (defuzzify-lom? C a F) & Defuzzify the value of F using the largest of the maxima \\
        (defuzzify-mom? C a F) & Defuzzify the value of F using the middle of the maxima \\
        (defuzzify-som? C a F) & Defuzzify the value of F using the smallest of the maxima \\
        (bnp? f) & Computes the Best Non-Fuzzy Performance (BNP) of a fuzzy number $f$ \\
        \bottomrule
    \end{tabularx}
\end{table}

%% file: appendix_fuzzydl_parser.tex
\begingroup

\renewcommand{\arraystretch}{1.5}
\rowcolors{2}{gray!10}{white}

\scriptsize
\begin{longtable}{%
    @{}>{\raggedright\arraybackslash}p{0.4\textwidth}>{\raggedright\arraybackslash}p{0.55\textwidth}@{}%
}
	\caption{\texttt{fuzzyDL} directives handled by the class \texttt{DLParser} and the corresponding parser actions or objects created in the current implementation. Here, \texttt{kb} denotes the instance of the class \texttt{KnowledgeBase} stored as an attribute of \texttt{DLParser}.}
	\label{tab:dl-parser-directives} \\ 
    
    \toprule
    \hiderowcolors
	\textbf{Directive} & \textbf{Actions / Objects Created} \\
    \showrowcolors
	\midrule
	\endfirsthead
	\toprule
    \hiderowcolors
	\textbf{Directive} & \textbf{Actions / Objects Created} \\
    \showrowcolors
	\midrule
	\endhead
    \bottomrule
    \hiderowcolors
    \multicolumn{2}{r}{\small Continue on next page}\\
    \showrowcolors
    \endfoot

	\bottomrule
    \endlastfoot
    
    Comments beginning with \texttt{\#} or \texttt{\%} &
    Recognised and suppressed. They do not create semantic objects and do not affect the Knowledge Base. \\

    \kw{(define-fuzzy-logic lukasiewicz)}, \linebreak \kw{(define-fuzzy-logic zadeh)}, \linebreak \kw{(define-fuzzy-logic classical)} &
    Parsed by \texttt{\_fuzzy\_logic\_parser}; the selected value is stored in the Knowledge Base through \texttt{kb.set\_logic(...)} and controls the semantics of conjunction, disjunction, negation, and implication. \\

    Atomic concept names: \kw{*top*}, \kw{*bottom*} &
    Mapped either to existing concepts stored in the Knowledge Base or to the distinguished truth concepts \texttt{TruthConcept.get\_top()} and \texttt{TruthConcept.get\_bottom()}. \\

    \kw{(define-truth-constant $c\ n$)} &
    Stored through \texttt{kb.set\_truth\_constants(c,n)}. The constant may later be reused inside degrees and constraints. \\

    \kw{(define-modifier $M$ linear-modifier($c$))} &
    Creates a \texttt{LinearModifier} object and adds it through \texttt{kb.add\_modifier}. \\

    \kw{(define-modifier $M$ triangular-modifier($a,b,c$))} &
    Creates a \texttt{TriangularModifier} object and adds it through \texttt{kb.add\_modifier}. \\

    \kw{(define-fuzzy-concept $C$ crisp($k_1,k_2,a,b$))} &
    Creates a \texttt{CrispConcreteConcept} and stores it in the Knowledge Base through \texttt{kb.add\_concept}. \\

    \kw{(define-fuzzy-concept $C$ left-shoulder($k_1,k_2,a,b$))} &
    Creates a \texttt{LeftConcreteConcept}, stores it in the Knowledge Base through \texttt{kb.add\_concept}, and sets \texttt{kb.concrete\_fuzzy\_concepts = True}. \\

    \kw{(define-fuzzy-concept $C$ right-shoulder($k_1,k_2,a,b$))} &
    Creates a \texttt{RightConcreteConcept}, stores it in the Knowledge Base through \texttt{kb.add\_concept}, and sets \texttt{kb.concrete\_fuzzy\_concepts = True}. \\

    \kw{(define-fuzzy-concept $C$ triangular($k_1,k_2,a,b,c$))} &
    Creates a \texttt{TriangularConcreteConcept}, stores it in the Knowledge Base through \texttt{kb.add\_concept}, and sets \texttt{kb.concrete\_fuzzy\_concepts = True}s. \\

    \kw{(define-fuzzy-concept $C$ trapezoidal($k_1,k_2,a,b,c,d$))} &
    Creates a \texttt{TrapezoidalConcreteConcept}, stores it in the Knowledge Base through \texttt{kb.add\_concept}, and sets \texttt{kb.concrete\_fuzzy\_concepts = True}. \\

    \kw{(define-fuzzy-concept $C$ linear($k_1,k_2,a,b$))} &
    Creates a \texttt{LinearConcreteConcept}, stores it in the Knowledge Base through \texttt{kb.add\_concept}, and sets \texttt{kb.concrete\_fuzzy\_concepts = True}. \\

    \kw{(define-fuzzy-concept $C$ modified($M,F$))} &
    Retrieves the modifier $M$ and the previously declared concrete concept $F$, creates a \texttt{ModifiedConcreteConcept}, stores it through \texttt{kb.add\_concept}, and sets \texttt{kb.concrete\_fuzzy\_concepts = True}. \\

    \kw{(define-fuzzy-number-range $k_1\ k_2$)} &
    Calls \texttt{TriangularFuzzyNumber.set\_range($k_1,k_2$)} and thereby fixes the admissible global domain for triangular fuzzy numbers. \\

    \kw{(define-fuzzy-number $N\ e$)},\linebreak with $e$ a constant, a triangular triple, \linebreak or an expression using \kw{f+}, \kw{f-}, \kw{f*}, \kw{f/} &
    Creates a \texttt{TriangularFuzzyNumber} named $N$ directly or computes one from previously defined fuzzy numbers, stores it through \texttt{kb.add\_fuzzy\_number}, and sets \texttt{kb.concrete\_fuzzy\_concepts = True}. \\

    \kw{(range $F$ *integer* $k_1\ k_2$)}, \linebreak \kw{(range $F$ *real* $k_1\ k_2$)}, \linebreak \kw{(range $F$ *boolean*)}, \linebreak \kw{(range $F$ *string*)} &
    Parsed as concrete-feature declarations named $F$. The parser invokes the appropriate feature constructor in the Knowledge Base, namely \texttt{define\_integer\_concrete\_feature}, \texttt{define\_real\_concrete\_feature}, \texttt{define\_boolean\_concrete\_feature}, or \texttt{define\_string\_concrete\_feature}. \\

    \kw{(crisp-concept $C_1 \ldots C_n$)} &
    For each listed symbol, it retrieves the corresponding concept and marks it as crisp through \texttt{kb.set\_crisp\_concept}. \\

    \kw{(crisp-role $R_1 \ldots R_n$)} &
    Marks each listed role as crisp through \texttt{kb.set\_crisp\_role}. \\

    \kw{(define-fuzzy-similarity $s$)} &
    Adds $s$ as a fuzzy similarity relation through \texttt{kb.add\_similarity\_relation}. \\

    \kw{(define-fuzzy-equivalence $s$)} &
    Adds $s$ as a fuzzy equivalence relation through \texttt{kb.add\_equivalence\_relation}. \\

    \kw{(not $C$)}, \linebreak \kw{(self $R$)} &
    Create, respectively, the negated concept $-C$ and a \texttt{SelfConcept}. In the \kw{self} case, the role is also added as abstract through \texttt{kb.abstract\_roles.add(R)}. \\

    \kw{(and $C_1\ \ldots\ C_n$)}, \linebreak \kw{(g-and $C_1\ \ldots\ C_n$)}, \linebreak \kw{(l-and $C_1\ \ldots\ C_n$)} &
    Create an \texttt{OperatorConcept} using the current logic or the explicitly requested G\"odel/{\L}ukasiewicz conjunction. The parser checks that all arguments are abstract concepts. \\

    \kw{(or $C_1\ \ldots\ C_n$)}, \linebreak \kw{(g-or $C_1\ \ldots\ C_n$)}, \linebreak \kw{(l-or $C_1\ \ldots\ C_n$)} &
    Create an \texttt{OperatorConcept} for disjunction, again using either the current logic or the explicitly requested fuzzy semantics. \\

    \kw{(implies $C_1\ C_2$)}, \linebreak \kw{(g-implies $C_1\ C_2$)}, \linebreak \kw{(l-implies $C_1\ C_2$)}, \linebreak \kw{(kd-implies $C_1\ C_2$)} &
    Create an \texttt{ImpliesConcept} with the appropriate implication operator. Under Zadeh logic, the generic implication is mapped to the Zadeh implication. \\

    \kw{(all $R\ C$)},\linebreak \kw{(some $R\ C$)} &
    Create the universal restriction \texttt{AllSomeConcept.all($R,C$)} and the existential restriction \texttt{AllSomeConcept.some($R,C$)} after validating the role $R$ by means of \texttt{kb.check\_role(R)}. This check ensures that $R$ is not simultaneously treated as both an abstract and a concrete role.\\

    \kw{(has-value $R\ a$)} &
    Create a \texttt{HasValueConcept.has\_value($R,a$)}, and marks the role $R$ as abstract. \\

    Approximation operators: \linebreak \kw{(tua $R\ C$)}, \linebreak \kw{(tla $R\ C$)}, \linebreak \kw{(ua $R\ C$)}, \linebreak \kw{(la $R\ C$)}, \linebreak \kw{(lua $R\ C$)}, \linebreak \kw{(lla $R\ C$)} &
    Create an \texttt{ApproximationConcept} with the requested approximation type and then convert it into an equivalent \texttt{AllSomeConcept} by means of the method \texttt{to\_all\_some\_concept()} of \texttt{ApproximationConcept}. The conversion is given by $\texttt{la} \mapsto \forall R.C$,
    $\texttt{tla} \mapsto \forall R.\forall R.C$,
    $\texttt{lla} \mapsto \exists R.\forall R.C$,
    $\texttt{ua} \mapsto \exists R.C$,
    $\texttt{tua} \mapsto \forall R.\exists R.C$, and
    $\texttt{lua} \mapsto \exists R.\exists R.C$.
    The role $R$ must already be declared as a fuzzy similarity relation; otherwise, the parser raises an error. \\

    \kw{([>= $v$] $C$)}, \linebreak \kw{([<= $v$] $C$)} &
    Create a \texttt{ThresholdConcept} when $v$ is numeric, or an \texttt{ExtThresholdConcept} when $v$ is parsed as a MILP variable. The parser checks that the concept $C$ is abstract. \\

    \kw{($M\ C$)}, with $M$ a previously declared modifier &
    Retrieves the modifier $M$ and the concept $C$, applies \texttt{mod.modify(concept)}, and returns an instance of \texttt{ModifiedConcept}. \\

    \kw{($n\ C$)} &
    Creates a \texttt{WeightedConcept} $C$ with weight $n$, used as an atomic component of weighted aggregates. \\

    \kw{(w-max $(v_1\ C_1)\ \ldots\ (v_k\ C_k)$)} &
    Creates a \texttt{WeightedMaxConcept}. The parser checks that, for $i = 1, \ldots, k $, $0 \leq v_i \leq 1$. \\

    \kw{(w-min $(v_1\ C_1)\ \ldots\ (v_k\ C_k)$)} &
    Creates a \texttt{WeightedMinConcept}. The parser checks that, for $i = 1, \ldots, k $, $0 \leq v_i \leq 1$. \\

    \kw{(w-sum $(n_1\ C_1)\ \ldots\ (n_k\ C_k)$)} &
    Creates a \texttt{WeightedSumConcept}. The parser checks that, for $i = 1, \ldots, k $, $0 \leq n_i \leq 1$, and $\sum_i n_i \leq 1$. \\

    \kw{(w-sum-zero $(n_1\ C_1)\ \ldots\ (n_k\ C_k)$)} &
    Creates a \texttt{WeightedSumZeroConcept}. The parser checks that, for $i = 1, \ldots, k $, $0 \leq n_i \leq 1$, and $\sum_i n_i \leq 1$. \\

    \kw{(q-owa $Q\ C_1\ \ldots\ C_n$)} &
    Retrieves $Q$ as a previously declared left-shoulder or right-shoulder concrete concept and creates a \texttt{QowaConcept} for concepts $C_1, \ldots, C_n$. \\

    \kw{(owa $(w_1\ \ldots\ w_n)$ $(C_1\ \ldots\ C_n)$)} &
    Creates an \texttt{OwaConcept}. The parser checks that $\sum_i w_i =  1$. \\

    \kw{(choquet $(v_1\ \ldots\ v_n)$ $(C_1\ \ldots\ C_n)$)} &
    Creates a \texttt{ChoquetIntegral}. The parser checks that $0\leq v_1 \leq \ldots \leq v_n = 1$. \\

    \kw{(sugeno $(v_1\ \ldots\ v_n)$ $(C_1\ \ldots\ C_n)$)} &
    Creates a \texttt{SugenoIntegral}. The parser checks that $0\leq v_1 \leq \ldots \leq v_n = 1$. \\

    \kw{(q-sugeno $(v_1\ \ldots\ v_n)$ $(C_1\ \ldots\ C_n)$)} &
    Creates a \texttt{QsugenoIntegral}. The parser checks that $0\leq v_1 \leq \ldots \leq v_n = 1$. \\

    \kw{(sigma-count $R\ C\ [a_1\ \ldots\ a_k]\ F_C$)} &
    Creates a \texttt{SigmaConcept} from the role $R$, base concept $C$, list of individuals $\{a_1, \ldots, a_k\}$, and fuzzy concrete concept $F_C$. The last argument must denote a previously declared left-shoulder, right-shoulder, or triangular concrete concept. \\

    Datatype restrictions with \kw{(<= $F\ x$)}, \kw{(>= $F\ x$)}, or \linebreak \kw{(= $F\ x$)} &
    Parsed by \texttt{\_parse\_datatype\_restriction}. In the current implementation, $F$ is treated as the concrete feature and $x$ is a bound, a fuzzy number, a variable, or a \texttt{FeatureFunction}. The parser creates and returns the corresponding datatype-restriction concept through \texttt{kb.add\_datatype\_restriction}. \\

    \kw{$(n [\ast] F)$}, \linebreak \kw{$(F_1 - F_2)$}, \linebreak \kw{$(F_1 + \ldots + F_n)$} &
    Converted into \texttt{FeatureFunction} objects by \texttt{\_parse\_restrictions}. These objects may then be embedded within datatype restrictions, and each $F_i$ may itself be a linear combination of \texttt{FeatureFunction} objects. For example, one may define \kw{(netPrice + ($0.2$ * netPrice) + ($0.3$ * netPrice))}, that is, a linear combination of the feature \texttt{netPrice}. \\
    
    \kw{(constraints ((expr >= $n$)) ...)},\linebreak \kw{(constraints ((expr <= $n$)) ...)},\linebreak \kw{(constraints ((expr = $n$)) ...)} &
    Expressions \texttt{expr} are translated into \texttt{Term}, \texttt{Expression}, and \texttt{Inequation} objects. The resulting linear constraints are injected into the MILP model through \texttt{kb.milp.add\_new\_constraint(...)}. \\

    \kw{(constraints ((binary $x$)) ...)} &
    Retrieves the MILP variable and sets its type to \texttt{VariableType.BINARY}. \\

    \kw{(constraints ((free $x$)) ...)} &
    Retrieves the MILP variable and sets its type to \texttt{VariableType.CONTINUOUS}. \\

    \kw{(show-concrete-fillers $F_1\ \ldots\ F_n$)} &
    Add each listed concrete role/feature in \texttt{kb.milp.show\_vars} so that its fillers are printed during inspection. \\

    \kw{(show-concrete-fillers-for ind  $F_1\ \ldots\ F_n$)} &
    Add the fillers of the listed concrete roles/features for the individual \texttt{ind}. \\

    \kw{(show-concrete-instance-for ind $F\ C_1\ \ldots\ C_n$)} &
    Checks that $F$ is concrete and that each $C_i$ is a concrete fuzzy concept or fuzzy number; then stores the corresponding inspection request in \texttt{show\_vars} for the individual \texttt{ind}. \\

    \kw{(show-abstract-fillers $R_1\ \ldots\ R_n$)} &
    Add each abstract role $R_i$ in \texttt{show\_vars} so that abstract fillers and their concept memberships can be displayed. \\

    \kw{(show-abstract-fillers-for ind $R_1\ \ldots\ R_n$)} &
    Add abstract fillers $R_i$ for the individual \texttt{ind}.  \\

    \kw{(show-concepts $a_1\ \ldots\ a_n$)} &
    Adds the listed individuals to the set of objects whose concept memberships must be displayed. \\

    \kw{(show-instances $C_1\ \ldots\ C_n$)} &
    Adds the listed concepts to the set of concepts whose instances must be displayed. \\

    \kw{(show-variables $x_1\ \ldots\ x_n$)} &
    Retrieves the corresponding MILP variables and marks them for display. \\

    \kw{(show-language)} &
    Sets \texttt{kb.show\_language = True}. This directive causes the parser to display the Description Logic expressivity of the Knowledge Base. The reported expressivity may range from $\mathcal{ALC}$ to $\mathcal{SHIF}(\mathrm{D})$. \\
    
    \kw{(instance $a\ C\ [d]$)} &
    Creates an assertion \texttt{kb.add\_assertion(a, C, d)} between individual $a$ and concept $C$ with optional degree $d$. When omitted, the degree defaults to \texttt{DegreeNumeric.get\_degree(1.0)}. \\

    \kw{(related $a\ b\ R\ [d]$)} &
    Creates an abstract-role assertion \texttt{kb.add\_relation(a,R,b,d)} between individuals $a$ and $b$, and role $R$ with optional degree $d$. Concrete roles are explicitly rejected in this context. When omitted, the degree defaults to \texttt{DegreeNumeric.get\_degree(1.0)}. \\

    \kw{(implies $C_1\ C_2\ [d]$)} &
    Adds a standard General Concept Inclusion (GCI) through \texttt{kb.implies($C_1,C_2,d$)} for concepts $C_1$ and $C_2$, with optional degree $d$. When omitted, the degree defaults to \texttt{DegreeNumeric.get\_degree(1.0)}. \\

    \kw{(g-implies $C_1\ C_2\ [d]$)} &
    Adds a G\"odel GCI through \texttt{kb.goedel\_implies($C_1,C_2,d$)}. \\

    \kw{(l-implies $C_1\ C_2\ [d]$)} &
    Adds a {\L}ukasiewicz GCI through \texttt{kb.lukasiewicz\_implies($C_1,C_2,d$)}. \\

    \kw{(kd-implies $C_1\ C_2\ [d]$)} &
    Adds a Kleene--Dienes GCI through \texttt{kb.kleene\_dienes\_implies($C_1,C_2,d$)}. \\

    \kw{(z-implies $C_1\ C_2$)} &
    Adds a Zadeh-style inclusion through \texttt{kb.zadeh\_implies($C_1,C_2$)}. \\

    \kw{(define-concept $A\ C$)} &
    Creates a concept $C$ with name $A$ through \texttt{kb.define\_concept}. \\

    \kw{(define-primitive-concept $A\ C$)} &
    Creates a primitive concept definition $C$ with name $A$ through \texttt{kb.define\_atomic\_concept}. \\

    \kw{(equivalent-concepts $C_1\ C_2$)} &
    Adds equivalence between concepts $C_1$ and $C_2$ through \texttt{kb.define\_equivalent\_concepts}. \\

    \kw{(disjoint $C_1\ \ldots\ C_k$)} &
    Add pairwise disjointness through \texttt{kb.add\_concepts\_disjoint}. \\

    \kw{(disjoint-union $C_1\ \ldots\ C_k$)} &
    Adds a disjoint union between $C_1$ and $C_i$, with $i=2, \ldots, n$, through \texttt{kb.add\_disjoint\_union\_concept}. \\

    \kw{(range $R\ C$)} as an axiom &
    Interpreted as a role-range axiom and stored through \texttt{kb.role\_range}, i.e., it stores the range of the role $R$ as the concept $C$. \\

    \kw{(domain $R\ C$)} &
    Interpreted as a role-domain axiom and stored through \texttt{kb.role\_domain}, i.e., it stores the domain of the role $R$ as the concept $C$. \\

    \kw{(functional $R$)}, \linebreak \kw{(inverse-functional $R$)}, \linebreak\kw{(reflexive $R$)}, \linebreak\kw{(symmetric $R$)}, \linebreak\kw{(transitive $R$)} &
    Set the corresponding role $R$ property in the Knowledge Base. Inverse functionality is rejected for concrete roles. \\

    \kw{(inverse $R_1\ R_2$)} &
    Add \texttt{$R_1$} and \texttt{$R_2$} as inverse roles through \texttt{kb.add\_inverse\_roles}. Concrete roles are rejected. \\

    \kw{(implies-role $R_1\ R_2\ [d]$)} &
    Adds a role-inclusion axiom through \texttt{kb.role\_implies}. If the degree is omitted, the implementation uses the crisp value $1.0$. \\

    \kw{(sat?)} &
    Creates a \texttt{KbSatisfiableQuery} and appends it to \texttt{DLParser.queries\_list}. \\

    \kw{(all-instances?~$C$)} &
    Creates an \texttt{AllInstancesQuery} and appends it to \texttt{DLParser.queries\_list}. \\

    \kw{(max-instance?~$a\ C$)}, \linebreak \kw{(min-instance?~$a\ C$)} &
    Create, respectively, a \texttt{MaxInstanceQuery} and a \texttt{MinInstanceQuery}, and append them to \texttt{DLParser.queries\_list}. \\

    \kw{(max-related?~$a\ b\ R$)}, \linebreak \kw{(min-related?~$a\ b\ R$)} &
    Create, respectively, a \texttt{MaxRelatedQuery} and a \texttt{MinRelatedQuery}, and append them to \texttt{DLParser.queries\_list}. The role $R$ is also recorded as abstract. \\

    \kw{(max-subs?~$C\ D$)}, \kw{(min-subs?~$C\ D$)},\linebreak
    \kw{(max-g-subs?~$C\ D$)}, \kw{(min-g-subs?~$C\ D$)}, \linebreak
    \kw{(max-l-subs?~$C\ D$)}, \kw{(min-l-subs?~$C\ D$)}, \linebreak
    \kw{(max-kd-subs?~$C\ D$)}, \kw{(min-kd-subs?~$C\ D$)} &
    Create a \texttt{MaxSubsumesQuery} or \texttt{MinSubsumesQuery}, and append them to \texttt{DLParser.queries\_list}. For the generic directives \kw{max-subs?} and \kw{min-subs?}, the implication operator is selected according to the current logic, namely {\L}ukasiewicz under \texttt{lukasiewicz} semantics and Zadeh otherwise. The remaining directives explicitly force the query to use, respectively, G\"odel, {\L}ukasiewicz, or Kleene--Dienes implication. \\

    \kw{(max-sat?~$C\ [a]$)}, \linebreak \kw{(min-sat?~$C\ [a]$)} &
    Create a \texttt{MaxSatisfiableQuery} or \texttt{MinSatisfiableQuery} for the concept $C$, optionally restricted to the supplied individual $a$. The resulting query is then appended to \texttt{DLParser.queries\_list}. \\

    \kw{(max-var?~expr)}, \linebreak \kw{(min-var?~expr)} &
    Create, respectively, a \texttt{MaxQuery} or \texttt{MinQuery} over the parsed MILP expression \texttt{expr}, and append them to \texttt{DLParser.queries\_list}. \\

    \kw{(defuzzify-lom?~$C\ a\ F$)}, \linebreak \kw{(defuzzify-mom?~$C\ a\ F$)}, \linebreak \kw{(defuzzify-som?~$C\ a\ F$)} &
    Create, respectively, a \texttt{LomDefuzzifyQuery}, \texttt{MomDefuzzifyQuery}, or \texttt{SomDefuzzifyQuery}, and append them to \texttt{DLParser.queries\_list}. The feature $F$ must already be declared as a concrete feature. \\

    \kw{(bnp?~$f$)} &
    Creates a \texttt{BnpQuery} from a fuzzy-number expression or a previously declared fuzzy number. The global fuzzy-number range must already be defined. The resulting query is then appended to \texttt{DLParser.queries\_list}. \\
\end{longtable}
\endgroup

%% file: main.arxiv.bbl
\begin{thebibliography}{10}

\bibitem{Baader_Calvanese_McGuinness_Nardi_Patel-Schneider_2007}
Franz Baader, Deborah~L. McGuinness, Daniele Nardi, and Peter~F.
  Patel-Schneider, editors.
\newblock {\em {The Description Logic Handbook: Theory, Implementation and
  Applications}}.
\newblock Cambridge University Press, 2 edition, 2007.

\bibitem{Cython}
Stefan Behnel, Robert Bradshaw, David Woods, Mat\'{u}\v{s} Valo, and Lisandro
  Dalc\'{\i}n.
\newblock {Cython -- C-Extensions for Python}.
\newblock Documentation: \url{https://cython.org/}, 2007.
\newblock Accessed: 2026-04-28.

\bibitem{Bezdek81}
James~C. Bezdek.
\newblock {\em {Pattern Recognition with Fuzzy Objective Function Algorithms}}.
\newblock Springer, 1981.

\bibitem{BOBILLO2012258}
Fernando Bobillo, Miguel Delgado, and Juan G\'{o}mez-Romero.
\newblock {DeLorean: A reasoner for fuzzy OWL 2}.
\newblock {\em Expert Systems with Applications}, 39(1):258--272, 2012.

\bibitem{fuzzydl_syntax}
Fernando Bobillo and Umberto Straccia.
\newblock {Syntax and Semantics of FuzzyDL}.
\newblock Webpage:
  \url{https://www.umbertostraccia.it/cs/software/fuzzyDL/fuzzyDL.html}.
\newblock Accessed: 2026-04-28.

\bibitem{bobillo2008fuzzydl}
Fernando Bobillo and Umberto Straccia.
\newblock {fuzzyDL: An Expressive Fuzzy Description Logic Reasoner}.
\newblock In {\em 2008 IEEE International Conference on Fuzzy Systems (IEEE
  World Congress on Computational Intelligence)}, pages 923--930, 2008.

\bibitem{bobillo2010representing}
Fernando Bobillo and Umberto Straccia.
\newblock {Representing fuzzy ontologies in OWL 2}.
\newblock In {\em International Conference on Fuzzy Systems}, pages 1--6, 2010.

\bibitem{BOBILLO20111073}
Fernando Bobillo and Umberto Straccia.
\newblock {Fuzzy ontology representation using {OWL} 2}.
\newblock {\em International Journal of Approximate Reasoning},
  52(7):1073--1094, 2011.

\bibitem{bobillo12}
Fernando Bobillo and Umberto Straccia.
\newblock {Generalized Fuzzy Rough Description Logics}.
\newblock {\em Information Sciences}, 189:43--62, 2012.

\bibitem{ASOC2013}
Fernando Bobillo and Umberto Straccia.
\newblock {Aggregation Operators for Fuzzy Ontologies}.
\newblock {\em Applied Soft Computing}, 13(9):3816--3830, 2013.

\bibitem{Bobillo16a}
Fernando Bobillo and Umberto Straccia.
\newblock Optimising fuzzy description logic reasoners with general concept
  inclusions absorption.
\newblock {\em Fuzzy Sets and Systems}, 292:98--129, 2016.

\bibitem{BobilloStraccia2016}
Fernando Bobillo and Umberto Straccia.
\newblock {The fuzzy ontology reasoner {fuzzyDL}}.
\newblock {\em Knowledge-Based Systems}, 95:12--34, 2016.

\bibitem{Bobillo18}
Fernando Bobillo and Umberto Straccia.
\newblock {Reasoning within Fuzzy OWL 2 EL Revisited}.
\newblock {\em Fuzzy Sets and Systems}, 351:1--40, 2018.

\bibitem{Penaloza2015}
Stefan Borgwardt, Felix Distel, and Rafael Pe{\~{n}}aloza.
\newblock {The limits of decidability in fuzzy description logics with general
  concept inclusions}.
\newblock {\em Artificial Intelligence}, 218:23--55, 2015.

\bibitem{Borgwardt2017}
Stefan Borgwardt and Rafael Pe{\~{n}}aloza.
\newblock {Fuzzy Description Logics -- A Survey}.
\newblock In {\em 11th International Conference on Scalable Uncertainty
  Management (SUM 2017)}, pages 31--45. Springer, 2017.

\bibitem{re2c}
Peter Bumbulis.
\newblock {re2c -- Regular Expressions to Code}.
\newblock Documentation: \url{https://re2c.org/}, 1993.
\newblock Accessed: 2026-04-28.

\bibitem{Calegari2008}
Silvia Calegari and Elie Sanchez.
\newblock {Object-fuzzy concept network: An enrichment of ontologies in
  semantic information retrieval}.
\newblock {\em Journal of the American Society for Information Science and
  Technology}, 59(13):2171--2185, 2008.

\bibitem{Cardillo24}
Franco~Alberto Cardillo, Franca Debole, and Umberto Straccia.
\newblock {PN-OWL: A Two Stage Algorithm to Learn Fuzzy Concept Inclusions from
  OWL Ontologies}.
\newblock {\em Fuzzy Sets and Systems}, 490(109048), 2024.

\bibitem{Cardillo22}
Franco~Alberto Cardillo and Umberto Straccia.
\newblock {{Fuzzy OWL-BOOST}: Learning Fuzzy Concept Inclusions via Real-Valued
  Boosting}.
\newblock {\em Fuzzy Sets and Systems}, 438:164--186, 2022.

\bibitem{CastronovoFilipponeGaliciLaRosaPavoneTabacchi:DSA-ISC:2025}
Lydia Castronovo, Giuseppe Filippone, Mario Galici, Gianmarco La~Rosa,
  Arianna~Maria Pavone, and Marco~Elio Tabacchi.
\newblock Sdf fuzzia: A fuzzy-based, ai support system for decision-making
  frameworks.
\newblock In Mario Pavone, Carlos~A. Coello~Coello, Raffaele Cerulli, Salvatore
  Greco, and El-Ghazali Talbi, editors, {\em Decision Sciences}, pages
  395--405, Cham, 2026. Springer Nature Switzerland.

\bibitem{CrossChen2018}
Valerie Cross and Shangye Chen.
\newblock {Fuzzy Ontologies: State of the Art Revisited}.
\newblock In {\em 37th Conference of the North American Fuzzy Information
  Processing Society (NAFIPS 2018)}, volume 831 of {\em Communications in
  Computer and Information Science}, pages 230--242. Springer, 2018.

\bibitem{CostaLaskeyLukasiewicz2008}
Paulo C.~G. da~Costa, Kathryn~B. Laskey, and Thomas Lukasiewicz.
\newblock {\em {Uncertainty Representation and Reasoning in the Semantic Web}},
  pages 315--340.
\newblock Information Science Reference, October 2008.

\bibitem{FilipponeLaRosaTabacchi:SUM:2024}
Giuseppe Filippone, Gianmarco La~Rosa, and Marco~Elio Tabacchi.
\newblock Sdf-fuzzia: A fuzzy-ontology based plug-in for the intelligent
  analysis of geo-thematic data.
\newblock In S{\'e}bastien Destercke, Maria~Vanina Martinez, and Giuseppe
  Sanfilippo, editors, {\em Scalable Uncertainty Management}, pages 163--169,
  Cham, 2025. Springer Nature Switzerland.

\bibitem{cbc_user_guide}
John Forrest and Robin Lougee-Heimer.
\newblock {CBC User Guide}.
\newblock COIN-OR Documentation:
  \url{https://www.coin-or.org/Cbc/cbcuserguide.html}.
\newblock Accessed: 2026-04-28.

\bibitem{glpk_gnu}
{Free Software Foundation}.
\newblock {GLPK (GNU Linear Programming Kit)}.
\newblock Website: \url{https://www.gnu.org/software/glpk/}.
\newblock Accessed: 2026-04-28.

\bibitem{Gao2005562}
Mingxia Gao and Chunnian Liu.
\newblock {Extending OWL by fuzzy description logic}.
\newblock In {\em 17th IEEE International Conference on Tools with Artificial
  Intelligence (ICTAI'05)}, volume 2005, page 562 – 567, 2005.

\bibitem{HermitReasonerWebsite}
Birte Glimm, Ian Horrocks, Boris Motik, Giorgos Stoilos, and Zhe Wang.
\newblock {HermiT OWL Reasoner}.
\newblock Website: \url{http://www.hermit-reasoner.com/}.
\newblock Accessed: 2026-04-28.

\bibitem{HermiT}
Birte Glimm, Ian Horrocks, Boris Motik, Giorgos Stoilos, and Zhe Wang.
\newblock {HermiT: An OWL 2 Reasoner}.
\newblock {\em J. Autom. Reason.}, 53(3):245–269, 10 2014.

\bibitem{grabisch2011aggregation}
Michel Grabisch, Jean-Luc Marichal, Radko Mesiar, and Endre Pap.
\newblock {Aggregation functions: Means}.
\newblock {\em Information Sciences}, 181(1):1--22, 2011.

\bibitem{Grau2008}
Bernardo~Cuenca Grau, Ian Horrocks, Boris Motik, Bijan Parsia, Peter
  Patel-Schneider, and Ulrike Sattler.
\newblock {OWL 2: The next step for OWL}.
\newblock {\em Journal of Web Semantics}, 6(4):309–322, 09 2008.

\bibitem{OWL2}
W3C OWL~Working Group.
\newblock {OWL 2 Web Ontology Language Document Overview (Second Edition)}.
\newblock Website:
  \url{https://www.w3.org/TR/2012/REC-owl2-overview-20121211/}, 12 2012.
\newblock Accessed: 2026-04-28.

\bibitem{OWL2Overview}
W3C OWL~Working Group.
\newblock {OWL 2 Web Ontology Language Document Overview (Second Edition)}.
\newblock Website: \url{https://www.w3.org/TR/owl2-overview/}, 12 2012.
\newblock Accessed: 2026-04-28.

\bibitem{GurobiOptimizer}
{Gurobi Optimization, LLC}.
\newblock {Gurobi Optimizer}.
\newblock Website: \url{https://www.gurobi.com/product}.
\newblock Accessed: 2026-04-28.

\bibitem{GurobiReferenceManual}
{Gurobi Optimization, LLC}.
\newblock {Gurobi Optimizer Reference Manual}.
\newblock Website:
  \url{https://docs.gurobi.com/projects/optimizer/en/current/index.html}.
\newblock Accessed: 2026-04-28.

\bibitem{gurobipy_12_0_0}
{Gurobi Optimization, LLC}.
\newblock {Python interface to Gurobi: gurobipy v12.0.0}.
\newblock PiPy: \url{https://pypi.org/project/gurobipy/}.
\newblock Accessed: 2026-04-28.

\bibitem{HajekP06}
Petr H{\'a}jek.
\newblock What does mathematical fuzzy logic offer to description logic?
\newblock In Elie Sanchez, editor, {\em Fuzzy Logic and the Semantic Web},
  Capturing Intelligence, chapter~5, pages 91--100. Elsevier, 2006.

\bibitem{highs_solver}
{Hall, J. and Galabova, I. and Gottwald, L. and Feldmeier, M.}
\newblock {HiGHS -- high performance software for linear optimization}.
\newblock Website: \url{https://highs.dev/}.
\newblock Accessed: 2026-04-28.

\bibitem{OWLAPI}
Matthew Horridge and Sean Bechhofer.
\newblock {The {OWL API}: {A Java API} for {OWL} Ontologies}.
\newblock {\em Semantic Web}, 2(1):11--21, 2011.

\bibitem{Datil}
Ignacio Huitzil and Fernando Bobillo.
\newblock {Fuzzy Ontology Datatype Learning using {Datil}}.
\newblock {\em Expert Systems With Applications}, 228:120299, 2023.

\bibitem{Huitzil20}
Ignacio Huitzil, Fernando Bobillo, Juan Gomez-Romero, and Umberto Straccia.
\newblock {Fudge: Fuzzy Ontology Building with Consensuated Fuzzy Datatypes}.
\newblock {\em Fuzzy Sets and Systems}, 401:91--112, 2020.

\bibitem{ibm_CPLEX_optimizer}
{IBM}.
\newblock {IBM ILOG CPLEX Optimizer}.
\newblock Website:
  \url{https://www.ibm.com/it-it/products/ilog-cplex-optimization-studio/cplex-optimizer}.
\newblock Accessed: 2026-04-28.

\bibitem{pulp_3_2_1}
{J.S. Roy}.
\newblock {PuLP is an LP modeler written in python. PuLP can generate MPS or LP
  files and call GLPK, COIN CLP/CBC, CPLEX, and GUROBI to solve linear
  problems: PuLP v3.2.1}.
\newblock PiPy: \url{https://pypi.org/project/PuLP/}.
\newblock Accessed: 2026-04-28.

\bibitem{Klir95}
George~J. Klir and Bo~Yuan.
\newblock {\em {Fuzzy sets and fuzzy logic: theory and applications}}.
\newblock Prentice-Hall, Inc., 1995.

\bibitem{Knappe2007}
Rasmus Knappe, Henrik Bulskov, and Troels Andreasen.
\newblock {Perspectives on ontology-based querying}.
\newblock {\em International Journal of Intelligent Systems}, 22(7):739--761,
  2007.

\bibitem{uncertainty}
Kenneth Laskey and Kathryn Laskey.
\newblock {Uncertainty Reasoning for the World Wide Web: Report on the URW3-XG
  Incubator Group}.
\newblock Website:
  \url{https://www.w3.org/2005/Incubator/urw3/XGR-urw3-20080331}, 2008.
\newblock Accessed: 2026-04-28.

\bibitem{LUKASIEWICZ2008291}
Thomas Lukasiewicz and Umberto Straccia.
\newblock {Managing uncertainty and vagueness in description logics for the
  Semantic Web}.
\newblock {\em Journal of Web Semantics}, 6(4):291--308, 2008.
\newblock Semantic Web Challenge 2006/2007.

\bibitem{pyparsing}
Paul McGuire.
\newblock {PyParsing -- A Python Parsing Module}.
\newblock PiPy: \url{https://pypi.org/project/pyparsing/}.
\newblock Accessed: 2026-04-28.

\bibitem{owl2_profiles_w3c}
Boris Motik, Bernardo Cuenca~Grau, Ian Horrocks, Zhe Wu, Achille Fokoue, and
  Carsten Lutz.
\newblock {OWL 2 Web Ontology Language Profiles (Second Edition)}.
\newblock Website: \url{https://www.w3.org/TR/owl2-profiles/}.
\newblock Accessed: 2026-04-28.

\bibitem{ProtegePaper}
Mark~A. Musen.
\newblock {The prot{\'{e}}g{\'{e}} project: a look back and a look forward}.
\newblock {\em {AI} Matters}, 1(4):4--12, 2015.

\bibitem{flex}
Vern Paxson.
\newblock {Fast Lexical Analyzer Generator}.
\newblock Documentation: \url{https://github.com/westes/flex}, 1987.
\newblock Accessed: 2026-04-28.

\bibitem{python_private_name_mangling}
{Python Software Foundation}.
\newblock {6. Expressions --- Private name mangling}.
\newblock Python Documentation:
  \url{https://docs.python.org/3/reference/expressions.html\#private-name-mangling}.
\newblock Accessed: 2026-04-28.

\bibitem{salkin1989foundations}
H.M. Salkin, K.~Mathur, and R.~Haas.
\newblock {\em {Foundations of Integer Programming}}.
\newblock North-Holland, 1989.

\bibitem{sanchez2006fuzzy}
Elie Sanchez, editor.
\newblock {\em {Fuzzy logic and the semantic web}}, volume~1 of {\em Capturing
  Intelligence}.
\newblock Elsevier, 2006.

\bibitem{fuzzydlowl2_github}
{SDF Team}.
\newblock {Fuzzy DL OWL 2}.
\newblock GitHub Repository:
  \url{https://github.com/SDF-Unipa/fuzzy\_dl\_owl2}, 2025.
\newblock Accessed: 2026-04-28.

\bibitem{fuzzydlowl2_docs}
{SDF Team}.
\newblock {Fuzzy DL OWL 2}.
\newblock ReadTheDocs: \url{https://fuzzy-dl-owl2.readthedocs.io/en/latest},
  2025.
\newblock Accessed: 2026-04-28.

\bibitem{Simou2008}
N.~Simou, Th. Athanasiadis, G.~Stoilos, and S.~Kollias.
\newblock {Image indexing and retrieval using expressive fuzzy description
  logics}.
\newblock {\em Signal, Image and Video Processing}, 2(4):321--335, 2008.

\bibitem{PelletReasonerGithub}
Evren Sirin, Bijan Parsia, Bernardo~Cuenca Grau, Aditya Kalyanpur, and Yarden
  Katz.
\newblock {Pellet: An Open Source OWL DL reasoner for Java}.
\newblock Github Repository: \url{https://github.com/stardog-union/pellet}.
\newblock Accessed: 2026-04-28.

\bibitem{SIRIN200751}
Evren Sirin, Bijan Parsia, Bernardo~Cuenca Grau, Aditya Kalyanpur, and Yarden
  Katz.
\newblock {Pellet: A practical OWL-DL reasoner}.
\newblock {\em Journal of Web Semantics}, 5(2):51--53, 2007.
\newblock Software Engineering and the Semantic Web.

\bibitem{Stoilos2006}
G.~Stoilos, N.~Simou, G.~Stamou, and S.~Kollias.
\newblock {Uncertainty and the Semantic Web}.
\newblock {\em IEEE Intelligent Systems}, 21(5):84--87, 2006.

\bibitem{Stoilos2010656}
G.~Stoilos, G.~Stamou, and J.Z. Pan.
\newblock {Fuzzy extensions of OWL: Logical properties and reduction to Fuzzy
  Description Logics}.
\newblock {\em International Journal of Approximate Reasoning}, 51(6):656 –
  679, 2010.

\bibitem{stoilos2007extending}
Giorgos Stoilos and Giorgos~B Stamou.
\newblock {Extending Fuzzy Description Logics for the Semantic Web.}
\newblock In {\em 3rd International Workshop on OWL: Experience and Directions
  (OWLED 2007)}, volume 258 of {\em CEUR Workshop Proceedings}, 2007.

\bibitem{fuzzydl_web}
Umberto Straccia.
\newblock {The fuzzyDL System}.
\newblock Webpage:
  \url{https://www.umbertostraccia.it/cs/software/fuzzyDL/fuzzyDL.html}.
\newblock Accessed: 2026-04-28.

\bibitem{straccia2014foundations}
Umberto Straccia.
\newblock {\em {Foundations of Fuzzy Logic and Semantic Web Languages}}.
\newblock Chapman and Hall/CRC, New York, 1 edition, 2014.

\bibitem{protege}
Prot\'{e}g\'{e} Team.
\newblock {Prot\'{e}g\'{e}}.
\newblock Website: \url{https://protege.stanford.edu/}.
\newblock Accessed: 2026-04-28.

\bibitem{pyowl2}
SDF Team.
\newblock {PyOWL2}.
\newblock PiPy: \url{https://pypi.org/project/pyowl2/}.
\newblock Accessed: 2026-04-28.

\bibitem{pyowl2_github}
SDF Team.
\newblock {PyOWL2}.
\newblock GitHub Repository: \url{https://github.com/SDF-Unipa/pyowl2/}.
\newblock Accessed: 2026-04-28.

\bibitem{pyowl2_docs}
SDF Team.
\newblock {PyOWL2}.
\newblock ReadTheDocs: \url{https://pyowl2.readthedocs.io/en/latest/}.
\newblock Accessed: 2026-04-28.

\bibitem{Foga2006}
Q.T. Tho, S.C. Hui, A.C.M. Fong, and Tru~Hoang Cao.
\newblock {Automatic fuzzy ontology generation for semantic Web}.
\newblock {\em IEEE Transactions on Knowledge and Data Engineering},
  18(6):842--856, 2006.

\bibitem{mip_1_16rc0}
Tulio A.~M. Toffolo.
\newblock {Python tools for Modeling and Solving Mixed-Integer Linear Programs
  (MIPs): mip v1.16rc0}.
\newblock PiPy: \url{https://pypi.org/project/mip/}.
\newblock Accessed: 2026-04-28.

\bibitem{Yen1991}
John Yen.
\newblock {Generalizing term subsumption languages to fuzzy logic}.
\newblock In {\em 12th International Joint Conference on Artificial
  Intelligence (IJCAI 1991) - Volume 1}, pages 472–--477. Morgan Kaufmann
  Publishers Inc., 1991.

\bibitem{zadeh1965fuzzy}
Lotfi~Asker Zadeh.
\newblock {Fuzzy sets}.
\newblock {\em Information and control}, 8(3):338--353, 1965.

\end{thebibliography}
